\documentclass[11pt]{article}

\usepackage{arxiv}

\usepackage{amsmath,amsfonts,latexsym,fullpage,graphicx,vmargin, mathrsfs}
\usepackage{pstricks}
\usepackage{comment}
\usepackage[inline]{enumitem}
\usepackage{cleveref}

\usepackage[ruled]{algorithm2e}
\usepackage{algpseudocode}

\RequirePackage[OT1]{fontenc}

\usepackage{amsmath}
\usepackage{amsfonts}
\usepackage{caption}
\usepackage{subcaption}
\usepackage{multirow}
\usepackage{tikz}
\usetikzlibrary{graphs}
\usetikzlibrary{cd}
\usepackage{bm}
\usepackage{url}
\usepackage{placeins}
\usepackage{float}

\RequirePackage{natbib}

\newtheorem{teo}{Theorem}
\newtheorem{prop}{Proposition}
\newtheorem{defi}{Definition}
\newtheorem{rmk}{Remark}
\newtheorem{lem}{Lemma}
\newtheorem{cor}{Corollary}

\newtheorem{assump}{Assumption}

\DeclareMathOperator*{\Crit}{Crit}
\DeclareMathOperator*{\CA}{CA}
\DeclareMathOperator*{\LCA}{LCA}

\DeclareMathOperator*{\supp}{supp}
\DeclareMathOperator*{\Mapp}{Mapp}

\DeclareMathOperator*{\Cov}{Cov}

\DeclareMathOperator*{\im}{Im}

\DeclareMathOperator*{\Top}{Top}
\DeclareMathOperator*{\FTop}{FTop}
\DeclareMathOperator*{\Sets}{Sets}
\DeclareMathOperator*{\FSet}{FSet}
\DeclareMathOperator*{\Vect}{Vec}

\usepackage{scalerel}
\DeclareMathOperator*{\Bigcdot}{\scalerel*{\cdot}{\bigodot}}

\newcommand{\X}{X_{\Bigcdot}}
\newcommand{\Y}{Y_{\Bigcdot}}

\newcommand{\T}{\pi_0(\X)}
\newcommand{\G}{\pi_0(\Y)}

\newcommand{\R}{\mathbb{R}}
\newcommand{\N}{\mathbb{N}}

\newcommand{\virgolette}[1]{``#1''}

\firstpageno{1}

\begin{document}

\title{A Finitely Stable Edit Distance 
        for Functions Defined on Merge Trees}

\author{Matteo Pegoraro\thanks{Department of Mathematics, KTH, Stockholm}}

%\author{\name Matteo Pegoraro \email matteo.pegoraro@polimi.it \\
%       \addr MOX --- Department of Mathematics \\
%       Politecnico di Milano \\
%       Milano, Italy
%       \AND
%       \name Mario Beraha \email mario.beraha@polimi.it \\
%       \addr Department of Computer Science \\
%       Universit\`a degli Studi di Bologna \\
%       Bologna, Italy}
%
\maketitle

\begin{abstract}
In this work we define a metric structure to compare functions defined on different merge trees. 
The metric introduced possesses some stability properties, which we illustrate within a standard topological data analysis (TDA) framework, and can be computed with a dynamical binary linear programming approach.
We showcase the effectiveness of the whole framework with simulated data sets. Using functions defined on merge trees proves to be very effective in situations where other topological data analysis tools, like persistence diagrams, cannot be used meaningfully.  
\end{abstract}


%\MPnote{mettere questione curvatura di alexandrov un po' meglio!}

\section{Introduction}

Topological Data Analysis (TDA) is the name given to an ensemble of 
techniques which are mainly focused on retrieving topological
information from different kinds of data \cite{insights_complex_data}. Consider for instance the case of point clouds: the (discrete) topology of a point cloud itself is quite poor and it would be much 
more interesting if, using the point cloud, one could gather information 
about the topological space data was sampled from.
Since, in practice, this is often not possible, one can still try to capture the \virgolette{shape} of the point cloud. The idea of \emph{persistent homology} (PH) \cite{PH_survey} is 
an attempt to do so: using the initial point cloud, a nested 
sequence of topological spaces is built, which are heavily dependent on the initial point cloud, and PH
tracks along this sequence the persistence of the different topological 
features which appear and disappear. As the name \emph{persistent homology} suggests, the topological features are understood in terms of generators of 
the homology groups \cite{hatcher} taken along the sequence of 
spaces.
One of the foundational results in TDA is that this information can be represented by a set of points on the plane \cite{PD_1, PD_2}, with a point of coordinates $(x,y)$ representing a topological feature being born at time $x$ along the sequence, and disappearing at time $y$. Such representation is called \emph{persistence diagram} (PD).
%Sometimes a representation called Persistence Barcode  is preferred to Persistence Diagrams, but they convey exactly the same information.
Persistence diagrams can be given a metric structure through the \emph{Bottleneck} and \emph{Wasserstein} metrics, which, despite having 
good properties in terms of continuity with respect to perturbation of the original data \cite{cohen_PD, lip}, provide badly behaved metric spaces - with non-unique geodesics arising in many situations.  
Various attempts to define tools to work in such spaces 
have been made \cite{confidence_PD, cuturi_OPT, prob_meas_PD, means_PD}, but this still proves to be a hard problem.
In order to obtain spaces with better properties - e.g. with unique means - and/or information which is vectorized, a number of topological summaries alternative to PDs have been proposed, such as: persistence landscapes \cite{landscapes}, persistence images \cite{pers_img} and persistence silhouettes \cite{silhouettes}.
For a review on TDA vectorization techniques see \cite{ali2023survey}.

All the aforementioned machinery has been successfully applied to a
great number of problems in a very diverse set of scientific fields:
complex shape analysis \cite{meas_shape},
sensor network coverage 
\cite{sensor_network}, 
protein structures 
\cite{protein_2, protein_1},
%, protein_3,protein_4},
DNA and RNA structures 
\cite{DNA,RNA},
robotics 
\cite{robotics_1, robotics_2},
%robotics_3},
signal analysis 
%\cite{signal_1,signal_2,signal_3},
%more general time series 
and dynamical systems 
\cite{time_series_3, time_series_2, time_series_1}, 
%cancer treatment 
%\cite{cancer_1, cancer_2, 
%cancer_3, 
%cancer_4},
materials science 
\cite{materials_2, materials_1}, 
%materials_3, materials_4, materials_5},
%finance 
%\cite{finance_1, finance_2, finance_3},
neuroscience 
\cite{neuro_1,neuro_2}, 
%neuro_3, 
%neuro_4, 
%neuro_5, neuro_6, neuro_7, neuro_8, neuro_9}, 
%more general 
network analysis 
\cite{net_2, net_1},
%, net_3, net_4}, 
and even deep learning theory \cite{deep_1, deep_2}.

\subsection*{Related Works}

Close to the definition of persistent homology for $0$ dimensional homology groups, lie the ideas of \emph{merge trees} of functions, \emph{phylogenetic trees} and \emph{hierarchical clustering dendrograms}.
Merge trees of functions \cite{merge_parall_1} describe the path-connected components of the sublevel sets of a real valued function and are obtained as a particular case of \emph{Reeb graphs} \cite{reeb, reeb_2}, representing the evolution of the level sets of a bounded Morse function \cite{morse} defined on a path-connected domain.  Phylogenetic trees  \cite{phylo,garba2021information} and clustering dendrograms \cite{dendro_1, dendro_2} are very similar objects which describe the evolution of a set of labels under some similarity measure or agglomerative criterion and can be framed as merge trees of some filtration on the set of labels. 
Informally speaking, from persistence diagrams, under general conditions, one can not reconstruct the merging pattern of path-connected components/labels/clusters along a filtration, information which, instead, is conveyed by merge trees, phylogenetic trees, and clustering dendrograms \cite{ curry2024trees, kanari2020trees}. 

Given the widespread use of all the aforementioned trees, in the last years a lot of research sparkled on such topics, with particular focus on defining metric structures, with the aim of employing populations of Reeb graphs or merge trees for data analysis. 
Different but related metrics have been proposed to compare Reeb graphs \cite{bauer2020reeb,   de2016categorified, di2016edit}, and merge trees \cite{merge, cardona2021universal,  cavinato2022imaging, merge_intrins, merge_interl,    pegoraro2021interleaving, pegoraro2024finitely,   merge_wass, merge_farlocca,  merge_frechet, interl_approx, merge_farlocca_2}. There is also some work investigating structures lying in between merge trees and persistence diagrams \cite{mergegrams}.

In this work, we are interested in investigating these topics from a different perspective: we want to exploit merge trees to build topological summaries which contain more information about the starting filtration w.r.t. persistence diagrams and merge trees themselves. To do so, we exploit that merge trees parametrize path-connected components, and evaluate some function on each path-connected component in the filtration, obtaining scalar fields defined on merge trees which describe measure-related information or additional homological features. 
A related idea can be found \cite{curry2021decorated}, in which the authors propose the idea of \emph{decorated merge trees}: 
they decorate each merge tree with a barcode, decomposing higher dimensional homological phenomena across path-connected components. A similar idea is also developed for Reeb graphs  in \cite{curry2023topologically}. See \Cref{sec:decorated} for more details.
We also highlight that \cite{biswas2022window} considers the idea of defining functions on graph-shaped objects, with the aim of studying their sublevel set filtration.

\subsection*{Main Contributions}

In this work, we aim to define a metric to compare functions defined on different merge trees. To this end, we build on the edit distance for weighted trees introduced in \cite{pegoraro2023edit}, using as edge weights the restriction of the functions themselves. We develop a formal framework to justify this construction, culminating in \Cref{teo:fun_stab} and \Cref{prop:stab_betti}, which establish the stability properties of the resulting metric in a setting that is standard within topological data analysis.

\begin{figure}
    	\centering
	    \includegraphics[width = 0.49\textwidth]{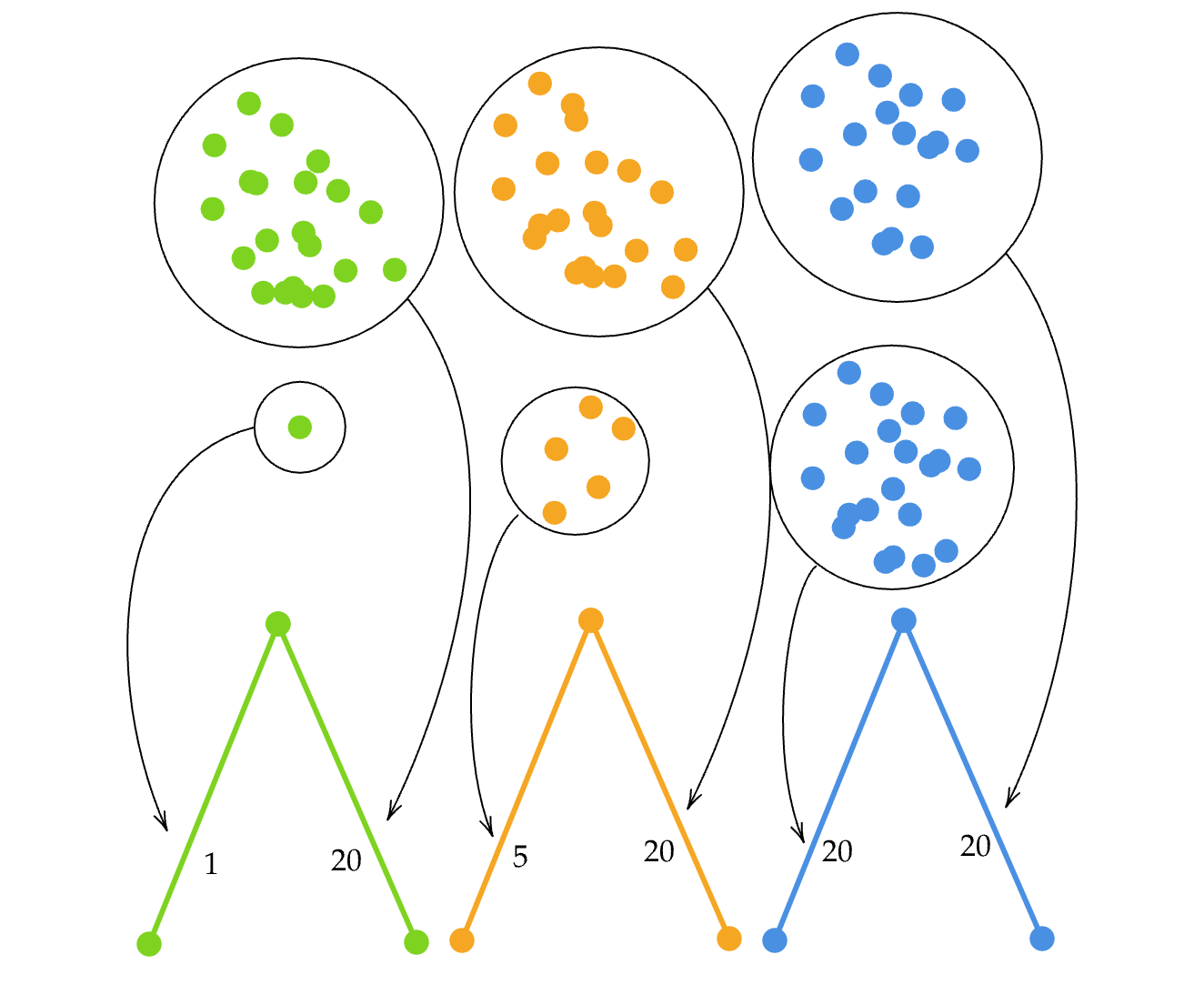}
    	\caption{Qualitative example to motivate the idea of considering functions defined on merge trees. The three point clouds, identified by the different colors, present very similar clustering structures from the point of view of TDA: two main clusters. As a result,  the merge trees (and so the persistence diagrams) associated to the Vietoris-Rips filtration, whose main branches are qualitatively depicted below each point cloud, are very similar. However, if we attach to each edge the number or the percentage of points in the associated path-connected component, we would be able to recognize three very distinct scenarios in terms of data-generating distributions.}
    	\label{fig:cluster_ex}
\end{figure}

The motivation behind this work stems from the observation that defining functions on merge trees can effectively extract additional information from data, offering valuable insights in a variety of data analysis scenarios, as already shown by \cite{curry2021decorated}.
Differently from \cite{curry2021decorated}, we take a broader approach to the topic, which we now motivate with the qualitative example depicted in \Cref{fig:cluster_ex}, which will be presented in a formal way in later sections. 

\Cref{fig:cluster_ex} displays three point clouds, each shown in a different color, originating from clearly distinct generative processes. In particular, the distribution of points between the two main clusters in each point cloud varies significantly. For instance, the smaller cluster in the green point cloud, could be considered an outlier.
However, both the merge trees and the persistence diagrams of the Vietoris-Rips filtration of the point clouds would be very close to each other, as, from a topological point of view, there are two main path-connected components, whose critical values are very similar across the filtration. Instead, if we  attach to each edge the number or the percentage of points which are contained in the associated path-connected component, we are able to clearly separate these scenarios. 

This can be seen as a naive formulation of the core idea we develop: more precisely, our approach leverages the parametrization of the path-connected components of a filtration provided by merge trees to extract, at each point in their geometric realization, information about the associated path-connected component.

\subsection*{Outline}

The paper is organized as follows. 
In \Cref{sec:tree_like} and \Cref{sec:intro_MT} we introduce most of the definitions needed for our dissertation, starting from most recent TDA literature, and we tackle the problem of representing with a discrete summary - a merge tree - the merging pattern of the path-connected components of a filtration of topological spaces. 
Once merge trees are introduced, we use  \Cref{sec:back_to_vec} to formally introduce spaces of
functions on a merge tree.
In \Cref{sec:edit_functions} we face the problem of finding a suitable metric structure to compare functions defined on different trees. In \Cref{sec:info_examples} we present some examples to showcase situations in which functions defined on merge trees can be useful. \Cref{sec:stability} contains the main theoretic results of the paper, which validate the whole framework. \Cref{sec:simulated_data} presents some simulated scenarios to test some of the functions defined in \Cref{sec:info_examples}.
We end up with some conclusions in  \Cref{sec:conclusions}.

The Appendix contains most of the proofs of some results, along with simulation studies and useful material which can help the reader in navigating through the content of the present work with multiple examples and additional details.
The outline of such contents is presented at the beginning of the Appendix and coherently referenced through the paper.

\section{Abstract Merge Trees}\label{sec:tree_like}

%In TDA the main sources of information are sequences of homology groups with field coefficients: using different pipelines a single datum is turned into a filtration of topological spaces $\{X_t\}_{t\in\mathbb{R}}$, which, in turn, induces - via some homology functor with coefficients in the field $\mathbb{K}$ - a family of vector spaces with linear maps which are usually all isomorphisms but for a finite set of points in the sequence. Such objects are called (one-dimensional) \emph{persistence modules} \cite{chazal2008persistencemodules}.
%Any persistence module is then turned into a topological summary, for instance a persistence diagram, which completely classifies such objects up to isomorphisms. That is, if for two persistence modules (satisfying certain finiteness conditions) there exists a family of linear isomorphisms giving a natural transformation between the two functors, then they are represented by the same persistence diagram. And viceversa.

In the first part of this work we describe a categorical approach to merge trees.

\subsection{Preliminary Definitions}
\label{sec:defi_intro}
\begin{figure}
    \begin{subfigure}[c]{0.49\textwidth}
    	\centering
    	\includegraphics[width = \textwidth]{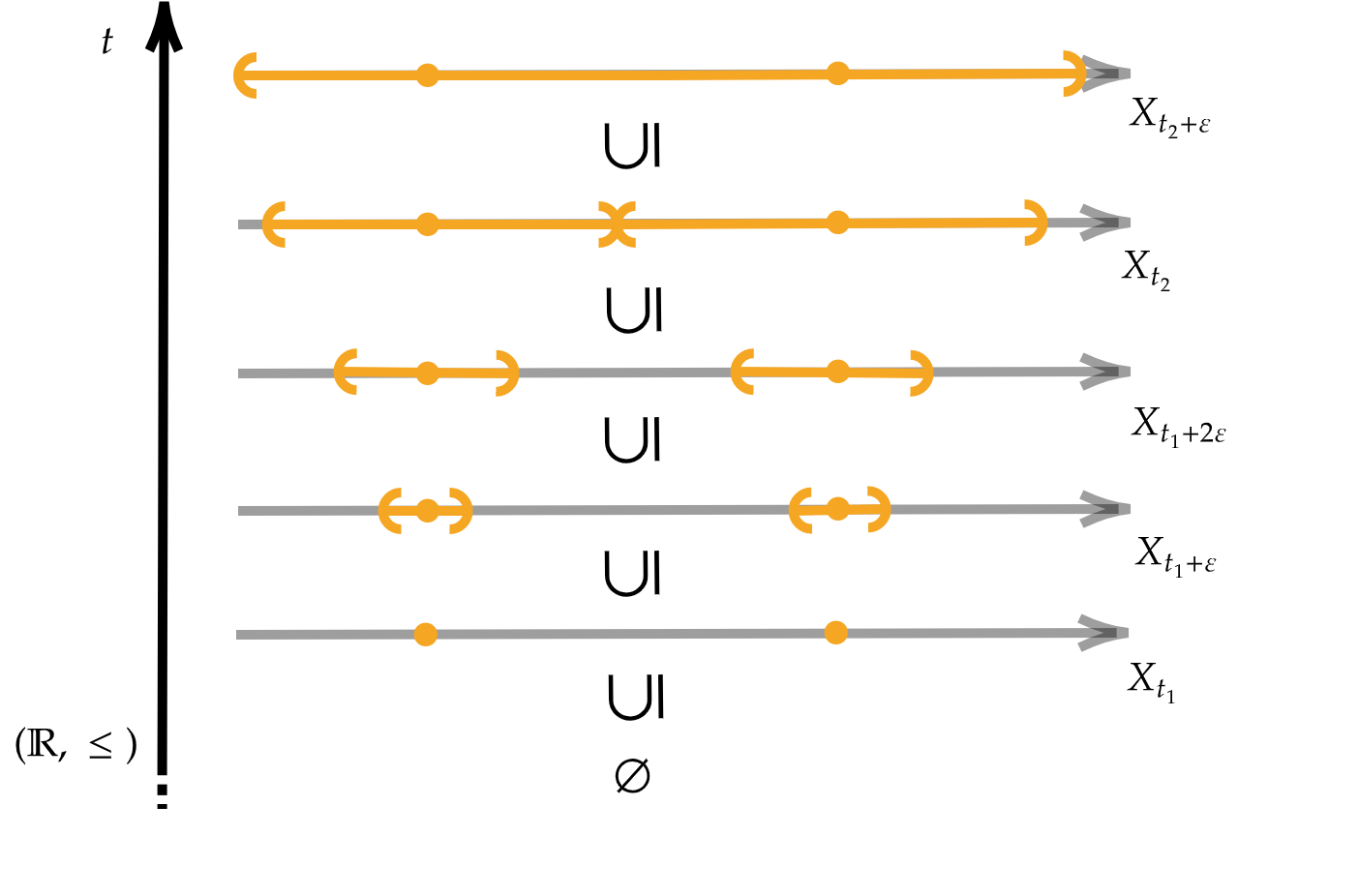}
    	\caption{A filtration $\X$.}
    	\label{fig:filtr}
	\end{subfigure}
%	\hspace{0.5 cm}
    \begin{subfigure}[c]{0.49\textwidth}
    	\centering
	    \includegraphics[width = \textwidth]{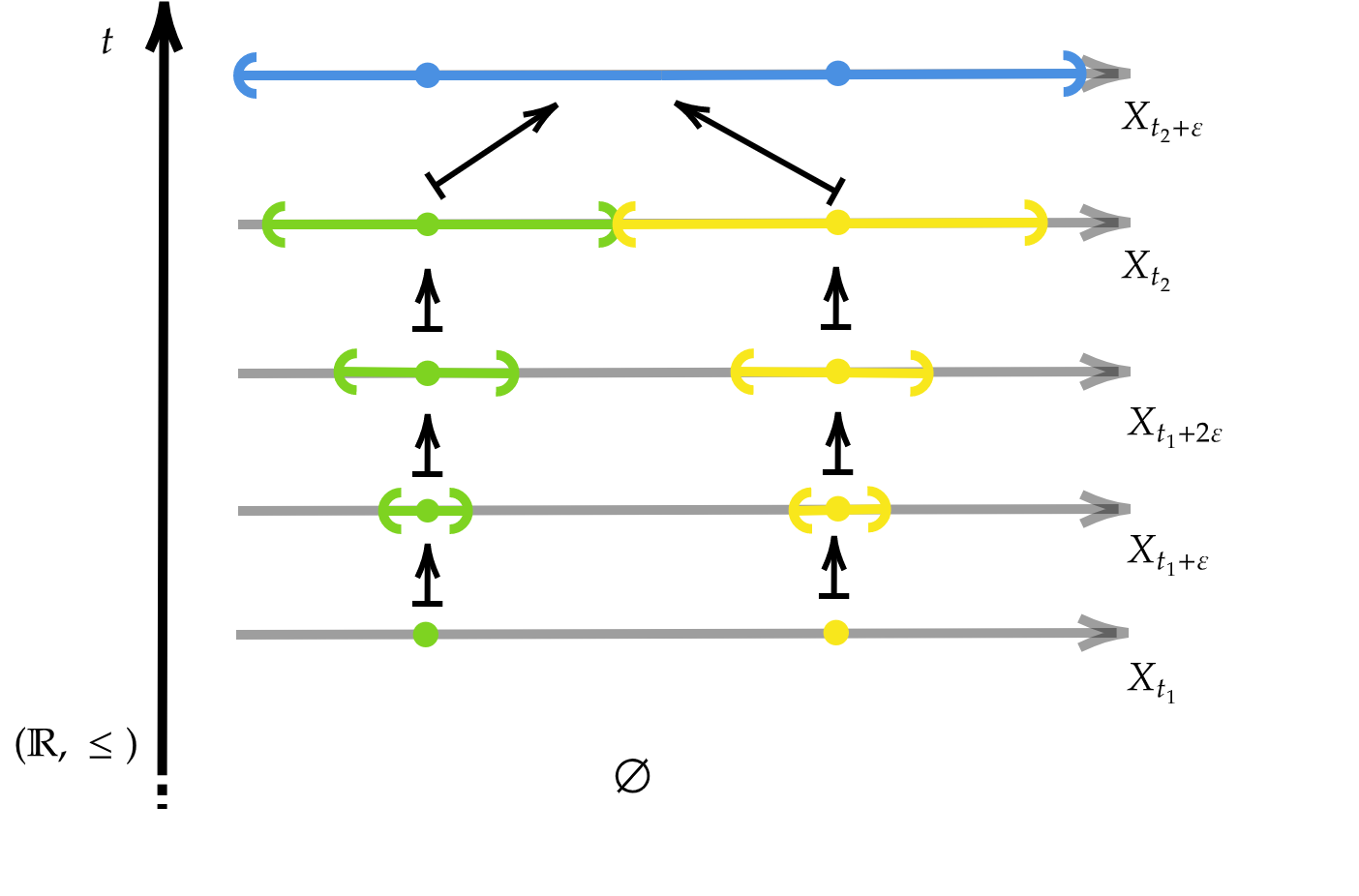}
    	\caption{The abstract merge tree $\T$.}
    	\label{fig:abs_MT}
	\end{subfigure}
%	\hspace{0.5 cm}

    \begin{subfigure}[c]{0.49\textwidth}
    	\centering
	    \includegraphics[width = \textwidth]{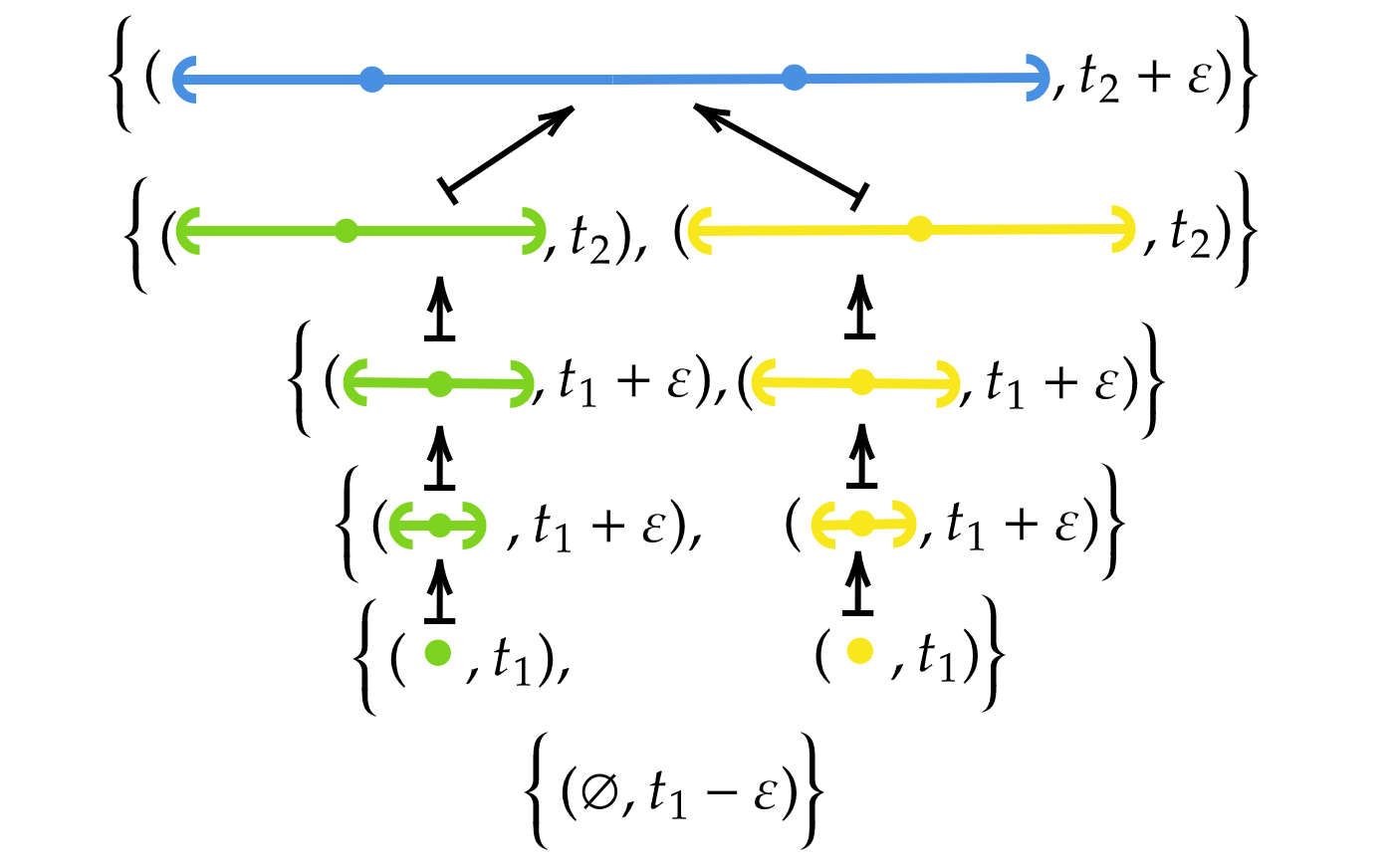}
    	\caption{The display poset $D_{\T}$.}
    	\label{fig:display}
	\end{subfigure}

\caption{An example of a filtration along with its abstract merge tree and its display poset. The colors are used throughout the plots to highlight the relations between the different objects.}
\label{fig:preliminary}
\end{figure}

We begin by introducing the main mathematical objects central to our research, grounding them in the existing scientific literature. Where standard notation is lacking or ambiguous, we propose new, well-motivated definitions to avoid confusion and prevent conflicts with established terminology.
 \Cref{fig:preliminary} illustrates some of the objects we introduce in this section.

\begin{defi}[\cite{curry2021decorated}]
A filtration of topological spaces is a (covariant) functor 
$X_{\Bigcdot}:\mathbb{R}\rightarrow \Top$ from the poset $(\mathbb{R},\leq )$ to $\Top$, the category of topological spaces with continuous functions, such that:  
$X_{t\leq t'}:X_t \rightarrow X_{t'}$, for $t\leq t'$,
are injective maps. 
\end{defi}

\subparagraph*{Example}
Given a real valued function $f:X\rightarrow \mathbb{R}$ the \emph{sublevel set} filtration is given by $X_{t}=f^{-1}((-\infty,t])$ and $X_{t<t'}=i:f^{-1}((-\infty,t])\hookrightarrow f^{-1}((-\infty,t'])$.

\subparagraph*{Example}
Given a finite set $C\subset \mathbb{R}^n$ its  offset filtration is given by $X_{t}=\bigcup_{c\in C}B_{t}(c)$, with $B_{t}(c)=\{x\in\mathbb{R}^n \mid \parallel c-x \parallel < t\}$. As before: $X_{t<t'}=i:\bigcup_{c\in C}B_{t}(c)\hookrightarrow \bigcup_{c\in C}B_{t'}(c)$.

\bigskip

Given a filtration $\X$ we can compose it with the functor $\pi_0$ sending each topological space into the set of its path-connected components, obtaining a functor with values in $\Sets$, the category of sets, with arrows being functions. We will also impose finiteness conditions that will require that the functor only takes values in $\FSet$, the category of finite sets.

We recall that, according to standard topological notation, $\pi_0(X)$ is the set of the path-connected components of $X$ and, given a continuous function $q:X\rightarrow Y$, $\pi_0(q):\pi_0(X)\rightarrow \pi_0(Y)$ is defined as:

\[
U\mapsto V  \text{ such that }q(U)\subset V.
\]

\begin{defi}[\cite{carlsson2013classify, curry2018fiber}]
A persistent set is a functor $S:\mathbb{R}\rightarrow \Sets$. In particular, given a filtration of topological spaces $\X$, the persistent set of components of $\X$ is $\pi_0 \circ \X$.
\end{defi}

Note that, by endowing a persistent set with the discrete topology, every persistent set can be seen as the persistent set of components of a filtration. Thus, a general persistent set $S$ can be written as $\pi_0(\X)$ for some filtration $\X$. 

Now we want to carry on, going towards the definition of merge trees.
The existing paths for giving such notion relying on the language of TDA split at the definition of persistence module. All such approaches however share similar notions of \emph{constructible} persistent sets \cite{patel2018generalized} or modules \cite{curry2021decorated}.
We give here the definition of constructible persistent sets adapted from \cite{patel2018generalized}. The original definition in \cite{patel2018generalized} is stated for persistence modules (as defined in \cite{patel2018generalized}) and it is slightly different - see \Cref{rmk:constructible_defi}. 

\begin{defi}[modified from \cite{patel2018generalized}]
\label{def:construct}
A persistent set $S:\mathbb{R}\rightarrow \Sets$ is constructible if there is a finite collection of real numbers $\{t_1<t_2<\ldots<t_n\}$ such that:
\begin{itemize}
\item $S(t<t')=\emptyset$ for all $t<t_1$;
%\item $S(t)=\{\star\}$ for all $t>t_n$;
\item for $t,t'\in (t_i,t_{i+1})$ or $t,t'>t_n$, with $t<t'$, then $S(t<t')$ is bijective.
\end{itemize}
The set $\{t_1<t_2<\ldots<t_n\}$ is called a critical set and $t_i$ are called critical values. If $S(t)$ is always a finite set, then $S$ is called a finite persistent set.
\end{defi}

\begin{rmk}\label{rmk:constructible_defi}
In the literature there is not an univocal  way to treat critical values: in \cite{de2016categorified}, Definition 3.3, constructibility conditions are stated in terms of open intervals (due to the use of cosheaves); in \cite{patel2018generalized}, Definition 2.2, all the conditions are stated in terms of half-closed intervals $[t_i,t_{i+1})$; while \cite{curry2021decorated} differentiates between  the open interval $(t_n,+\infty)$ i.e. $t,t'>t_n$, and the half closed intervals $[t_i,t_{i+1})$. For reasons which will be detailed in \Cref{sec:critical_values}, we stated all the conditions following \cite{de2016categorified}, with open intervals, thus relaxing the conditions presented in \cite{curry2021decorated}.
\end{rmk}

At this point we highlight two different categorical approaches to obtain merge trees. The author of 
\cite{patel2018generalized} requires a persistence module to be a functor
$F:\mathbb{R}\rightarrow C$ with $C$ being an essentially small symmetric monoidal category with images (see \cite{patel2018generalized} and references therein). 
If then one wants to work with values in some category of vector spaces over some field $\mathbb{K}$, it is required that $F(t)$ is always finite dimensional. A merge tree, see \cite[Example 2.1]{patel2018generalized}, is then a constructible persistence module with values in $\FSet$, the category of finite sets.

The authors of \cite{curry2021decorated} instead, state that a persistence module is a functor 
$F:\mathbb{R}\rightarrow \Vect_\mathbb{K}$, with 
$\Vect_\mathbb{K}$ being the category of vector spaces over the field $\mathbb{K}$. This definition seems to be in line with the ones given by other works, especially in multidimensional persistence (see for instance  \cite{scola2017noise} and references therein). On top of that, \cite{curry2021decorated} obtains a (generalized) merge tree as the \emph{display poset} (see \Cref{defi:display}) of a persistent set. The constructibility condition on the persistent set then implies the merge tree to be  \emph{tame}. 

We find natural to work with objects which are functors, as the merge trees defined in \cite{patel2018generalized}, but we require some properties which are closer to the ones of constructible persistent sets, as in \Cref{def:construct}.
Thus, mixing those definitions, we give the notion of an \emph{abstract merge tree}.

\begin{defi}\label{def:abstract_mt}
An abstract merge tree is a persistent set $S:\mathbb{R}\rightarrow \Sets$ such that there is a finite collection of real numbers $\{t_1<t_2<\ldots<t_n\}$ which satisfies:
\begin{itemize}
\item $S(t)=\emptyset$ for all $t<t_1$;
\item $S(t)=\{\star\}$ for all $t>t_n$;
\item if $t,t'\in (t_i,t_{i+1})$, with $t<t'$, then $S(t<t')$ is bijective.
\end{itemize}
The values $\{t_1<t_2<\ldots<t_n\}$ are called critical values of the tree.

If $S(t)$ is always a finite set, $S$ is a finite abstract merge tree. 
\end{defi}

\begin{assump}
From now on we will be always working with finite abstract merge trees and, to lighten the notation, we assume any abstract merge tree to be finite, without explicit reference to its finiteness.
\end{assump}

We point out that two abstract merge trees $\pi_0(\X)$ and $\pi_0(\Y)$ are isomorphic if there is a natural transformation $\alpha_t:\pi_0(X_t)\rightarrow \pi_0(Y_t)$ which is bijective for every $t$. This is equivalent to having the same number of path-connected components for every $t$ and having bijections which make the following square commute:

\[
\begin{tikzcd}
X_{t}\ar[d,"\alpha_{t}"]\ar[r]&X_{t'}\ar[d,"\alpha_{t'}"]\\
Y_{t}\ar[r]&Y_{t'}
\end{tikzcd}
\]
 for all $t<t'$.

We take one last definition from  \cite{curry2021decorated} which will be needed in later sections. 

\begin{defi}[\cite{curry2021decorated}]\label{defi:display}
Given a persistent set $S:\mathbb{R}\rightarrow \Sets$ we define its display poset as:
\[
D_S:=\bigcup_{t\in\mathbb{R}}S(t)\times\{t\}. 
\]
The set $D_S$ can be given a partial order with $(a,t)\leq (b,t')$ if $S(t\leq t')(a)=b$. 
\end{defi}

Given a persistent set $S$ and its display poset $D_S$ we define $h((a,t))=t$ and $\pi((a,t))=a$ for every $(a,t)\in D_S$. From $D_S$ we can clearly recover $S$ via $S(t)=\pi(h^{-1}(t))$ and $S(t\leq t')(a)=b$ with $a\leq b$.
Thus the two representations are equivalent and, at any time, we will use the one which is more 
convenient for our purposes.
Note that this construction is functorial:  any natural transformation $\eta:S\rightarrow S'$ between persistent sets, gives a map of sets $D_\eta: D_S\rightarrow D_{S'}$ with 
$D_\eta((a,t))=(\eta_t(a),t)$. Clearly $D_{\eta\circ \nu}=D_\eta \circ D_\nu$.

\subsection{Critical Values}
\label{sec:critical_values}

Before bridging between abstract merge trees and merge trees, we need to focus on some subtle facts about critical values: 1) neither in \Cref{def:construct}
nor in \Cref{def:abstract_mt} critical values are uniquely defined 2) we decided to relax the constructibility conditions of \cite{curry2021decorated} to account for scenarios which would be otherwise excluded from their framework. In particular, we will show that with a coherent set of definitions we can meaningfully reduce to the setting of \cite{curry2021decorated} even in more general situations.

Thanks to the functoriality of persistent sets, we immediately solve the first point: we can  take the intersection of all the possible sets of critical values to obtain a minimal one, which is possibly empty.

\begin{prop}\label{prop:intersect_critical}
Let $S$ be a constructible persistent set and let $\{C_{i}\}_{i\in I}$ be a family of finite critical sets of $S$. Then $C=\bigcap_{i\in I} C_i$ is a critical set. 
\begin{proof}
Clearly $C$ is a finite set, possibly empty. The thesis is then a consequece of the following fact: if $t\notin C$ then there is $\varepsilon>0$ such that $S(t-\varepsilon<t+\varepsilon)$ is bijective. So we can remove $t$ from any critical set of $S$ and still obtain a critical set.
\end{proof}
\end{prop}

\begin{assump}
Leveraging on \Cref{prop:intersect_critical}, any time we take any abstract merge tree or a constructible persistent set and consider its critical values, we mean the elements of the minimal critical set.
\end{assump}

Consider an abstract merge tree $\T$ and let $t_1<t_2<\ldots<t_n$ be its (minimal set of) critical values.
Let $i_{t}^{t'}:= X_{t\leq t'}:X_t\rightarrow X_{t'}$.
Given a critical value $t_j$, due to the minimality condition, we know that for $\varepsilon>0$ small enough, at least one between  $\pi_0(i_{t_j-\varepsilon}^{t_j})$ and $\pi_0(i_{t_j}^{t_j+\varepsilon})$ is not bijective.

We want to distinguish between two scenarios: 
\begin{itemize}
\item if $\pi_0(i_{t_j}^{t_j+\varepsilon})$ is bijective, we say that topological changes in the persistent set (and in the filtration) happen \emph{at} $t_j$;
\item  if $\pi_0(i_{t_j}^{t_j+\varepsilon})$ is not bijective, we say that topological changes in the persistent set (and in the filtration) happen \emph{across} $t_j$.
\end{itemize}

\begin{figure}
    	\centering
    	\includegraphics[width = \textwidth]{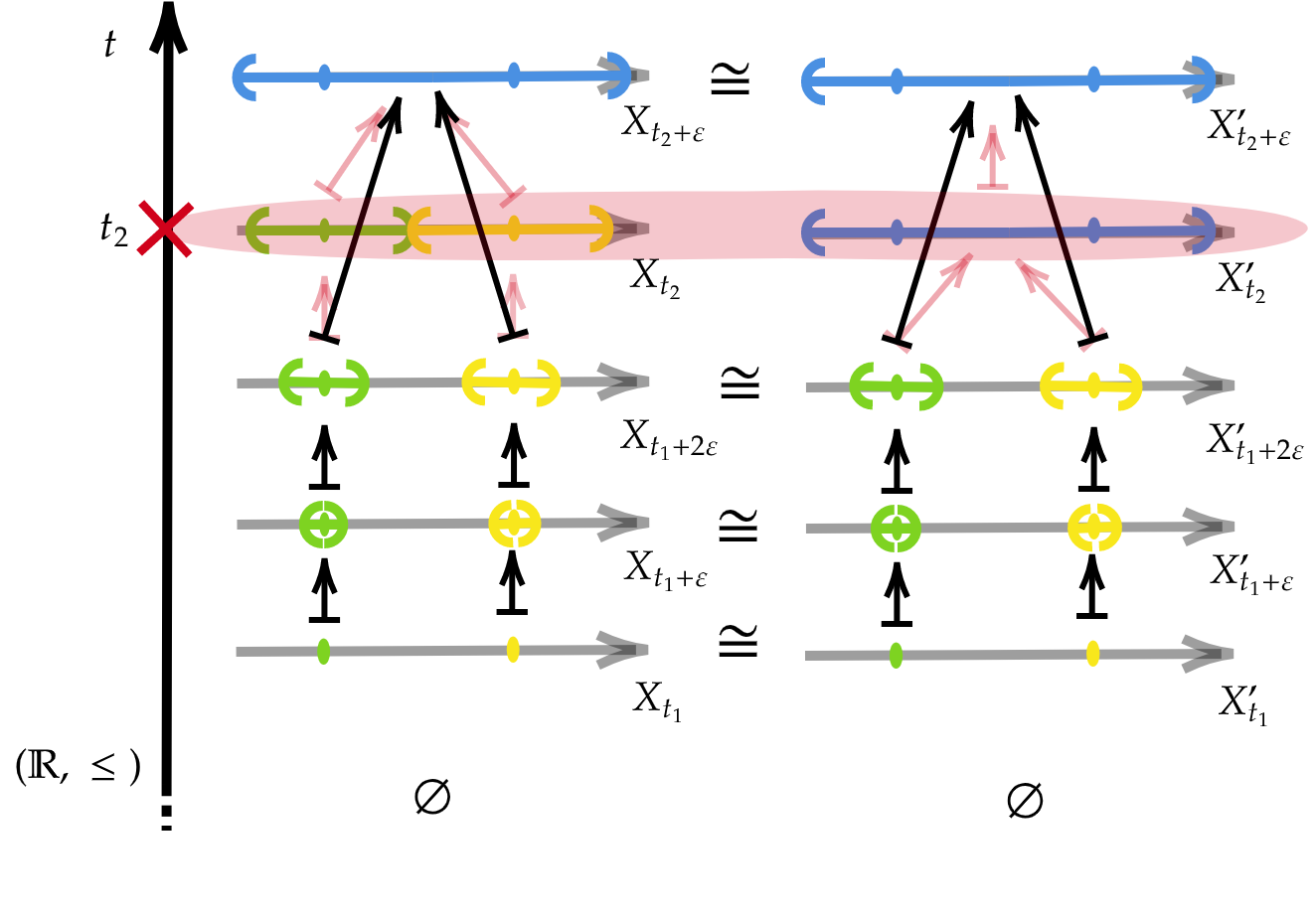}
\caption{An abstract merge tree $\T$ which is not regular (left), along with the regular abstract merge tree $R(\T)$ (right), obtained as in \Cref{prop:regular}. The horizontal isomorphism signs represent the a.e. isomorphism between them: they are isomorphic on $\mathbb{R}-\{t_2\}$. The colors are again used to highlight the relations between the different objects.}
\label{fig:regular}
\end{figure}

The constructibility conditions in \cite{curry2021decorated} are stated so that topological changes always happen at critical values. Accordingly, we give the following definition.
See also \Cref{fig:regular}.

\begin{defi}
A constructible persistent set $\T$ is said to be regular if all topological changes happen at its critical values.
\end{defi}

Consider the following filtrations
of topological spaces: $X_t = (-t,t)\bigcup (1-t,1+t)$ and $Y_t = [-t,t]\bigcup [1-t,1+t]$ for $t>0$ and $X_0=Y_0=\{0,1\}$. For $t<0$ the filtrations are empty.
The situation arising from this filtration is similar to what is depicted in \Cref{fig:regular}: the persistent sets $\pi_0(X_t)$ and $\pi_0(Y_t)$ are two abstract merge trees and they share the same set of critical values, namely $\{0,1/2\}$. They only differ at the critical value $1/2$: $\pi_0(X_{1/2})=\{(-1/2,1/2),(1/2,1)\}$, while 
$\pi_0(Y_{1/2})=\{(-1/2,1)\}$. In $\X$ changes happen across the critical values - $\pi_0(X_{1/3})\cong \pi_0(X_{1/2})$ and $\pi_0(X_{1/2})\ncong \pi_0(X_{1})$, while in $\Y$ changes happen at the critical values - $\pi_0(Y_{1/3})\ncong \pi_0(Y_{1/2})$ and $Y_{1/2}\cong Y_{1}$.

It is then clear that $\X$ and $\Y$ are not isomorphic as abstract merge trees (as $\pi_0(X_{1/2})\ncong \pi_0(Y_{1/2})$), but, at the same time, they differ only by their behavior at critical values. 
Filtrations like $\X$ appear in many interesting situations (see also  \cite{pegoraro2024functional}, Appendix A) which we don't want to exclude from our framework. 
At the same time, we are not interested in distinguishing between $\X$ and $\Y$, and for this reason, we ask for a weaker notion of equivalence between abstract merge trees.

Given $Z\subset \mathbb{R}$, clearly $Z$ inherits an ordering from the one in $\mathbb{R}$ and we can consider $Z$ as a poset category. Thus, we can take the restriction to $Z$ of any filtration of topological spaces $\X$ (and similarly of any persistent set) via the inclusion $Z\hookrightarrow \mathbb{R}$. We indicate this restriction as $X_{\Bigcdot\mid Z}$. Moreover, let $\mathcal{L}$ be the Lebesgue measure on $\mathbb{R}$. Refer to \Cref{fig:regular} for a visual interpretation of the following definitions and propositions.

\begin{defi}
Two persistent sets $\pi_0(\X)$ and $\pi_0(\Y)$ are almost everywhere (a.e.) isomorphic if there is a Lebesgue measurable set $Z\subset \mathbb{R}$ such that $\mathcal{L}(\mathbb{R}-Z)=0$ and there is a natural isomorphism 
$\alpha: \pi_0(X_{\Bigcdot\mid Z}) \rightarrow \pi_0(Y_{\Bigcdot\mid Z})$. We write $\pi_0(\X)\cong_{a.e}\pi_0(\Y)$.
\end{defi}

\begin{prop}\label{prop:a.e._eq_rel}
Being a.e. isomorphic is an equivalence relation between persistent sets.
\begin{proof}
Reflexivity and symmetry are trivial: the first one holds with $Z=\emptyset$ and the second one holds by definition of natural isomorphism. Lastly, transitivity holds because any finite union of measure zero sets is a measure zero set.
\end{proof}
\end{prop}

Now we prove that in each equivalence class of a.e. isomorphic abstract merge trees we can always pick a regular abstract merge tree, which is unique up to isomorphism.

\begin{prop}\label{prop:regular}
For every abstract merge tree $\pi_0(\X)$ there is a unique (up to isomorphism) abstract merge tree $R(\pi_0(\X))$ such that:
\begin{enumerate}
\item $\pi_0(\X)\cong_{a.e.} R(\pi_0(\X))$;
\item $R(\pi_0(\X))$ is regular.
\end{enumerate}
\begin{proof}
    See \Cref{sec:proofs_TDA}.
\end{proof}
\end{prop}

Regular abstract merge trees make many upcoming definitions and results more natural and straightforward.
With \Cref{prop:regular} we formally state that this choice is indeed consistent with the equivalence relation previously established, and we can resort to regular abstract merge trees without excluding non-regular scenarios.

\begin{figure}
    	\centering
    \begin{subfigure}[c]{0.49\textwidth}
	    \includegraphics[width = \textwidth]{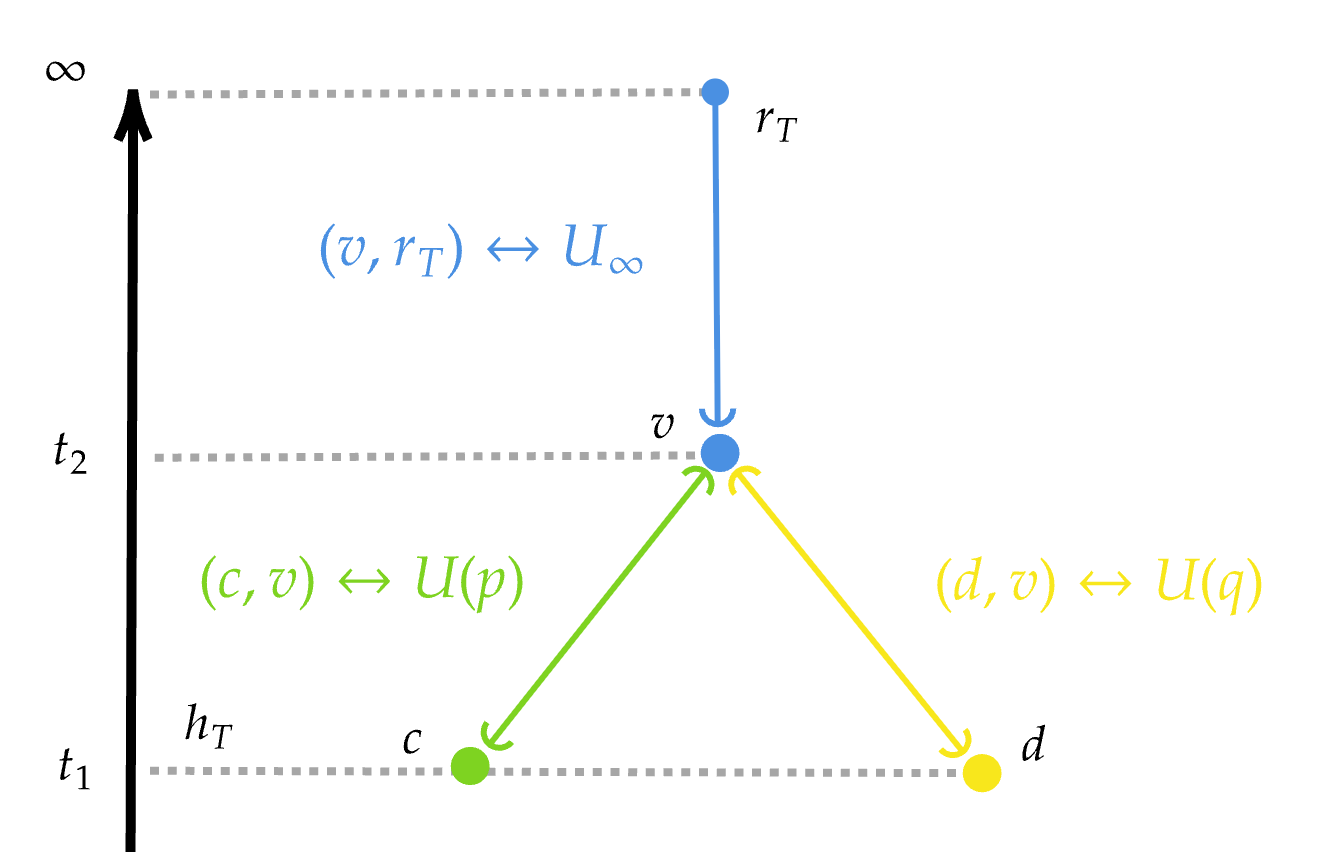}
            \captionsetup{singlelinecheck=off, margin={0.3cm, 0.1cm}}
    	\caption{The merge tree $\mathcal{M}(R(\T))$ associated to the abstract merge trees in \Cref{fig:regular}. The brackets at the end of the edges and the labels $U(p),U(q),U_\infty$ refer to the canonical a.e. covering defined in \Cref{sec:metric_display}.}
    	\label{fig:MT}
	\end{subfigure}		
    \begin{subfigure}[c]{0.49\textwidth}
    	\centering
	    \includegraphics[width = \textwidth]{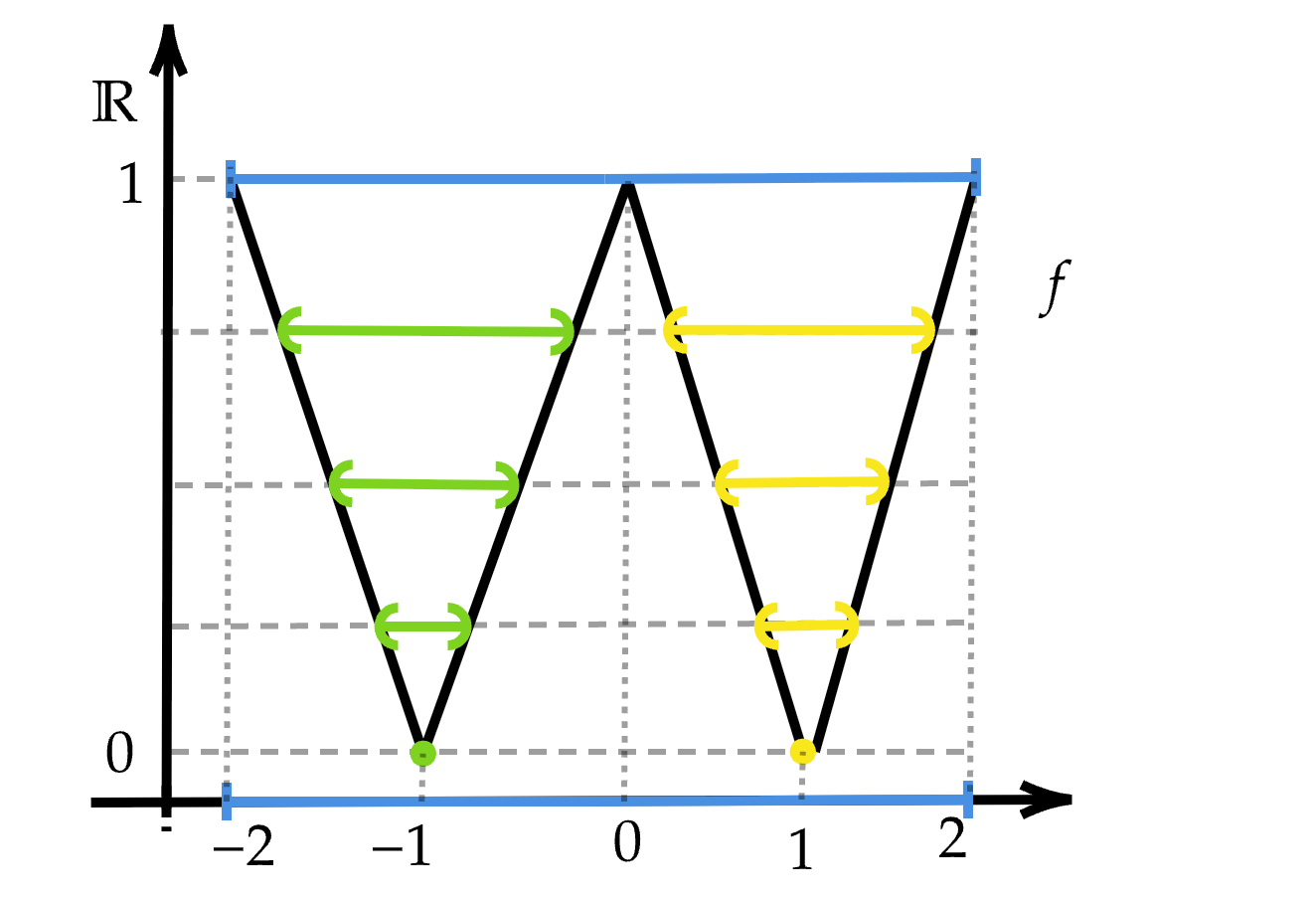}
            \captionsetup{singlelinecheck=off, margin={0.3cm, 0.1cm}}
    	\caption{The function $f(x)=\mid \mid x\mid -1 \mid$, along with its sublevel set filtration $\X$ and the associated persistent set $\T$. The merge tree $\mathcal{M}(R(\T))$ is isomorphic to the one in \Cref{fig:MT} upon taking $t_1=0$ and $t_2=1$. }
    	\label{fig:example}
	\end{subfigure}	

    \caption{Additional figures to visually interpret the content of \Cref{sec:intro_MT} and \Cref{sec:back_to_vec}.}
\end{figure}

\section{Merge Trees}
\label{sec:intro_MT}

We now introduce the discrete counterpart of abstract merge trees. These structures are referred to simply as \emph{merge trees} in much of the relevant literature \cite{merge_intrins, merge_farlocca}, while \cite{curry2021decorated} uses the term \emph{computational} merge trees. Although we appreciate the motivation behind the latter terminology, we adopt the former for consistency with the majority of the literature and for simplicity, given that these objects form the main focus of the theoretical developments in this paper.
%Before proceeding, we point out that there is a third approach - on top of the categorical and the computational ones - to the definition of merge trees, followed for instance by \cite{merge_interl}, which in \cite{curry2021decorated} is referred to as \emph{classical} merge trees. 
%We avoid dealing with such objects in our dissertation and any interested reader can find in \cite{curry2021decorated} how that definition relates with the other ones we report. 

%Now we need some graph-related definitions.

\begin{defi}\label{def:tree_struct}
A tree structure $T$ is given by a set of vertices $V_T$ and a set of edges $E_T\subset V_T\times V_T$ which form a connected rooted acyclic graph.  We indicate the root of the tree with  $r_T$. We say that \(T\) is finite if \(V_T\) is finite. The degree of a vertex $v\in V_T$ is the number of edges which have that vertex as one of the endpoints, and is denoted $ord_T(v)$. 
Any vertex with an edge connecting it to the root is its child and the root is its parent: this is the first step of a recursion which defines the parent and children relation for all vertices in \(V_T.\)
%In this way we recursively define parent and children (possibly none) for any vertex on the tree. 
The vertices with no children are called leaves  or taxa and are collected in the set $L_T$. The relation $child < parent $ generates a partial order on $V_T$. The edges in $E_T$ are identified by ordered pairs $(a,b)$ with $a<b$.
A subtree of a vertex $v$, called $sub_T(v)$, is the tree structure whose set of vertices is $\{x \in V_T\mid  x\leq v\}$. 
\end{defi}

Note that, given a tree structure $T$, identifying an edge $(v,v')$ with its lower vertex $v$, gives a bijection between $V_T-\{r_T\}$ and $E_T$, that is $E_T\cong V_T-\{r_T\}$ as sets. 
Given this bijection, we often use $E_T$ to indicate the vertices $v\in V_T-\{r_T\}$,  to simplify the notation.

We want to identify merge trees independently of their vertex set, and thus we introduce the following isomorphism classes.

\begin{defi}
Two tree structures $T$ and $T'$ are isomorphic if there exists a bijection $\eta:V_T\rightarrow V_{T'}$ that induces a bijection between the edge sets $E_T$ and $E_{T'}$: $(a,b)\mapsto (\eta(a),\eta(b))$. Such $\eta$ is an isomorphism of tree structures.
\end{defi}

Finally, we give the definition of a merge tree, slightly adapted from \cite{merge_intrins}.

\begin{defi}\label{defi:merge_TDA}
A merge tree is a finite tree structure $T$ with a monotone increasing (w.r.t. the poset structure of $V_T$) height function $h_T:V_T\rightarrow \mathbb{R}\cup\{+\infty\}$ and such that 1) $ord_T(r_T)=1$ 2) $h_T(r_T)=+\infty$ 3) $h_T(v)\in\mathbb{R}$ for every $v<r_T$.

Two merge trees $(T,h_T)$ and $(T',h_{T'})$ are isomorphic if $T$ and $T'$ are isomorphic as tree structures and the isomorphism $\eta:V_T\rightarrow V_{T'}$ is such that $h_T = h_{T'} \circ \eta$. Such $\eta$ is an isomorphism of merge trees. We use the notation $(T,h_T)\cong (T',h_{T'})$.
\end{defi}

With some slight abuse of notation we set $\max h_T =\max_{v\in V_{T}\mid v<r_T}h_{T}(v)$ and $\arg \max h_T =\max\{v\in V_{T}\mid v<r_T\}$. Note that, given $(T,h_T)$ merge tree, there is only one edge of the form $(v,r_T)$ and we have $v=\arg\max h_T$.

The connection between abstract merge trees and merge trees is clarified in \Cref{sec:ab_mt_vs_mt}; however, before proceeding, we need to introduce an additional equivalence relation on merge trees.

\begin{defi}\label{defi:ghosting}
Given a tree structure $T$, we can eliminate a  degree two vertex, connecting the two adjacent edges  which share such vertex as endpoint. 
Suppose we have two edges $e=(v_1,v_2)$ and $e'=(v_2,v_3)$, with $v_1<v_2<v_3$. And suppose $v_2$ is of degree two. Then, we can remove $v_2$ and merge $e$ and $e'$ into a new edge $e''=(v_1,v_3)$.
This operation is called the \emph{ghosting} of the vertex $v_2$. Its inverse transformation, which restores the original tree, is called a \emph{splitting} of the edge $e''$.
Similarly, given a merge tree, by ghosting vertices one obtains a new merge tree  with the height function on the new merge tree being obtained by restricting the height function of the old tree to the remaining vertices.
\end{defi}

Now we can state the following definition.

\begin{defi}
Merge trees are equal up to degree $2$ vertices if they become isomorphic after applying a finite number of ghostings or splittings. We write $(T,h_T)\cong_2(T',h_{T'})$. 
\end{defi}

\subsection{Regular Abstract Merge Trees and Merge Trees}
\label{sec:ab_mt_vs_mt}

In this section we study the relation between abstract merge trees and merge trees. We collect all the important facts on this topic in the following proposition. 
\Cref{fig:abs_MT} and \Cref{fig:MT} can help the reader going through the upcoming results.

\begin{prop}\label{prop:equivalence}
The following hold:
\begin{enumerate}
\item we can associate to a regular abstract merge tree $R(\T)$, a merge tree without degree $2$ vertices $\mathcal{M}(R(\T))$;
\item we can associate to a merge tree $(T,h_T)$, a regular abstract merge tree $\mathcal{F}((T,h_T))$. Moreover, we have $\mathcal{M}(\mathcal{F}((T,h_T)))\cong_2 (T,h_T)$ and $\mathcal{F}(\mathcal{M}(R(\T))\cong_{a.e.} \T$;
\item given two abstract merge trees $\X$ and $\Y$, $\mathcal{M}(R(\T))\cong \mathcal{M}(R(\G))$ if and only if $\T\cong_{a.e} \G$.
\item given two merge trees $(T,h_T)$ and $(T',h_{T'})$, we have $\mathcal{F}((T,h_T))\cong \mathcal{F}((T,h_T))$ if and only if $(T,h_T)\cong_2 (T',h_{T'})$.
\end{enumerate}
\begin{proof}
    See \Cref{sec:proofs_TDA}.
\end{proof}
\end{prop}

We point out an additional fact about degree $2$ vertices.
Suppose that we were to remove a leaf in a merge tree, the parent of the deleted vertex may become a degree two vertex. In case that happens, such vertex carries no topological information, since the merging that the point
was representing, is no longer happening (was indeed removed). And in fact the abstract merge tree associated to the merge tree with the degree $2$ vertex and to the merge tree with the degree $2$ vertex ghosted are the same
by \Cref{prop:equivalence}. 
Thus working up to degree two vertices is a very natural framework to work with merge trees. And this must be taken into consideration when setting up the framework to deal with functions defined on merge trees.

The proof of \Cref{prop:equivalence} carries this important corollary.

\begin{cor}\label{cor:embed}
Given a merge tree $(T,h_T)$ and the abstract merge tree $\T=\mathcal{F}((T,h_T))$, we have 
$E_T\hookrightarrow D_{\T}$ (recall that $E_T\cong V_T-\{r_T\}$) induced by the map $v\mapsto (v,h_T(v))$.
\end{cor}

With the help of \Cref{cor:embed}, we also give the following definition.

\begin{defi}\label{defi:crit_pts}
Given an abstract merge tree $\T=\mathcal{F}((T,h_T))$, the critical points of $D_{\T}$ are defined as the image of 
$E_T\hookrightarrow D_{\T}$, with 
    $ T = \mathcal{M}(R(\T))$.
\end{defi}

\subsection{Example of Merge Tree}
\label{sec:easy_example}

Now we briefly give an example of a merge tree representing the merging structure of path-connected components along the sublevel set filtration of a function. The reader should refer to \Cref{fig:example} for a visual interpetation and to \Cref{sec:main_ideas} 
for more examples, which also propel the use of merge trees over persistence diagrams.

%\subsubsection{Point Clouds}
%
%Consider $X_t:= (-t,t)\bigcup (1-t,1+t)$ for $t> 0$. We have $X_t: =H_0(X_t)= \mathbb{K}^2$ for $t\in [0,1/2]$ and $X_t =H_0(X_t)= \mathbb{K}$ for $t>1/2$. If we consider $Y_t:= [-t,t]\bigcup [1-t,1+t]$ and $Y_t: =H_0(Y_t)$, we have $Y_t= \mathbb{K}^2$ for $t\in [0,1/2)$ and $Y_t =\mathbb{K}$ for $t\geq 1/2$. The two filtrations  
%$\{X_t\}$ (whenever we build a filtration with $H_0(\cdot)$ and do not specify a basis, we imply that we fix the one given by path-connected components) and $\{Y_t\}$ share the same set of critical values, namely $\{0,1/2\}$ and they only differ by the vector spaces at the critical value $1/2$: $X_{1/2}=\mathbb{K}^2$, while 
%$Y_{1/2}=\mathbb{K}$. In $\{X_t\}$ changes happen across the critical values - $X_{1/3}= X_{1/2}$ and $X_{1/2}\ncong X_{1}$, while in $\{Y_t\}$ changes happen at the critical values - $Y_{1/3}\ncong Y_{1/2}$ and $Y_{1/2}= Y_{1}$.

Consider the function $f=\mid \mid x\mid -1\mid $ defined on the interval $[-2,2]$. Consider the sublevel set filtration $X_t=f^{-1}((-\infty,t])$.
The sublevel set $X_t$ is an interval of the form $[-1-t,-1+t]\bigcup [1-t,1+t]$, for $t\in [0,1]$.

Consider then the abstract merge tree $\T$.
For any $t\in [0,1)$, the path-connected components are $a_t=\{a^t_1,a^t_{-1}\}$, with $a^t_1=[1-t,1+t]$ and $a^t_{-1}=[-1-t,-1+t]$ and for $t\geq 1$, $a^{t}_{2}=\{[-2,2]\}$.
The critical values of the filtration are $t_1=0$ and $t_2=1$. The maps are $a_i^t\mapsto a_i^{t'}$ and for $i=-1,1$, for $t\leq t' <1$; $a_1^t,a_{-1}^t\mapsto a_2^{t'}$ for $t<1\leq t' $ and the identity for $t,t'\geq 1$.  

The merge tree $\mathcal{M}(\T)=(T,h_T)$ associated to $\T$ has a tree structure given by a root, an internal vertex and two leaves - as in \Cref{fig:MT}: if we call $v_1:=a^0_{1}$, $v_{-1}:=a^0_{-1}$ and $v_{2}:=a^2_{2}$, the merge tree $\mathcal{M}(\T)$ is given by the vertex set $\{v_1,v_{-1},v_2,r_T\}$ and edges $e_1=(v_1,v_2)$, $e_2=(v_{-1},v_2)$ and $e_3=(v_2,r_T)$. The height function has values $h_T(v_1)=h_T(v_{-1})=t^-=0$, $h_T(v_2)=1$ and $h_T(r_T)=+\infty$.

\section{Functions Defined on Display Posets}\label{sec:back_to_vec}

Now we formalize how we want to deal with functions defined on merge trees, devoting much care to setting up a framework in accordance with the equivalence relations introduced in \Cref{sec:critical_values}. 
We build a topology (in fact, a pseudo-metric) on merge trees and a measure. These are very natural constructions and can be identified, respectively, with the quotient topology when merge trees are built from functions as in \cite{merge_interl}, and the pullback of the Lebesgue measure on $\R$ via the map $D_{\T}\rightarrow \R$.

\subsection{Metric Spaces}
\label{sec:defi_metric_spaces}

Following \cite{burago2022course}, we briefly state the definitions related to metric geometry that we need in the present work.

\begin{defi}
Let $X$ be an arbitrary set. A non-negative function $d:X\times X \rightarrow \R$ is a (finite) pseudo metric if for all $x,y,z\in X$ we have:
\begin{enumerate}
\item $d(x,x)=0$
\item $d(x,y)=d(y,x)$
\item $d(x,y)\leq d(x,z)+d(z,y)$.
\end{enumerate}
The space $(X,d)$ is called a pseudo metric space. A pseudo metric space is a topological space with the topology generated by the open balls $B_\varepsilon(x):=\{y\in X \mid d(x,y)<\varepsilon\}$.

Given a pseudo metric $d$ on $X$, if for all $x,y\in X$, $x\neq y$, we have $d(x,y)>0$ then $d$ is called a metric or a distance and $(X,d)$ is a metric space.
\end{defi}

\begin{prop}[Proposition 1.1.5 \cite{burago2022course}]
For a pseudo metric space $(X,d)$, the following is an equivalence relation: $x\sim y$ iff $d(x,y)=0$. Moreover, the quotient space $(X,d)/\sim$ is a metric space.
\end{prop}

\begin{defi}
An isometric embedding between two metric spaces $(X,d_X)$ and $(Y,d_Y)$ is an injective function $f:X\rightarrow Y$ such that $d_X(x,y)=d_Y(f(x),f(y))$. If $f$ is also bijective then it is called an isometry (or an isometric isomorphism).  
\end{defi}

\subsection{The Display Poset as a Pseudo-Metric Space}
\label{sec:metric_display}

Now we build function spaces on display posets.
We begin by giving the notion of \emph{common ancestors} for subsets of the display poset of an abstract merge tree.

\begin{defi}
Given $Q\subset D_{\T}$, with $\sup h(Q)<\infty$, the common ancestors of $Q$ is the set $\CA(Q)$ defined as:
\[
\CA(Q) = \{p\in D_{\T}\mid p\geq Q\}
\]
If $\T$ is regular then we have a well defined element $\min \CA(Q)$ which we call the least common ancestor $\LCA(Q)$.
\end{defi}

The definition is well posed since $\{p\in D_{\T}\mid p\geq Q\}$ is non empty if $\sup h(Q)<\infty$. Moreover it is bounded from below in terms of $h$. 
%Clearly, if $Q\subset Q' \subset D_{\T}$ then $\inf \CA(Q)\leq \inf \CA(Q')$. 
%Leveraging on the notion of least common ancestors we give the display poset of an abstract merge tree a pseudo-metric structure. , generalizing Definition $16$ of \cite{curry2021decorated} also to 

\begin{prop}\label{prop:metric}
The display poset $D_{\T}$ of any abstract merge tree can be given a pseudo-metric structure with the following formula:
\[
d((a,t),(b,t'))= (\tilde{t} - t)  + (\tilde{t} - t')  
\]
with $\tilde{t}= \inf \{h(p)\mid p \in \CA(\{(a,t),(b,t')\})\}$. If $\T$ is regular, then, $d$ is a metric.
\begin{proof}
    See \Cref{sec:proofs_TDA}.
\end{proof}
\end{prop}

\begin{figure}	
    \begin{subfigure}[c]{\textwidth}
    	\centering
    	\includegraphics[width = \textwidth]{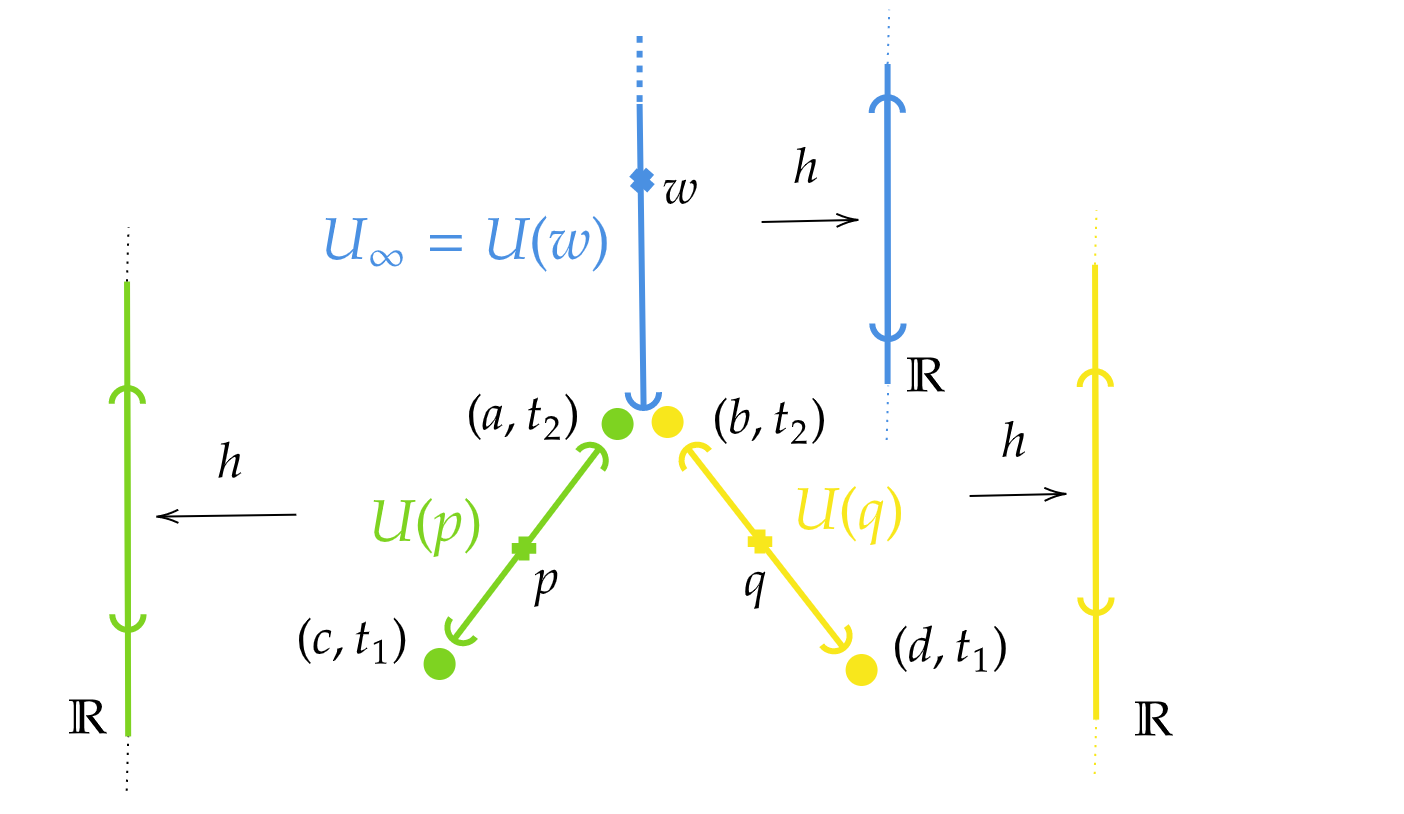}
	\end{subfigure}
%	\hspace{0.5 cm}

\caption{A graphical representation of the display poset, with its a.e. covering - see \Cref{sec:metric_display} - highlighted by the brackets at the endpoints  of the edges. Each such covering is the mapped homeomorphically to $\mathbb{R}$ via the height function $h$.
 Note that $d((a,t_2),(b,t_2))=0$ and $\{(a,t_2),(b,t_2)\}=CA((c,t_1),(d,t_1))$. The color scheme is coherent with the one in \Cref{fig:preliminary}.}
\label{fig:a.e. covering}
\end{figure}

We point out that a  similar definition has already been considered in Definition $16$ of \cite{curry2021decorated},  but stated only for regular merge trees, obtaining a (family of) metric(s) indexed on $p$. We exploit some results therein obtained to prove the following proposition, showing that the pseudo-metric we define is indeed a natural choice.

\begin{prop}\label{prop:just}
Given a compact smooth manifold $X$, if $\T$ is the merge tree associated to the sublevel set filtration of a Morse function $f:X\rightarrow \R$, then $D_{\T}$ is homeomorphic to the merge tree of $f$ as defined in \cite{merge_interl}.

\begin{proof}
The merge tree defined in \cite{merge_interl}, also referred to as classical merge tree in \cite{curry2021decorated}, is given by the Reeb graph of $\pi_f :E_f\rightarrow \R$ where: 
$E_f:=\{(x,r)\in X\times \R\mid f(x)\leq r\}$ is the epigraph of $f$ and
$\pi_f$ is the projection on the second component. The Reeb graph is then defined as $M_f:=E_f/\sim$ with $p=(x,r)\sim q=(x',r')$ iff $r=r'$ and $p,q\in A$, $A\in \pi_0(\pi_f^{-1}(r) )$.
The epigraph has the subset topology from the product topology and so $M_f$ inherits a quotient topology. 

First we easily see the bijection between the display poset and the classical merge tree. The key observation is that $x\in f^{-1}((-\infty,t])$ if and only if $f(x)\leq t$. That is, if and only if $(x,r)\in E_f$ for all $r\leq t$. 
Thus we have $\pi_f^{-1}(t)=\{(x,t)\mid f(x)\leq t\}$  and its image under $\pi_X:E_f\rightarrow X$ is just $f^{-1}((-\infty,t])$. 

Now, closed sets in $E_f$ are generated by sets of the form $C:=(C_X \times I)\cap E_f$, with $C_X$ closed in $X$, and $I$ closed in $\R$. Being $f$ continuous, $E_f$ is closed (see for instance \cite{rockafellar2009variational}, Ch. 1), and so $C$ is closed also in $X\times \R$. We know that $\pi_X:X\times\R\rightarrow X $ preserves closed sets, and thus $\pi_X(C)$ is closed in $X$. As a consequence, also $\pi_X:E_f\rightarrow X $ preserves closed sets.
Similarly, we have that the bijection $(A,r)\mapsto (\pi_X(A),r)$ between $\pi_f^{-1}(t)$  and $f^{-1}((-\infty,t])$ is an homeomorphism.  
Thus, $\pi_0(\pi_f^{-1}(t))\cong \pi_0(f^{-1}((-\infty,t]))$ via $A\mapsto \pi_X(A)$.
So, we also have a bijection between $M_f$ and $D_{\T}$. 
Which means that we can transfer on $M_f$ the topology defined (by the pseudo-metric) on $D_{\T}$.

Since $f$ is Morse, $\T$ is regular, and we can apply Proposition 2 of \cite{curry2021decorated} for $p=1$ which states that the two topologies on $M_f$ coincide. 
 
\end{proof}
\end{prop}

See \Cref{fig:a.e. covering} for an example of a display poset with its pseudo metric structure.

\begin{rmk}\label{rmk:metric_quotient}
\Cref{prop:metric} states that if $\T$ is a regular abstract merge tree, then we can induce a metric on $\mathcal{M}(\T)=(T,h_T)$ via $E_T\hookrightarrow D_{\T}$. It is not hard to see that this is the shortest path metric on $E_T$, with the length of an edge $e=(v,v')$ being given by $h_T(v')-h_T(v)$.
\end{rmk}

\begin{rmk}
Given $\T$ abstract merge tree, we have that the quotient of $D_{\T}$ under the relation $x\sim y$ iff $d(x,y)=0$, is isometric as a metric space (with the induced metric on the quotient) to
 $D_{R(\T)}$. 
\end{rmk}

\subsection{Functions Spaces on the Display Posets}

Thanks to \Cref{prop:metric} any display poset of an abstract merge tree can be given the topology generated by the open balls of the (pseudo) metric.

Consider now an abstract merge tree $\T$, with the critical points (see \Cref{defi:crit_pts}) of $D_{\T}$ being $\{v_1,\ldots,v_n\}$.
Without loss of generality, in this section, we suppose $h(v_i)\neq h(v_j)$ for $i\neq j$. 
So, we can set $t_i=h(v_i)$, and have
$t_i<t_{i+1}$.

Take $t\notin \{t_1,\ldots,t_n\}$. 
For any point $p=(a,t)\in D_{\T}$, we call:
\[
v_p = \max \{ v_i \in \{v_1,\ldots,v_n\} \mid v_i<p\}, \qquad t_p=h(v_p),
\]
\[
v^p = \min \{ v_i \in \{v_1,\ldots,v_n\} \mid v_i>p\}, \qquad t^p=h(v^p).
\]

An open ball of radius $\varepsilon>0$ is by definition: 
\[B_\varepsilon(p):=\{q \in D_{\T}\mid d(p,q)<\varepsilon\}.\]

Consider now $\varepsilon>0$, with $t_p\leq t-\varepsilon < t+\varepsilon \leq t^p$. Let $p=(a,t)$ be a point such that for every $\eta>0$ small enough: 
\[
\#\left( \pi_0(X_{t-\eta<t})^{-1}(a)\right)=1 \text{ and }\# \left(\pi_0(X_{t<t+\eta})^{-1}(\pi_0(X_{t<t+\eta})(a))\right)=1.
\]

The ball of radius $\varepsilon$ around $p$ is: 
\[
B_\varepsilon(p):=\{q \in \CA(\{p\})\mid h(q)<t+\varepsilon\} \cup
\{q \mid p \in \CA(\{q\})\text{ and }h(q)>t-\varepsilon\}.\]
Thus, for any such point $p=(a,t)$
we can define the set:
\begin{align*}
   U(p):=&\{q \in \CA(\{p\})\mid h(q)<t^p\} \cup \{q \mid p \in \CA(\{q\})\text{ and }h(q)>t_p\} \\
   =&\{q\in D_{\T}\mid  q < (v^p,t^p) \} \cap \{q\in D_{\T}\mid  q > (v_p,t_p) \} .
\end{align*}

which is an open neighbor of $p$. If 
$t> t_n$, then $t^p=\infty$ and so we have:

\[
U_{\infty}:=U(p)=\CA(\{(v_n,t_n))\}).
\]

Refer to \Cref{fig:a.e. covering} to have a visual intuition for the following proposition.

\begin{prop}\label{prop:h_omeo}
The map $h:D_{\T}\rightarrow \mathbb{R}$ is monotone, continuous and
 $h_{\mid U(p)}:U(p)\rightarrow(t_p,t^p)$ is an homeomorphism and an isometry.    

\begin{proof}
    Using the same notation of \Cref{prop:metric}, we have:
    \[
\mid h((a,t))-h((b,t´))\mid = \mid t-t'\mid \leq \tilde{t} - t  + \tilde{t} - t'  = d((a,t),(b,t')).  
    \]
    Thus $h$ is continuous. Monotonicity is trivial. Suppose now we have $p=(a,t)$ and $(b,t'),(c,t')\in U(p)$ such that $b\neq c$. This implies that, for every $\varepsilon>0$ small enough we have:
\[
\#\left(\pi_0(X_{t<t+\varepsilon})^{-1}(\pi_0(X_{t<t+\varepsilon})(b))\right)>1.
\]    
By construction, this means that, for every $\varepsilon>0$, there must be $v_i$ such that $(b,t'),(c,t') <(v_i,t_i) < q$, with $q=(\pi_0(X_{t<t+\varepsilon})(b),h(\pi_0(X_{t<t+\varepsilon})(b)))$. In particular, $(v_i,t_i) < (v^{(b,t')},t^{(b,t')})$.
Which is absurd.
Moreover, $h_{\mid U(p)}$ is clearly surjective for $h(\pi_0(X_{t<t'})(a))=t'$. Thus $h_{\mid U(p)}$ is a bijective map. If $(b,t'),(c,t'')\in U(p)$, then:
\[
\tilde{t}= \inf \{h(q)\mid q \in \CA(\{(b,t'),(c,t'')\})\}= \min \{t',t''\},
\]
which implies that 
    $h_{\mid U(p)}$ is an isometry. And thus an homeomorphism. 
\end{proof}
\end{prop}

\begin{defi}
The set $\mathcal{U}(D_{\T}):=\{U\subset D_{\T}\mid U = U(p)\text{ for some }p\in D_{\T}\}$ is called the a.e. canonical covering of $D_{\T}$. 
\end{defi}

\begin{rmk}
    Recall that the sets $U(p)$ are defined only for points $p=(a,t)$ for which there is $K>0$ such that for every $0<\varepsilon<K$, we have $\#\left(\pi_0(X_{t-\varepsilon<t})^{-1}(a)\right)=1$ and $\#\left(\pi_0(X_{t<t+\varepsilon})^{-1}(\pi_0(X_{t<t+\varepsilon})(a))\right)=1$.
\end{rmk}

\begin{prop}
    The set $\mathcal{U}(D_{\T})$ is finite by the finiteness of $\T$. Moreover, for every $q\in U(p)$, $U(p)=U(q)$.
\begin{proof}
    We prove only the second part of the statement.

    Consider $q\in U_p$.  
    Since $q\in U_p$, we have either $(v_p,t_p)<p<q<(v_p,t_p)$ or $(v_p,t_p)<q<p<(v_p,t_p)$. 
    Since $\{v \in D_{\T} \mid v>q\}$ is a totally ordered set, we must have $(v^p,t^p)\leq(v^q,t^q)$ or $(v^q,t^q)\leq (v^p,t^p)$.
    Which implies $v^p=v^q$. 

    Lastly, suppose 
     $v_p \neq v_q$. Having  $(v_p,t_p)\leq(v_q,t_q)$ or $(v_q,t_q)\leq (v_p,t_p)$ would imply $v_p=v_q$. So, we have $(v_p,t_p)\nleq(v_q,t_q)$ and $(v_q,t_q)\nleq (v_p,t_p)$. But this, in turn, implies the existence of $(v_i,t_i) > (v_p,t_p),(v_q,l_q)$ such that $(v_i,t_i) < p,q$, which is absurd.
\end{proof}
\end{prop}

With the help of $\mathcal{U}(D_{\T})$ we want
to induce a measure on the sigma algebra generated by the open sets of $D_{\T}$. 
This measure is inspired by the fact that Reeb graphs (and so merge trees) are stratified covering of the real line (see \cite{de2016categorified} and references therein) and thus we want to locally pull back a measure from the real line.

For a display poset $D_{\T}$ we define the measure $\mu_{\T}$ as:
\[
\mu_{\T}(Q)=\sum_{U\in \mathcal{U}(D_{\T})} \mathcal{L}(h(U\cap Q))
\]
A graphical representation of such measure can be found in \Cref{fig:display_measure}. Note that, if we call $D_{\T}^\circ = \bigcup_{U\in \mathcal{U}(D_{\T})}U$, we have $\mu_{\T}(D_{\T}-D_{\T}^\circ)=0$.

\begin{prop}
The expression:
\[
\mu_{\T}(Q)=\sum_{U\in \mathcal{U}(D_{\T})} \mathcal{L}(h(U\cap Q)),
\]
induces a measure on the sigma algebra generated by the open sets of $D_{\T}$.
\begin{proof}

We prove that $\mu_{\T}$ is $\sigma$-additive.
Let $X_i$, $i\in\mathbb{N}$, be disjoint sets in the Borel sigma algebra of $D_{\T}$; we need to prove that $\mu_{\T}(\bigcup_{i\in\mathbb{N}}X_i)=\sum_{i\in\mathbb{N}}\mu_{\T}(X_i)$.

We have: 
\[
(\bigcup_{i\in\mathbb{N}} X_i)\cap U = \{p \in D_{\T}\mid p \in X_i \text{ for some }i\text{ and }p \in U\} = \bigcup_{i\in\mathbb{N}} (X_i\cap U)
\]
 and so we are finished since $\mathcal{L}$ is $\sigma$-additive on $h(U\cap X_i)$.
Note that, if $Q$ is in the Borel sigma algebra of $D_{\T}$, being $h$ an homeomorphism on $U$ (due to \Cref{prop:h_omeo}),
$h(U\cap Q)$ is always Lebesgue measurable in $\mathbb{R}$.

%All other properties are verified analogously, exploiting the properties of the intersection and $\mathcal{L}$. 
%\begin{itemize}
%\item non-negativity: it is clearly satisfied;
%\item null empty set: $\emptyset \bigcap U=\emptyset$ and thus is it satisfied;
%\item countable additivity: let $A_i$, $i\in\mathbb{N}$ open sets in $D_{\T}$; then $(\bigcup_{i\in\mathbb{N}} A_i)\bigcap U = \bigcup_{i\in\mathbb{N}} (A_i\bigcap U)$ and so we are finished since $U\bigcap A_i$ is homeomorphic to an open set of $\mathbb{R}$ and countable additivity holds for the Lebesgue measure.

\end{proof}
\end{prop}

Note that the quotient map $D_{\T}\rightarrow D_{R(\T)}$ mentioned in \Cref{rmk:metric_quotient}, preserves measures.

Now, consider a function $f:D_{\T}\rightarrow \mathbb{R}$: by definition we have that $f$ is $\mu_{\T}$-\emph{measurable} if $f \circ (h_{\mid U})^{-1}$ is $\mathcal{L}$-measurable on $\mathbb{R}$ for every $U\in\mathcal{U}(D_{\T})$.
So, given a $\mu_{\T}$-measurable function $f:D_{\T}\rightarrow \mathbb{R}$ we can define:
\[
\int_{D_{\T}}f d\mu_{\T}= \sum_{U(p)\in \mathcal{U}(D_{\T})} \int_{t_p}^{t^p} f \circ (h_{\mid U(p)})^{-1} d\mathcal{L}.
\]

Leveraging on this definition, we want to define a framework to work with functions defined in some metric space $(E,d_e)$. For reasons which will be clarified in the next section, we want that inside the metric space $E$ there is a reference element $0$ such that the amount of information contained in the value $f(p)$ can in some sense be quantified as the distance $d_e(f(p),0)$. So we
 make the following assumption.
 
\begin{assump}
We always assume that $(E,d_e)$ is a metric space and that $(E,\ast, 0)$ is a monoid, i.e. that $\ast$ is an associative operation with neutral element $0$.
\end{assump}

We establish the following notation for any measure space $(M,\mu)$:
\[
L_p(M,E):=\Big \lbrace f:M\rightarrow E\mid d(f(\cdot),0):M\rightarrow \mathbb{R} \text{ measurable and }\int_{M}d_e(f(\cdot),0)^p d\mu<\infty\Big \rbrace /\sim
\]
with $\sim$ being the usual equivalence relation between functions identifying functions up to $\mu$-zero measure sets. This space becomes a monoid and a metric space with $(f+g)(p):=f(p)\ast g(p)$ and:
\[
d_{L_p}(f,g)=\int_{M}d_e(f(\cdot),g(\cdot))^p d\mu.
\]
To verify that $d_{L_p}$ is a metric is enough to see that $d_{L_p}(f,g)=0$ if and only if $f$ and $g$ differ on $\mu$-zero measure sets and prove the triangle inequality using that $L_p(M,\mathbb{R})$ is a normed space.

For the sake of brevity, in the following we do not write explicitly the request that $d(f(\cdot),0)$ is measurable and we imply it in the existence of its integral.
Thus, we are interested in the spaces:
\[
L_p(D_{\T},E):=\Big \lbrace f:D_{\T}\rightarrow E\mid \int_{D_{\T}}d_e(f(\cdot),0)^p d\mu_{\T}<\infty\Big \rbrace /\sim
\]

Consider now $\T$ and $\G$ such that $\T \cong_{a.e.} \G$. Let $Z\subset \mathbb{R}$ such that $\alpha:\pi_0(X_{\Bigcdot \mid Z})\rightarrow \pi_0(Y_{\Bigcdot \mid Z})$ is a natural isomorphism and $\mathcal{L}(\mathbb{R}-Z)=0$. Then $\alpha$ induces a bijection between the display posets:
\[
D_{\pi_0(X_{\Bigcdot \mid Z})}:=\bigcup_{t\in Z}\pi_0(X_t)\times\{t\} 
\]
and 
\[
D_{\pi_0(Y_{\Bigcdot \mid Z})}:=\bigcup_{t\in Z}\pi_0(Y_t)\times\{t\}. 
\]
With an abuse of notation we call such bijection $\alpha:D_{\pi_0(X_{\Bigcdot \mid Z})}\rightarrow D_{\pi_0(Y_{\Bigcdot \mid Z})}$.

Given $f:D_{\G}\rightarrow E$ we can clearly restrict it to $D_{\pi_0(Y_{\Bigcdot \mid Z})}$ and thus we can pull it back on $D_{\pi_0(X_{\Bigcdot \mid Z})}$ with $\alpha$:
\[
D_{\pi_0(X_{\Bigcdot \mid Z})}\xrightarrow{\alpha}
D_{\pi_0(Y_{\Bigcdot \mid Z})}\hookrightarrow D_{\G}\xrightarrow{f}E
\]
We call such function $\alpha^*f$. 

\begin{prop}\label{prop:fun_iso}
The rule $f\mapsto \alpha^*f$ described above induces map $\alpha^*:L_p(D_{\G},E)\rightarrow L_p(D_{\T},E)$ which is an isometry and a map of monoids.
\begin{proof}
Since 
$\mathcal{L}(\mathbb{R}-Z)=0$ then both $f\in L_p(D_{\pi_0(Y_{\Bigcdot \mid Z})},E)$ and $\alpha^*f\in L_p(D_{\pi_0(X_{\Bigcdot \mid Z})},E)$ identify a unique equivalence class, respectively,
in $L_p(D_{\G},E)$ and $L_p(D_{\T},E)$.
Moreover, it is easy to see that the map $\alpha^*$ is such that $\alpha^*(f+g)=\alpha^*f + \alpha^*g$ and $d_{L_p}(f,g)=d_{L_p}(\alpha^*f,\alpha^*g)$.
Lastly, because $\alpha$ is a natural isomorphism, then $(\alpha^{-1})^*$ yields the opposite correspondence.
\end{proof}
\end{prop}

\Cref{prop:fun_iso} implies that, for our purposes, we can always restrict ourselves to considering regular abstract merge trees. Thus we make the following assumption.

\begin{assump}
From now on we will always suppose that any abstract merge tree we consider is regular.
\end{assump}

\subsection{Local Representations of Functions}
\label{sec:local_rep}

When comparing two functions $f$, 
$g$ defined on different display posets, we face the problem of combining together two kinds of variability: using language borrowed from functional data analysis (see the Special Section on Time Warpings and Phase Variation on the Electronic Journal of Statistics, Vol 8 (2), and references therein) and shape analysis \cite{kendall_1977, kendall_1984, stat_shape} we have an \virgolette{horizontal} variability, due to the different domains (i.e. display posets), and a \virgolette{vertical} variability which depends on the actual values that the functions assume. It is reasonable that both kinds of variability contribute to the final distance value: we have a cost given by aligning the two display posets - horizontal variability - and a cost arising from the different amplitudes of the functions - vertical variability. In particular, we would like the horizontal variability to be measured in a way which is suitable for abstract merge trees (for instance, it should posses some kind of stability properties) and, similarly, the way in which the amplitude variability is measured should take a natural form, related to the spaces $L_p(D_{\T},E)$.

%In scientific literature, alignment (of shapes, functions, etc.) is carried out as a well-behaved transformation of the initial datum (or domain in the case of functions) in order to minimize differences (or maximize similarities) between different objects.
%Once the alignment is completed, the residual variability is measured and used for the analysis.  
%Usually data are aligned in a pairwise fashion
%to compute distances or to a reference object - typically some notion of average. 
%Setting up a coherent framework in which functions on different merge trees can be \virgolette{continuously} deformed one in the other presents at least two kinds of problems: 1) carrying out deformations (and so optimization) in the space of display posets 2) measuring the amplitude variability taking into account the measure-related changes in the domain due to deformations. Both, if approached in this way, are challenges are daunting problems, as alignment of trees is already very complicated when trees are embedded into $\mathbb{R}^3$ \cite{sriva_trees}.  
 
In other words, given $f:D_{\T}\rightarrow E$ and
$g:D_{\G}\rightarrow E$ we want to align, deform the display posets by locally comparing the information given by $f$ and $g$ and matching the display posets in a convenient way. The word \emph{locally} is on purpose vague at this stage of the discussion and should be thought as in some neighborhood of points of the posets. To compare local information carried by functions, we need to embed such objects in a common space so that differences can be measured.

\begin{figure}
    \begin{subfigure}[c]{0.49\textwidth}
    	\centering
    	\includegraphics[width =\textwidth]{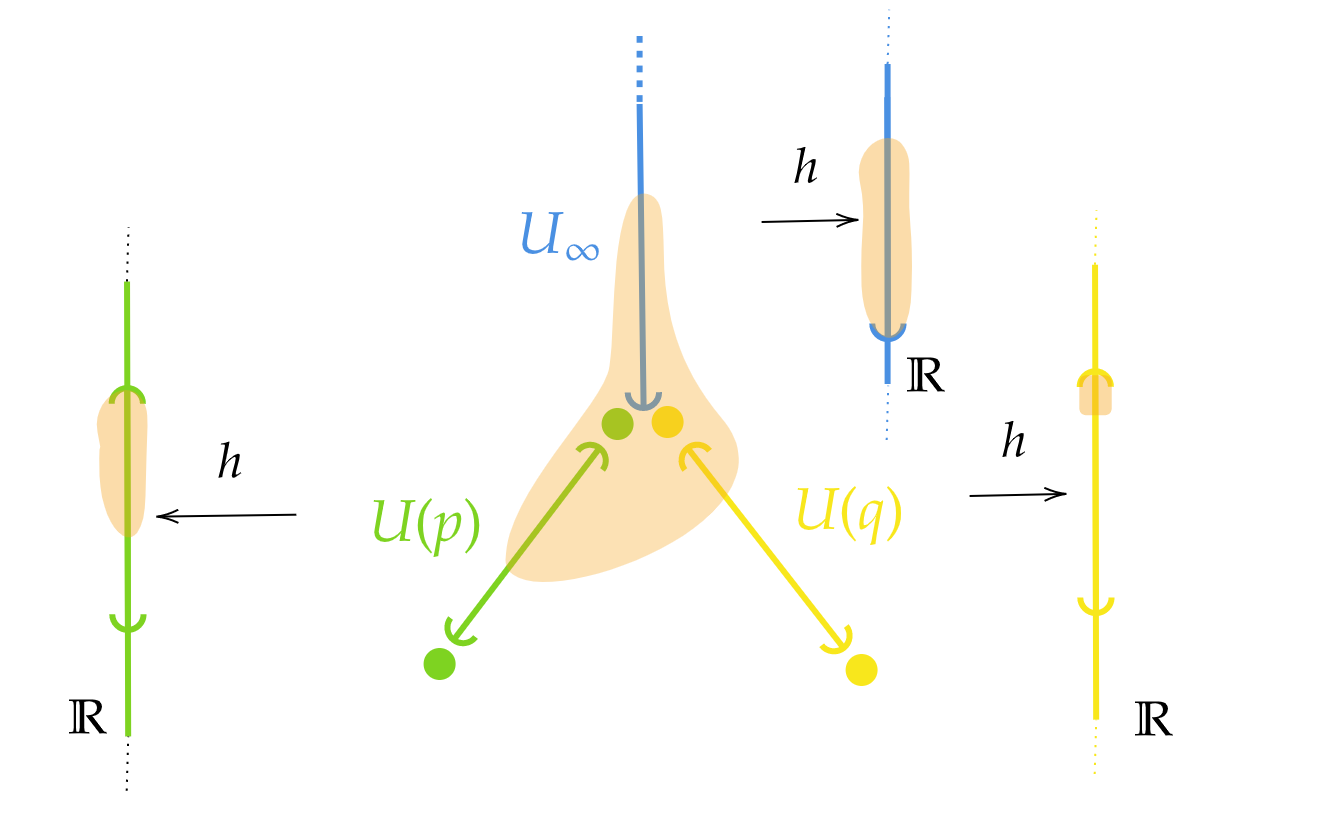}
            \captionsetup{singlelinecheck=off, margin={0.3cm, 0.1cm}}
    	\caption{A display poset $D_{\T}$ with the measure $\mu_{D_{\T}}$. The orange shaded set is first intersected with the open sets $U_\infty$, $U(p)$ and $U(q)$ and then its Lebesgue measure is taken in $\mathbb{R}$ via the height function $h$.}
    	\label{fig:display_measure}
	\end{subfigure}
    \begin{subfigure}[c]{0.49\textwidth}
    	\centering
    	\includegraphics[width =\textwidth]{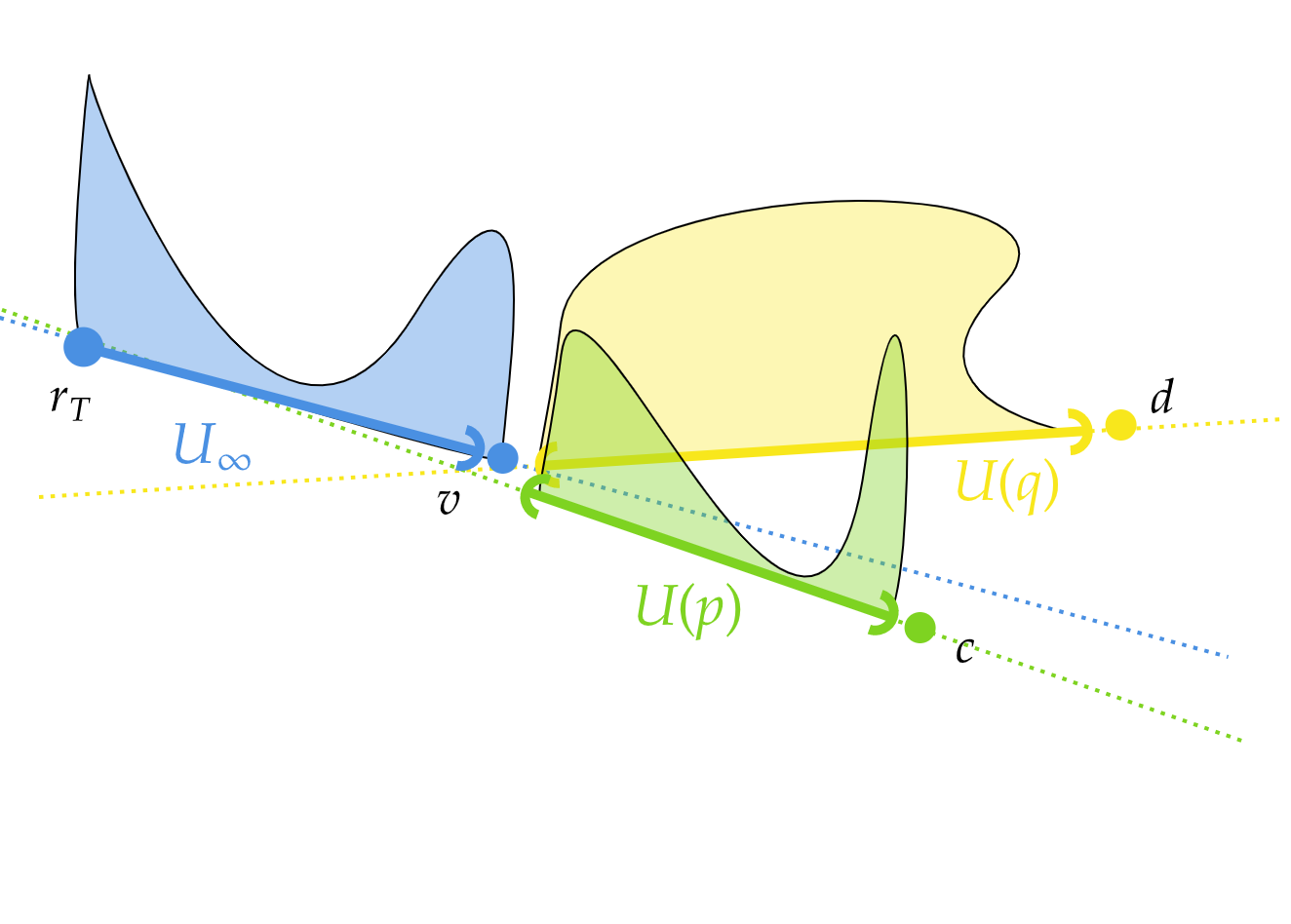}
            \captionsetup{singlelinecheck=off, margin={0.3cm, 0.1cm}}
    	\caption{A function $f:D_{\T}\rightarrow \mathbb{R}$ defined on the display poset $D_{\T}$. With different colors we have highlighted the restrictions of the function on the different open sets of the canonical a.e. covering.}
    	\label{fig:MT_f}
	\end{subfigure}

    \begin{subfigure}[c]{0.49\textwidth}
    	\centering
	    \includegraphics[width = \textwidth]{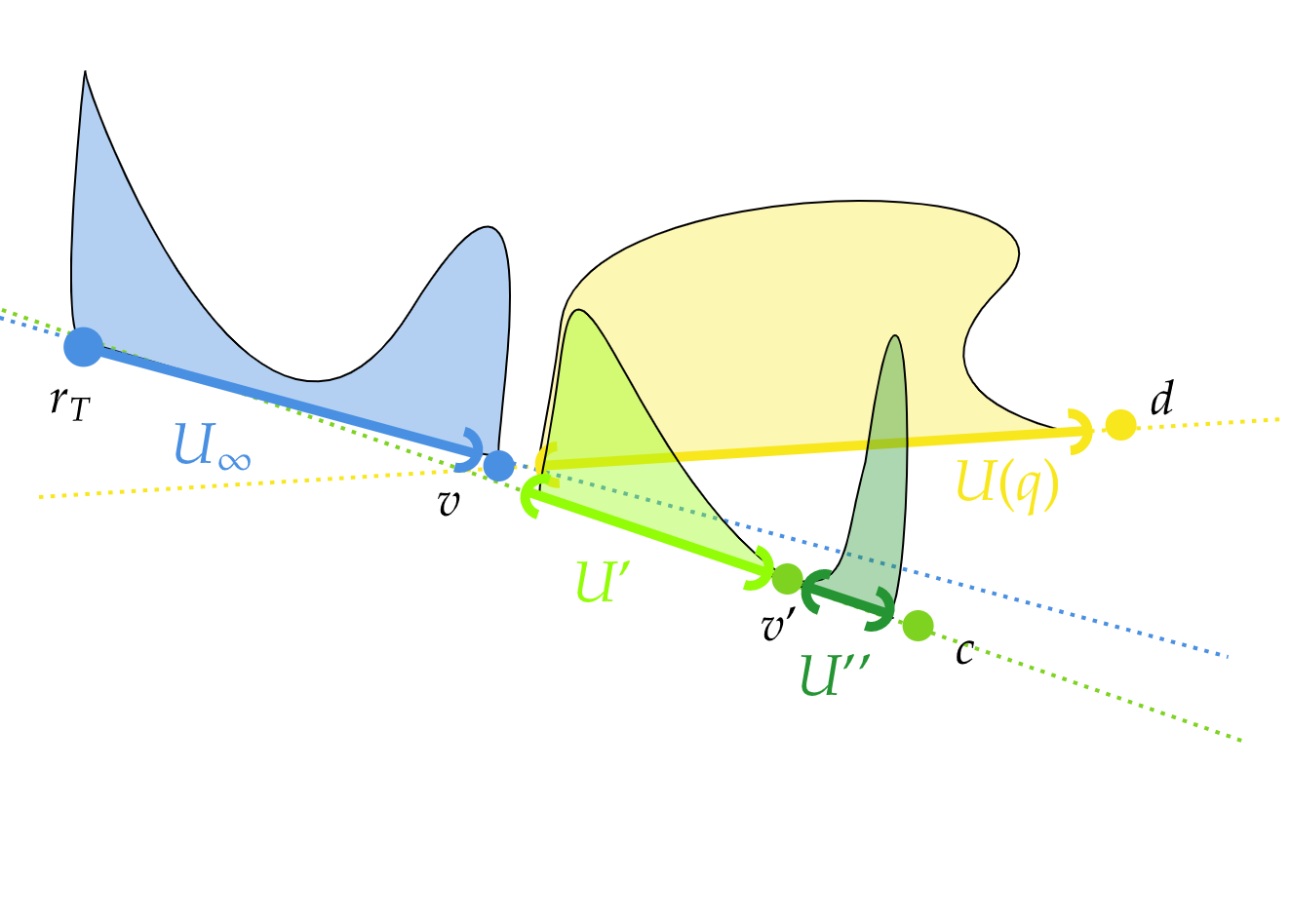}
            \captionsetup{singlelinecheck=off, margin={0.3cm, 0.1cm}}
    	\caption{A function $f:D_{\T}\rightarrow \mathbb{R}$ defined on the display poset $D_{\T}$  represented with the restrictions on a regular a.e. covering which refines the canonical one.}
    	\label{fig:f_ref}
	\end{subfigure}
%	\hspace{0.5 cm}
    \begin{subfigure}[c]{0.49\textwidth}
    	\centering
	    \includegraphics[width = \textwidth]{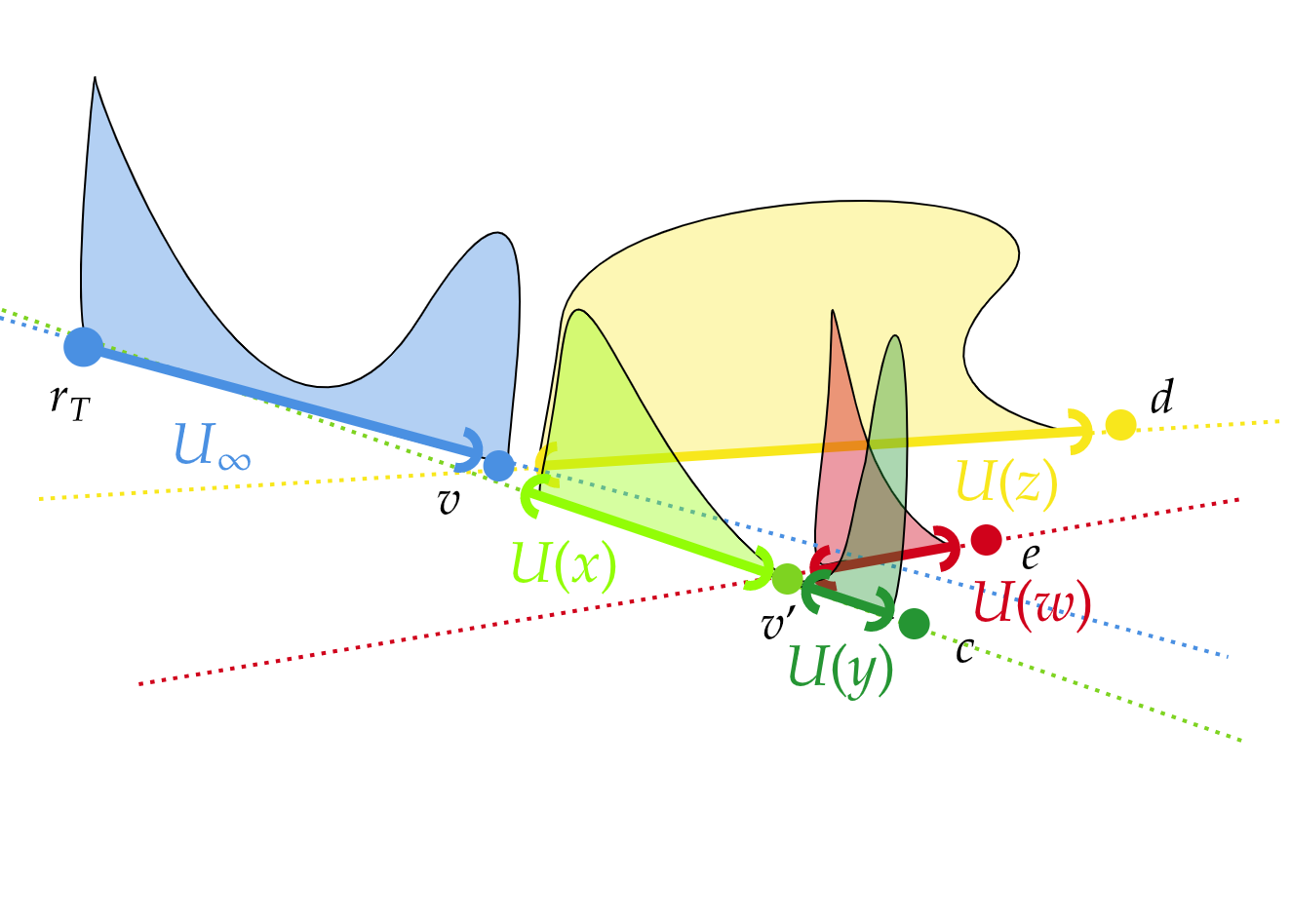}
            \captionsetup{singlelinecheck=off, margin={0.3cm, 0.1cm}}
    	\caption{A function $g:D_{\G}\rightarrow \mathbb{R}$ defined on the display poset $D_{\G}$ along with its restrictions on the canonical a.e. covering of $D_{\G}$. The refinement of the canonical a.e. covering of $D_{\T}$ which appears in \Cref{fig:f_ref} is much more suited than the canonical a.e. covering in \Cref{fig:MT_f} to compare the two functions: on $U_\infty$ and $U_\infty$ the functions are very similar, as the are on 
$U'$ and $U(x)$, on 
$U''$ and $U(y)$ and on 
$U(q)$ and $U(z)$.}
    	\label{fig:MT_g}
	\end{subfigure}

\caption{Measures and real valued functions defined on display posets. In every plot but the upper left one, for visualization purposes the posets are represented as embedded on the horizontal plane in $\mathbb{R}^3$ and plotted with thick lines. The vertical axis represents the value of the functions. With different colors we have highlighted the restrictions of the functions on different open sets. The colored dotted lines are are a qualitative visual representation of the embedding $\left( f:(a,b)\rightarrow \mathbb{R}\right) \mapsto \left(f':\mathbb{R}\rightarrow \mathbb{R}\right)$ where $f'$ extends $f$ with $0$ outside $(a,b)$.}

\label{fig:refinement}
\end{figure}

First we formalize the procedure of obtaining local information from a function $f:D_{\T}\rightarrow E$ - \Cref{fig:MT_f} can help in the visualization of such idea.
Given $D_{\T}$ display poset and its a.e. canonical covering, we have an isomorphism of metric spaces and monoids:
\[
L_p(D_{\T},E) \cong \bigoplus^p_{U\in \mathcal{U}(D_{\T})} L_p(h(U),E)
\]
where $\bigoplus^p$ means that the norm of the direct sum is the $p$-th root of the sum of the $p$-th powers of the elements in the direct sum.

In this way we split up a function $f$ on open disjoint subsets, without losing any information. However, as in \Cref{fig:refinement}, to compare different functions one may need to represent this information on a finer scale and thus
$\mathcal{U}_{D_{\T}}$ may not be the correct way to split up $f$, which may need to be partitioned in smaller pieces. Thus we allow $\mathcal{U}_{D_{\T}}$ to be refined with particular collections of open sets.

\begin{defi}
A collection of open sets of $D_{\T}$ is an a.e. covering of $D_{\T}$ if it covers $D_{\T}$ up to $\mu_{\T}$-zero measure set. An a.e. covering of $D_{\T}$ is regular
if it is made by disjoint, path-connected open sets, each contained in some $U\in \mathcal{U}(D_{\T})$.
   
Given $\mathcal{O}'$ regular a.e. covering of $D_{\T}$, a refinement of $\mathcal{O}'$ is a regular a.e. covering $\mathcal{O}$ such that
 for every $U\in\mathcal{O}$ there is $U'\in \mathcal{O}'$ such that $U\subset U'$.
\end{defi}

Given the display poset $D_{\T}$ of an abstract merge tree $\T$ we collect all the refinements of its a.e. canonical covering in the set $\Cov(\T)$. Note that, by definition, these are all regular coverings.

\begin{prop}\label{prop:lattice}
The set $\Cov(\T)$ is a lattice. It is a poset
with the relation $\mathcal{O}<\mathcal{O}'$ if $\mathcal{O}$ is a refinement of $\mathcal{O}'$ and for every pair of elements $\mathcal{O}$, $\mathcal{O}'$ there is a unique least upper bound $\mathcal{O}\vee \mathcal{O}'$ and a unique greater lower bound $\mathcal{O}\wedge \mathcal{O}'$. The operations are defined as follows:
\[
\mathcal{O}\vee \mathcal{O}':= \pi_0\left( \bigcup_{U\in\mathcal{O}' \text{ or }U\in\mathcal{O}} U\right)
\]

\[
\mathcal{O}\wedge \mathcal{O}':= \{U\cap U' \mid U'\in\mathcal{O}' \text{ and }U\in\mathcal{O}\}. 
\]
\begin{proof}
    See \Cref{sec:proofs_TDA}.
\end{proof}
\end{prop}

Given $\mathcal{O}\in \Cov(\T)$ we have:
\[
L_p(\T,E) \cong \bigoplus^p_{U\in \mathcal{O}} L_p(h(U),E)
\]

As already mentioned, to compare functions defined on different abstract merge trees we want to embed all these representations of functions into one common metric space, shared by all abstract merge trees. What we do is to consider $L_p((a,b),E)$, for some $(a,b)\subset\mathbb{R}$ and embed it into $L_p(\mathbb{R},E)$ by extending $f:(a,b)\rightarrow E$ to $\mathbb{R}$ with $0\in E$ outside $(a,b)$. In this way we have an isometric embedding $L_p((a,b),E)\hookrightarrow L_p(\mathbb{R},E)$.

In the next definition we need the notion of the essential support of a function $f:(M,\mu)\rightarrow E$ defined on a measure topological space $(M,\mu)$ and with values in $(E,\ast,0)$: 
\[
\supp(f)= M-\bigcup \{U\subset M \text{ open }\mid f_{\mid U}=0\text{ }\mu-\text{a.e.}\} 
\]

\begin{defi}\label{defi:local_repr}
Given $D_{\T}$ and $\mathcal{O}\in \Cov(\T)$, a local representation of a function in $L_p(D_{\T},E)$ on $\mathcal{O}$ is a function $\varphi_{\mathcal{O}}:\mathcal{O}\rightarrow L_p(\mathbb{R},E)$
such that $\supp(\varphi_{\mathcal{O}}(U))\subset h(U)$ for every $U\in\mathcal{O}$. 
\end{defi}

Note that if, instead of splitting $f$ on a finer scale, we want to look at the function on a coarser level, we can do that.
Consider $\mathcal{O}'$ refinement of $\mathcal{O}$; then for every $V\in\mathcal{O}$:
\[
\varphi_{\mathcal{O}}(V)=\sum_{U\in \mathcal{O}'\text{ such that } U\subset V} \varphi_{\mathcal{O}'}(U)
\]

\subsection{Regular Coverings and Merge Trees Up to Degree $2$ Vertices}
\label{sec:fun_on_MT}

Thanks to \Cref{prop:fun_iso} we have seen that to work with functions defined on display posets we can reduce to the case of regular abstract merge trees.
This makes the upcoming discussion much easier since, thanks to \Cref{prop:equivalence}, we can associate a merge tree to any regular
abstract merge tree. In particular, in this section we deal with the problem of associating functional weights to the edges of a merge tree, so that this becomes a combinatorial representation of a function defined on a display poset.

We have already seen that the metric $d$ defined on the display poset $D_{\T}$ induces the shortest path metric on the graph $(T,h_T)=\mathcal{M}(\T)$ via the inclusion $E_T\hookrightarrow D_{\T}$ - see \Cref{prop:equivalence}. 
Similarly, we can establish a correspondence between the edges $E_T$ and the
 a.e. canonical covering $\mathcal{U}(D_{\T})$: each edge $(v,v')\in E_T$
corresponds to the open set $U=\{p\in D_{\T}\mid v<p<v'\}$ or $U_{\infty}=\{p\in D_{\T}\mid v<p\}$ if $v'=r_T$ - as in \Cref{fig:MT}. This correspondence can be extended to a bijection between the equivalence class of merge trees up to degree $2$ vertices.

\begin{prop}\label{prop:regular_and_order_2}
Consider $T=\mathcal{M}(\T)$ and call $[T]$ the equivalence class of $T$ up to degree $2$ vertices. Then the set $\Cov(\T)$ and $[T]$ are in bijection $\mathcal{O}\mapsto T_{\mathcal{O}}$, with $T$, the only merge tree in $[T]$ without degree $2$ vertices, being mapped to the
 a.e. canonical covering $\mathcal{U}(D_{\T})$.
Moreover $\mathcal{O}<\mathcal{O'}$ if and only if $T_\mathcal{O}$ can be obtained from $T_\mathcal{O'}$ via ghostings.
\begin{proof}
The map is induced by each edge $(v,v')\in E_T$
being sent into the open set $U=\{p\in D_{\T}\mid v<p<v'\}$ or $U_{\infty}=\{p\in D_{\T}\mid v<p\}$ if $v'=r_T$.
The result then follows from \Cref{prop:h_omeo} plus the fact that path-connected subsets of $\R$ are connected intervals. 
\end{proof}
\end{prop}

As a consequence, we also have the following corollary, finally bridging between functions defined on display posets and weighted trees.

\begin{cor}\label{cor:local_repr}
Given an abstract merge tree $\T$ and the merge tree $T=\mathcal{M}(\T)$, we have a bijection between the following sets:
\[
\lbrace\varphi_{\mathcal{O}}:\mathcal{O} \rightarrow L_p(\R,E) \mid \mathcal{O}\in \Cov(\T) \text{ and } \supp(\varphi_{\mathcal{O}}(U))\subset h(U), \forall U \in \mathcal{O}\rbrace
\]
and
\[  \lbrace \varphi_{T'}:E_{T'} \rightarrow L_p(\R,E) \mid T'\in [T] \text{ and } \supp(\varphi_{T'}((v,v')))\subset [h_{T'}(v),h_{T'}(v')], \forall (v,v')\in E_{T'}\rbrace.
\]

\end{cor}

To sum up, we have proven that the local representation of a function on the display poset of an abstract merge tree is equivalent to a weighted tree, equal up to degree $2$ vertices to the merge tree representing the regular abstract merge tree, with the weights being the restriction of the original function to a suitable open set.  
For notational convenience, from now on, we may confuse the two sets in \Cref{cor:local_repr}, calling local representation of function also 
$\varphi_{T'}:E_{T'} \rightarrow L_p(\R,E)$ satisfying the requested properties.

\section{Edit Distance Between Local Representation of Functions}
\label{sec:edit_functions}

At this point we face the problem of defining a suitable (pseudo) metric framework for objects of the form $f\in L_p(D_{\T},E)$
and $g\in L_p(D_{\G},E)$, knowing that each of such objects can be represented by some $(T,h_T,\varphi_T)$, with $\varphi_T:E_T\rightarrow L_p(\mathbb{R},E)$ such that for each edge $e=(a,b)$: $\supp(\varphi_T(e))\subset[h_T(a),h_T(b)]$.

\subsection{Editing Local Representations of Functions}
\label{sec:edits_local_rep}

We want to exploit the edit distance defined \cite{pegoraro2023edit}. Such distance is inspired by tree edit distances \cite{gao2010survey, Tai}, but with key differences in the edit operations.
The philosophy of edit distances is to allow certain modifications of the base object, called edits, each being associated to a cost, and to define the distance between two objects as the minimal cost that is needed to transform the first object into the second with a finite sequence of edits.
In this way, up to properly setting up a set of edits, one can formalize the deformation of a tree comparing the local information induced by the weights of the trees, which, in our case, are the restrictions on the edges of a function defined on the display poset.  

The framework developed in \cite{pegoraro2023edit} requires that codomain of $\varphi_T:E_T\rightarrow W$ must satisfy certain properties.

\begin{defi}
A set $W$ is called editable if the following conditions are satisfied: 
\begin{itemize}
\item[(P1)] $(W,d)$ is a metric space
\item[(P2)] $(W,\ast,0)$ is a monoid (that is $W$ has an associative operation $\ast$ with zero element $0$)
\item[(P3)] the map $d(0,\cdot):W\rightarrow \mathbb{R}$
is a map of monoids between $(W,\ast)$ and $(\mathbb{R},+)$:
$d(0,x\ast y)= d(0,x)+d(0,y)$.
\item[(P4)] $d$ is $\ast$ invariant, that is: $d(x,y)=d(z\ast x,z\ast y)=d(x\ast z,y\ast z)$
\end{itemize}
\end{defi}

In \cite{pegoraro2023edit} it is shown that if $E$ is an editable space, then, also $W=L_1(\R,E)$ is an editable space. So, local representations of functions defined on a display poset fit into this framework as long as we take $p=1$ and $f:D_{\T}\rightarrow E$ takes values in an editable space. Moreover all the sets $\R_{\geq 0}$, $\mathbb{N}_{\geq 0}$ and their finite sums are editable spaces.

There are however situations which we want to avoid because they represent \virgolette{degenerate} functions which introduce formal complications.

\begin{defi}\label{defi:proper_info}
Given an editable space $E$ and a tree-structure $T$, a weight function $\varphi_T:E_T\rightarrow L_1(\R,E)$ is proper if we have $0\in \varphi(E_T)$ if and only if $E_T=\emptyset$ and $V_T=\{\star\}$. Analogously to \cite{pegoraro2023edit}, for the sake of brevity we call dendrogram the datum of a merge tree with a proper weight function $\varphi_T:E_T\rightarrow L_1(\R,E)$. 
\end{defi}

%\begin{assump}
%From now on we only work with proper local representations of functions.  To lighten the notation, however, we omit to write \virgolette{proper} explicitly. Similarly we omit repeating that all function spaces $L_1(\R,E)$ we consider are obtained with $E$ editable.
%\end{assump}

\begin{defi}
Given an (editable) space $L_1(\R,E)$ the dendrogram space $(\mathcal{T}, L_1(\R,E))$ is given by the set of dendrograms $(T,\varphi_T)$.
\end{defi}

Given an editable dendrogram space $(\mathcal{T}, L_1(\R,E))$, we can define our edits. 

\begin{itemize}

\item We call \emph{shrinking} of an edge 
a change of the local representation of a function associated to the edge. The new local representation function must be equal to the previous one on all edges, apart from the \virgolette{shrunk} one. In other words,  for an edge $e$, this means changing the value $\varphi(e)$ with another non zero function in $L_1(\R,E)$. Note that, in general, 
shrinkings do not preserve local representations of functions.
 
\item A \emph{deletion} is an edit with which an edge is deleted from the dendrogram. Consider an edge $(v_1,v_2)$. The result of deleting $v_1$ is a new tree structure, with the same vertices and edges a part from $v_1$ (the smaller one) and $(v_1,v_2)$, and with the parent of the deleted vertex which gains all of its children. Note that, if we start from a local representation of a function, the result of a deletion is always a local representation of a function.
The inverse of the deletion is the \emph{insertion} of an edge along with its lower vertex. 
We can insert an edge at a vertex \(v\) specifying the name of the new child of \(v\), the children of the newly added vertex (that can be either none, or any portion of the children of \(v\)), and the value of the function on the new edge. Again, insertions do not preserve local representations of functions.

\item Lastly, we generalize \Cref{defi:ghosting}, defining a transformation which eliminates a degree two vertex in a dendrogram, changing the local representation of a function. 
Suppose we have two edges $e=(v_1,v_2)$ and $e'=(v_2,v_3)$, with $v_1<v_2<v_3$. And suppose $v_2$ is of degree two. Then, we can remove $v_2$ and merge $e$ and $e'$ into a new edge $e''=(v_1,v_3)$, with $\varphi(e''):= \varphi(e)+ \varphi(e')$.
This transformation is called the \emph{ghosting} of the vertex and preserves local representation of functions. Its inverse transformation is called the \emph{splitting} of an edge. Splittings do not preserve local representations of functions.
\end{itemize}

An edit on a dendrogram $(T,\varphi_T)$ needs to be thought as a map $\{(T,\varphi_T)\}\rightarrow (\mathcal{T}, L_1(\R,E))$. 
So, if we have $e_0:\{(T,\varphi_T)\}\rightarrow (\mathcal{T}, L_1(\R,E))$ and $e_1:\{e_0(T)\}\rightarrow (\mathcal{T}, L_1(\R,E))$, we can compose the two edits.
This is what we mean by composition of edits.
Any finite composition of edits is referred to as an \emph{edit path}. 

Exploiting the definitions we have just given we can add some other details to the correspondence established by \Cref{cor:local_repr}, studying the relations between the ghosting defined in \Cref{defi:ghosting} and the one in \Cref{sec:edits_local_rep}.

\begin{defi}\label{defi:equiv_order_2}
Dendrograms are called equal up to degree $2$ vertices if they become isomorphic after applying a finite number of ghostings or splittings. We write $(T,\varphi_T)\cong_2(T',\varphi_{T'})$.  We call $(\mathcal{T}_2,L_1(\R,E))$ the space of equivalence classes of dendrograms in $(\mathcal{T},L_1(\R,E))$, equal up to degree $2$ vertices. 
\end{defi}

%Consider $\mathcal{M}(\T)=(T,h_T)$ and a proper weight function $\varphi_T$ with values in some editable space $L_1(\mathbb{R},E)$. We know that $\varphi_T$ is equivalent to a local representation of a function on $\mathcal{U}(D_{\T})$ which, in turns, amounts to the datum of a function $f$ in $L_1(D_{\T},E)$. By construction $\varphi_T(e)= f \circ (h_{\mid U})^{-1}:(t,t')\rightarrow E$, with $U\in\mathcal{U}(D_{\T})$ being associated to the edge $e$ and $h$ being the height function of $D_{\T}$. 

Looking at the definition of the ghosting edit, we can easily extend \Cref{cor:local_repr} to dendrograms up to degree $2$ vertices: if two dendrograms are local representations of the same function then they are equivalent up to degree $2$ vertices. 

\begin{cor}\label{cor:local_repr_2}
Consider $T=\mathcal{M}(\T)$ and $f\in L_1(D_{\T},E)$. Let the dendrogram $(T,\varphi_T)$ be any local representation of $f$, and call $[(T,\varphi_T)]$ the equivalence class of $(T,\varphi_T)$ up to degree $2$ vertices. We have:
\[
[(T,\varphi_T)]\cap L_1(D_{\T},E)\cong\lbrace\varphi_{\mathcal{O}}:\mathcal{O} \rightarrow L_1(\R,E) \mid \mathcal{O}\in \Cov(\T)\text{ and } \varphi_{\mathcal{O}}\text{ loc. repr. of }f\rbrace.
\]
Thus, for all regular $\T$ we have an injective map $L_1(D_{\T},E)\hookrightarrow (\mathcal{T}_2,L_1(\R,E))$.
\end{cor}

As a consequence of this result, we can say that the ghosting and splitting edits for local representation of functions represent the combinatorial equivalent of the lattice operations in $\Cov(\T)$: with a splitting we are refining the local representation and with the ghosting we are looking at the function on a coarser a.e. covering. All these dendrograms identify a unique function in $L_1(D_{\T},E)$ and a unique equivalence class in $\mathcal{T}_2$.

\begin{rmk}
Note that not all dendrograms in $(\mathcal{T}, L_1(\R,E))$ are local representation of functions. In fact, in general, we do not have: $\supp(\varphi_{T}((v,v')))\subset [h_{T}(v),h_{T'}(v)], \forall (v,v')\in E_{T}$. Upon collecting (the ismorphism classes of) all regular abstract merge trees into the set $\mathcal{RAMT}$, when can write this fact as the following map being injective but not surjective:
\[
\bigcup_{\T \in \mathcal{RAMT}}L_1(D_{\T},E) \rightarrow (\mathcal{T}_2, L_1(\R,E)).
\]
\end{rmk}

\subsection{Costs of Edit Operations}

Now we associate to every edit a cost so that we can measure distances between objects in $(\mathcal{T},L_1(\R,E))$. 
The costs of the edit operations are defined as follows:
\begin{itemize}
\item if, via shrinking, an edge goes from weight $f$ to weight $g$, then the cost of such operation is $d_{L_1}(f,g)$;
\item for any deletion/insertion of an edge with local function equal to $f$, the cost is equal to $d_{L_1}(f,0)$;
\item the cost of ghosting operations is 
$\mid d_{L_1}(f+ g,0)-d_{L_1}(f,0)-d_{L_1}(g,0)\mid =0$.
\end{itemize}

%The cost of an edit path $e_0\circ\ldots\circ e_n(T)$ is given by the sum of the costs of the single edit operations.

\begin{defi}
Given two dendrograms $T$ and $T'$ in $(\mathcal{T},L_1(\R,E))$,  define:
\begin{itemize}
\item $\Gamma(T,T')$ as the set of all finite edit paths between $T$ and $T'$;
\item $cost(\gamma)$ as the sum of the costs of the edits for any $\gamma\in\Gamma(T,T')$;
\item the dendrogram edit distance as:
\[
d_E(T,T')=\inf_{\gamma\in\Gamma(T,T')} cost(\gamma)
\]
\end{itemize}
\end{defi}

The following result, adapted from \cite{pegoraro2023edit}, together with \Cref{cor:local_repr_2}, says that $d_E$ is a metric for functions defined on display posets.

\begin{teo}[adapted from \cite{pegoraro2023edit}]\label{teo:metric} 
Given $E$ editable space, 
$((\mathcal{T}_2,L_1(\R,E)),d_E)$ is a metric space.
\end{teo}

Putting together \Cref{cor:local_repr_2} and 
\Cref{teo:metric} we have thus obtained a metric to compare $f\in L_1(D_{\T},E)$
and $g\in L_1(D_{\G},E)$.

\subsection{Mappings}
\label{sec:mappings}

We give the definition of a combinatorial tool used to compute $d_E$, called \emph{mapping}, taken from \cite{pegoraro2023edit}, as it is needed in \Cref{sec:optima_paths} and \Cref{sec:stability}.

\begin{defi}
A mapping between two dendrograms $T$ and $T'$ is a set 
$M\subset (E_T \cup \{\mathfrak{D},\mathfrak{G}\})\times (E_{T'} \cup \{\mathfrak{D},\mathfrak{G}\})$ 
satisfying:

\begin{itemize}
\item[(M1)] consider the projection of the Cartesian product $(E_T \cup \{\mathfrak{D},\mathfrak{G}\})\times (E_{T'} \cup \{\mathfrak{D},\mathfrak{G}\})\rightarrow (E_{T} \cup \{\mathfrak{D},\mathfrak{G}\})$; we can restrict this map to $M$ obtaining $\pi_T:M\rightarrow (E_T \cup \{\mathfrak{D},\mathfrak{G}\})$. The maps $\pi_T$ and $\pi_{T'}$ are surjective on $E_T$ and 
$E_{T'}$, i.e. $E_T\subset \im(\pi_T)$ and $E_{T'}\subset \im(\pi_{T'})$;
\item[(M2)]$\pi_T$ and $\pi_{T'}$ are injective on $M\cap (E_T\times E_{T'})$;
\item[(M3)]  given $(a,b)$ and $(c,d)\in M\cap (V_T\times V_{T'})$, $a>c$, if and only if $b>d$;
\item[(M4)] if $(a,\mathfrak{G})\in M$ (or analogously $(\mathfrak{G},a)$), 
then after applying all deletions of the form $(v,\mathfrak{D})\in M$, the vertex $a$ becomes a degree $2$ vertex. In other words: let $child(a)=\{b_1,..,b_n\}$. Then there is exactly one $i$ such that for all $j \neq i$, for all $v \in V_{sub(b_j)}$, we have $(v,\mathfrak{D})\in M$; and there is one and only one $c$ such that $c=\max\{x<b_i\mid (x,y)\in M$ for any $y \in V_{T'}\}$. 
\end{itemize}  

We call $\Mapp(T,T')$ the set of all mappings between $T$ and $T'$. 
\end{defi}

We may refer to edges which appear in the pairs in $M\cap (V_T\times V_{T'})$ as the \emph{paired} or \emph{matched} edges/vertices. 
Every $M\in \Mapp(T,T')$ parametrizes a set of edit paths, with identical cost, as follows:
\begin{itemize}
\item $\gamma_{d}^T$ is made by the deletions to be done on $T$, that is, the pairs $(v,\mathfrak{D})$, executed in any order. So we obtain $T^M_d=\gamma_{d}^T(T)$, which is well defined and does not depend on the order of the deletions.
Similarly, we define $\gamma_{d}^{T'}$ as a path made by the deletions to be done on $T'$, that is, the pairs $(\mathfrak{D},w)$, executed in any order, and obtain $T'^M_d=\gamma_{d}^{T'}(T')$. 
\item One then proceeds editing $T^M_d$ by ghosting all the vertices $(v,\mathfrak{G})$ in $M$, in any order, getting a path $\gamma^T_g$ and the dendrogram $T_M:= \gamma^T_g\circ \gamma^T_d (T)$. As before, we can do an analogous procedure on $T'^M_d$, ghosting all the vertices $(\mathfrak{G},w)$ in $M$, in any order, and building a path $\gamma^{T'}_g$ along with the dendrogram $T'_M:= \gamma^{T'}_g\circ \gamma^{T'}_d (T')$.
\item Since all the remaining points in $M$ are paired, the two dendrograms $T'_M$ and $T_M$ must be isomorphic as tree structures. This is guaranteed by the properties of $M$. So one can shrink $T_M$ onto $T'_M$, and the composition of the shrinkings, executed in any order is an edit path $\gamma_s^T$.
\end{itemize}

By construction $\gamma_s^T\circ \gamma_g^T\circ \gamma_d^T(T)=T'_M$,
and $(\gamma_d^{T'})^{-1}\circ (\gamma_g^{T'})^{-1}\circ \gamma_s^T\circ \gamma_g^T\circ \gamma_d^T (T)=T'$.
Where the inverse of an edit path is thought as the composition of the inverses of the single edit operations, taken in the inverse order.

Lastly, we call $\gamma_M$ the set of all possible edit paths of the form $(\gamma_d^{T'})^{-1}\circ (\gamma_g^{T'})^{-1}\circ \gamma_s^T\circ \gamma_g^T\circ \gamma_d^T$, 
obtained by changing the order in which the edit operations are executed inside $\gamma_d$, $\gamma_g$ and $\gamma_s$.
Even if $\gamma_M$ is a set of paths, its cost is well defined:

\[
cost(M):= cost(\gamma_M)=cost(\gamma_d^T)+cost(\gamma_s^T)+cost(\gamma_d^{T'}).
\]

In \cite[Main Theorem]{pegoraro2023edit}, it is proven that 
the $d_E(T,T')$ can always be realized with a mapping. Lastly, we also recall the following definition and result.

\begin{defi}[\cite{pegoraro2023edit}]\label{defi:M_2_mapp}
A mapping $M\in \Mapp(T,T')$ has maximal ghostings if the following hold: $(v,\mathfrak{G})\in M$ if, and only if, $v$ is of degree $2$ after the deletions in $T$ and, similarly, $(\mathfrak{G},w)\in M$ if, and only if, $w$ is of degree $2$ after the deletions in $T'$.

A mapping $M\in \Mapp(T,T')$ has minimal deletions if the following hold: $(v,\mathfrak{D})\in M$ implies that neither $v$ nor $parent(v)$ are of degree $2$ after applying all the other deletions in $T$ and, similarly, $(\mathfrak{D},w)\in M$ implies that neither $w$ nor $parent(w)$ are of degree $2$ after applying all the other deletions in $T'$.

We collect all mappings with maximal ghostings and minimal deletions in the set $M_2(T,T')$.
\end{defi}

\begin{lem}[\cite{pegoraro2023edit}]\label{lemma:M_2}
\[
\min\{cost(M)\mid  M\in \Mapp(T,T') \}=\min\{cost(M)\mid M\in M_2(T,T')\}
\]
\end{lem}

\subsection{Optimal Edit Paths}
\label{sec:optima_paths}

\begin{figure}	
    \centering
    \begin{subfigure}{\textwidth}
    \centering
    	\includegraphics[width = 0.8\textwidth]{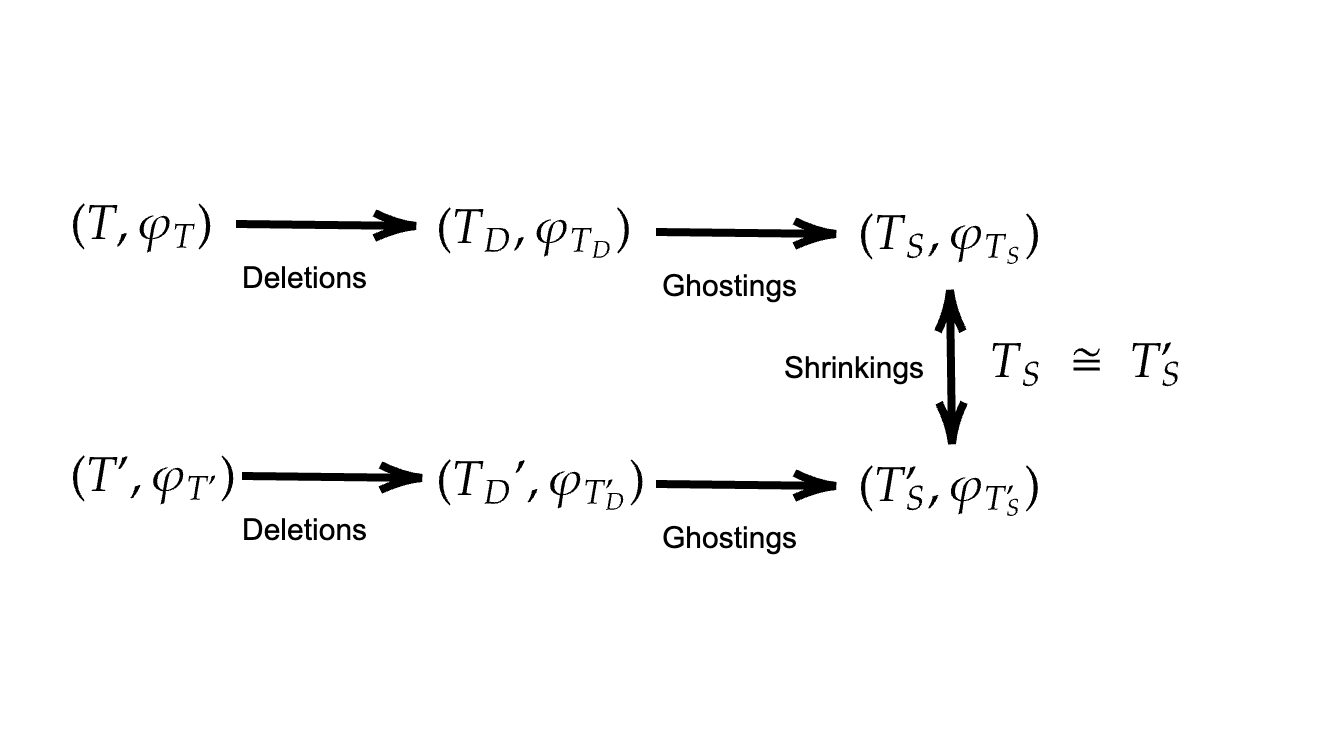}
	\end{subfigure}
%	\hspace{0.5 cm}
\caption{A representation of particular optimal edit paths between dendrograms. If $(T,\varphi_T),(T',\varphi_{T'})$ are local representations of functions, then so are $(T_S,\varphi_{T_S})$ and $(T'_S,\varphi_{T'_S})$.}
\label{fig:pipeline_edit}
\end{figure}

In \Cref{sec:edits_local_rep}, we noted that starting from a local representation of $f \in L_1(D_{\T}, E)$, the application of shrinkings, splittings, or insertions generally does not yield another local representation of a function. 

While this limitation could potentially be addressed in future work, we argue that it does not hinder the definition of a meaningful metric structure for comparing functions $f \in L_1(D_{\T}, E)$ and $g \in L_1(D_{\G}, E)$.
As shown in \Cref{sec:mappings}, there always exists a minimal edit path that proceeds as illustrated in \Cref{fig:pipeline_edit}. Given a starting dendrogram $(T, \varphi_T)$ and a target dendrogram $(T’, \varphi_{T’})$, one can first apply all deletions and ghostings to $T$ and $T’$, respectively, resulting in the intermediate dendrograms $(T_S, \varphi_{T_S})$ and $(T’S, \varphi{T’_S})$. If the original dendrograms are local representations of functions, then each dendrogram along these edit paths also remains a local representation. Hence, no metric artifact is introduced at this stage.

Furthermore, the properties of these optimal edit paths ensure that the resulting tree structures $T_S$ and $T’S$ are isomorphic, and the shrinkings precisely map each edge value $\varphi{T_S}(e)$ to $\varphi_{T’S}(e’)$ for some $e \in E{T_S}$ and $e’ \in E_{T’_S}$. In this way, the comparison between $f$ and $g$ reduces to a direct, edge-by-edge comparison of their values over isomorphic trees.

\section{Examples}
\label{sec:info_examples}

\begin{figure}
	\begin{subfigure}[c]{0.47\textwidth}
    	\centering
    	\includegraphics[ width = \textwidth]{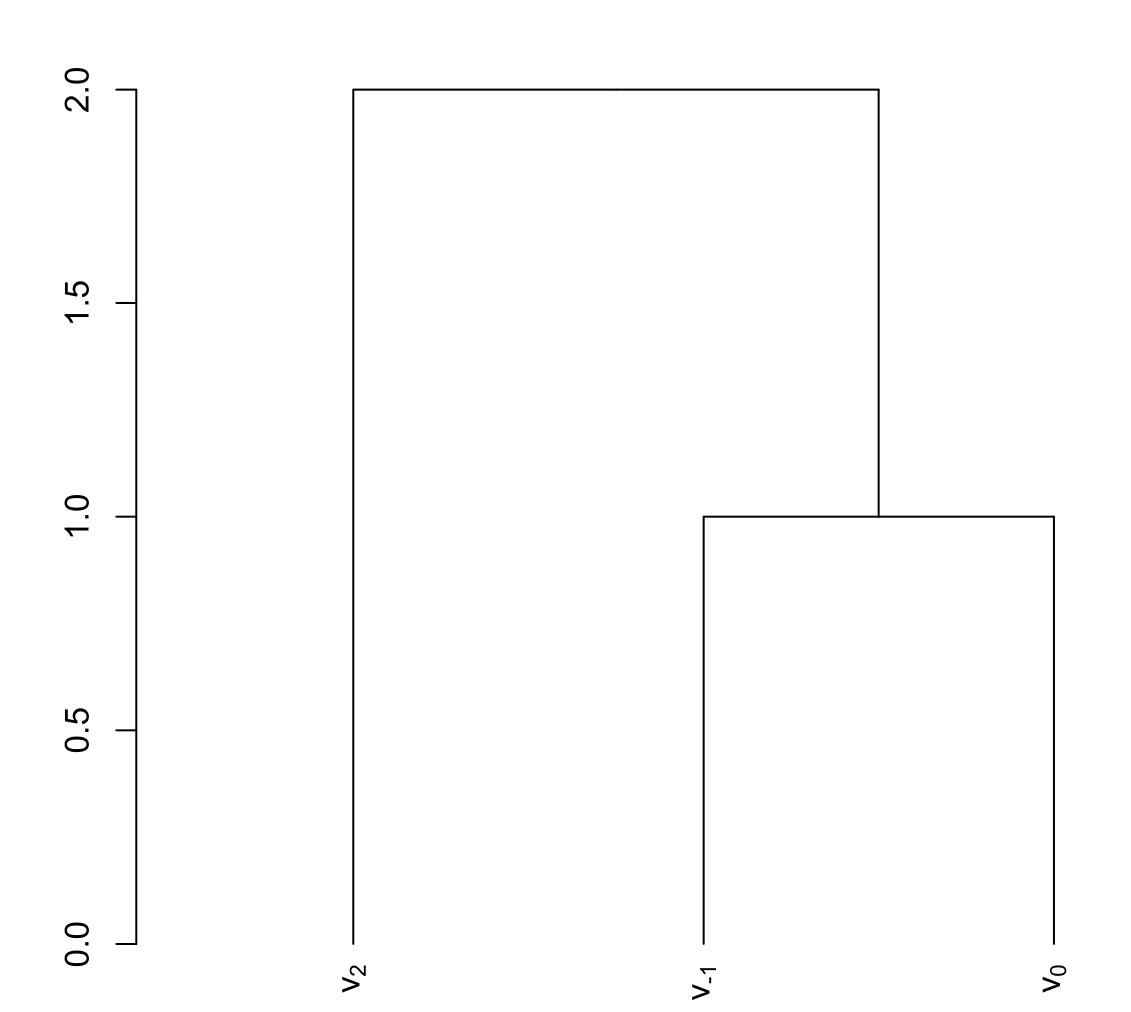}
            \captionsetup{singlelinecheck=off, margin={0.3cm, 0.1cm}}
    	\caption{Single linkage clustering dendrogram referring to the example in \Cref{sec:dendro_clus}.}
    	\label{fig:clus_dendro}
    \end{subfigure}
%    \hspace{1 cm}
    \begin{subfigure}[c]{0.47\textwidth}
    	\centering
    	\includegraphics[ width = \textwidth]{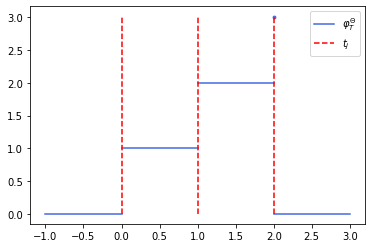}
            \captionsetup{singlelinecheck=off, margin={0.3cm, 0.1cm}}
    	\caption{In the context of the example in \Cref{sec:dendro_clus}, we see the sum of the weight functions of the vertices going from $v_0$ to the root $r_T$: $\varphi^{\Theta_c}_T(\{v_0\})+\varphi^{\Theta_c}_T(\{v_0,v_{-1}\})$. The dotted lines represent critical values.}
    	\label{fig:clus_mult}
    \end{subfigure}

%    \hspace{1 cm}   
    \centering
    \begin{subfigure}[c]{0.47\textwidth}
    	\centering
    	\includegraphics[ width = \textwidth]{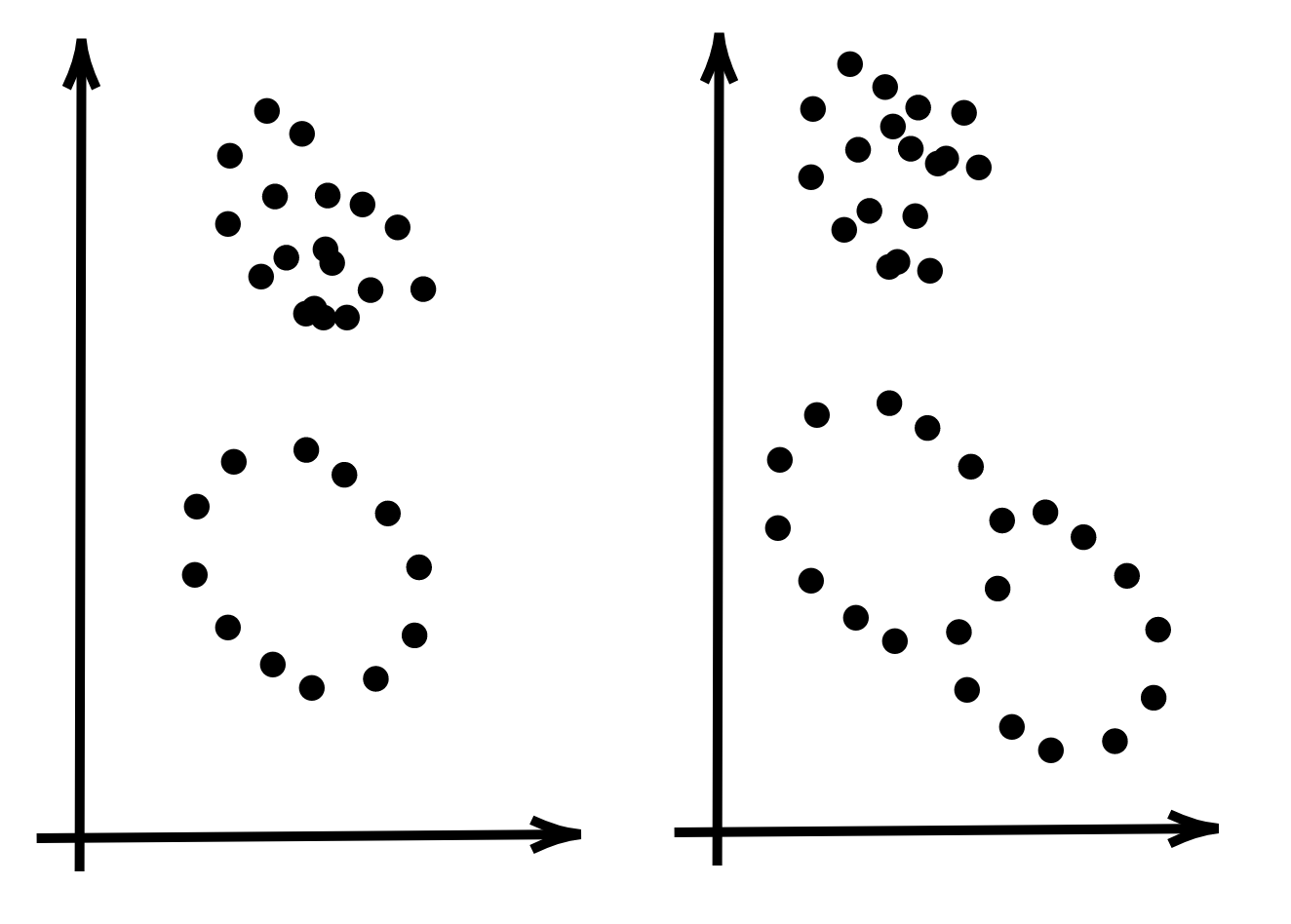}
    	\caption{Two point clouds made by two clusters each which cannot be separated by zero dimensional homology, but present different within-cluster homological information and can be distinguished by $\Theta_1$ defined in \Cref{sec:decorated}.}
    	\label{fig:homology}
    \end{subfigure}

\caption{Plots referring to the examples in \Cref{sec:info_examples}.}
\label{fig:back_to_vec}
\end{figure}

In this section we want to give some concrete examples of functions defined on display posets, to show how they can be used to capture useful information about a filtration $\X$. This section is complemented by \Cref{sec:simulated_data}, where we present two simulated scenarios in which we test some of the upcoming ideas, and 
\Cref{sec:example_edit_fn}, which shows the problems arising when trying to replicate this framework for PDs.

The general structure of the following examples is to consider a subcategory $\mathcal{B}$ of $\Top$ and pick a function $\Theta:\mathcal{B}\rightarrow E$. Then, $f:D_{\T}\rightarrow E$ is obtained as $f(a,t)=\Theta(a)$. We call $\varphi^{\Theta}_T$ the local representation of such function, and we prove that in all our examples the information contained in the functions generalizes, in some sense, the notion of merge trees.
More formally, having $(T,\varphi^{\Theta}_T)\cong (T',\varphi^\Theta_{T'})$ implies $(T,h_T)\cong (T',h_{T'})$. Under such hypotheses a metric to compare $(T,\varphi^{\Theta}_T)$ and $(T',\varphi^\Theta_{T'})$ can be pulled back to compare objects of the form $(T,h_T,\varphi_T)$ - or, equivalently, $f\in L_1(D_{\T},E)$
and $g\in L_1(D_{\G},E)$.

We immediately stress that many of the upcoming functions do not lie in $L_p(D_{\T},E)$, for some $D_{\T}$, as: 

\[
\lim_{x\rightarrow +\infty} d(f \circ (h_{\mid U_{\infty}})^{-1}(x),0)>0.
\] 

In \Cref{sec:normalizing_functions} we discuss how these examples can be modified to fit into the proposed framework.

\subsection{Cardinality of Clusters}
\label{sec:dendro_clus}

Consider the case of a merge tree $\mathcal{M}(\T)=(T,h_T)$, with $\X$ being the $\check{C}$ech filtration of the point cloud $\{x_1,\ldots, x_n\}$. 
Sensible information that one may want to track down along 
 $\T$ is the cardinality of the clusters. 
 Thus, we can take $\Theta_c:\FSet \rightarrow \mathbb{R}_{\geq 0}$, defined on all finite sets considered with the discrete topology, defined as $\Theta_c(A)=\# A$.
As a consequence, we have $\varphi^\Theta_T(e) = m\chi_{[t_i,t_j)}$, for some positive cardinality $m$ and some critical values $t_i,t_j$. Note that, clearly, 
$\supp(\varphi^{\Theta_{c}}_T(e))=[t_i,t_j]$.
Thus if we have $(T,\varphi^{\Theta_c}_T)\cong (T',\varphi^{\Theta_c}_{T'})$ then $(T,h_T)\cong (T',h_{T'})$. 

\bigskip

We now make a concrete example - see \Cref{fig:clus_dendro} and \Cref{fig:clus_mult}. Consider the finite set $\{v_{-1},v_0,v_2\}$ with $v_{-1}=-1,v_0=0,v_2=2\in \R$, and build the single linkage hierarchical clustering dendrogram along with the associated abstract merge tree.
The abstract merge tree is given by $a_t=\{\{v_{-1}\},\{v_0\},\{v_2\}\}$ for $t\in [0,1)$, 
$a_t=\{\{v_{-1},v_0\},\{v_2\}\}$ for $t\in [1,2)$ and 
$a_t=\{\{v_{-1},v_0,v_2\}\}$ for $t\geq t^+=2$.
With maps given by $x \mapsto y$ if $x\subset y$.

The merge tree $(T,h_T)$ one obtains - see \Cref{fig:clus_dendro} -  can be represented with the vertex set $V_T=\{\{v_{-1}\},\{v_0\},\{v_2\},\{v_{-1},v_0\},\{v_{-1},v_0,v_2\},r_T\}$.
The leaves are $\{v_{-1}\},\{v_0\}$ and $\{v_2\}$; the children of $\{v_{-1},v_0\}$ are $\{v_{-1}\}$ and $\{v_0\}$, 
and the ones of $\{v_{-1},v_0,v_2\}$ are $\{v_{-1},v_0\}$ and $\{v_2\}$.
The height function $h_T$ is given by $h_T(\{v_i\})=0$ for $i=-1,0,2$, $h_T(\{v_{-1},v_0\})=1$, $h_T(\{v_{-1},v_0,v_2\})=2$ and $h_T(r_T)=+\infty$.

Consider $\Theta_c$.
The local representation $\varphi^{\Theta_c}_T$ of the induced function is thus the following: 
$\varphi^{\Theta_c}_T(\{v_i\})=\chi_{[0,1)}$ for $i=-1,0$,
$\varphi^{\Theta_c}_T(\{v_2\})=\chi_{[0,2)}$, 
$\varphi^{\Theta_c}_T(\{v_{-1},v_0\})=2\chi_{[1,2)}$ and 
$\varphi^{\Theta_c}_T(\{v_{-1},v_0,v_2\})=3\chi_{[2,+\infty)\}}$. See \Cref{fig:clus_mult}.

\subsection{Measure of Sublevel Sets}
\label{sec:dendro_fun}

Now consider $U\subset\mathbb{R}^m$ convex bounded open set, with $\overline{U}$ being its topological closure, 
and let $\mathcal{L}$ be the Lebesgue measure in $\mathbb{R}^m$. Let $f:\overline{U}\rightarrow \mathbb{R}$ be a tame continuous function. Consider the sublevel set filtration $X_t = f^{-1}((-\infty,t])$ with $\pi_0(X_t) =\{U^t_1,\ldots, U^t_n \}$. Call $\psi_t^{t'}$ the functions 
$\psi_t^{t'}=X_{t\leq t'}$.
We set $\Theta_{\mathcal{L}}=\mathcal{L}$,
that is:
\begin{align*}
\Theta_{\mathcal{L}}:\mathcal{B}(\mathbb{R}^n)&\rightarrow \mathbb{R}_{\geq 0}   \\
U &\mapsto \mathcal{L}(U),
\end{align*}
with $\mathcal{B}(\mathbb{R}^n)$ being the Borel $\sigma$-algebra of $\mathbb{R}^n$. By construction, we can always consider: $\Theta_{\mathcal{L}}(U^t_i)=\mathcal{L}(U^t_i)$.

\begin{prop}
If we have $(T,\varphi^{\Theta_{\mathcal{L}}}_T)\cong (T',\varphi^{\Theta_{\mathcal{L}}}_{T'})$ then $(T,h_T)\cong (T',h_{T'})$.
\begin{proof}
Let $(T,h_T)$ being the merge tree representing  $\T$, and $\varphi^{\Theta_{\mathcal{L}}}_T$ the local representation of the associated function.
Since $f$ is continuous, for and edge $e=(v,v')\in E_T$ spanning from height $h_T(v)=t_i$ to $h_T(v')=t_j$, we can prove that $\supp(\varphi^{\Theta_{\mathcal{L}}}_T(e))=[t_i,t_j]$. We know that $v$ is associated to a connected component $U^{t_i}_k$, for some $k$. If $v$ represents the merging of two or more path-connected components $U^{t_i-\varepsilon}_{k_1}$ and $U^{t_i-\varepsilon}_{k_2}$, for some small $\varepsilon>0$, with $\mathcal{L}(U^{t_i-\varepsilon}_{k_1}),\mathcal{L}(U^{t_i-\varepsilon}_{k_2})>0$, then, since $U^{t_i-\varepsilon}_{k_1},U^{t_i-\varepsilon}_{k_2}\subset U^{t_i}_{k}$, we have $\mathcal{L}(U^{t_i}_{k})>0$. Thus if we prove the statement for $v$ leaf, we are done.
So, suppose $v$ is a leaf and consider $x_0\in U^{t_i}_{k}$. We know $f(x_0)=t_i$. By the continuity of $f$, for every $\varepsilon>0$ there is $\delta>0$ such that if $\mid \mid x-x_0\mid \mid <\delta$, then $f(x_0)\leq f(x)<f(x_0)+\varepsilon$. Since $\{x\in \overline{U} \mid  \mid \mid x-x_0\mid \mid <\delta \}$ is convex (and so path-connected), then it is contained in $\psi_{t_i}^{t_i+\varepsilon}(U^{t_i}_{k})$.
Moreover, since it contains the non-empty open set $\{x\in U \mid \text{ } \mid \mid x-x_0\mid \mid <\delta \}$, we have 
$\mathcal{L}(\psi_{t_i}^{t_i+\varepsilon}(U^{t_i}_{k}))>0$ for every $\varepsilon>0$. As a consequence, $\supp(\varphi^{\Theta_{\mathcal{L}}}_T(e))=[t_i,t_j]$.
\end{proof}
\end{prop}

Again we make a quick hands-on example.
Consider the function $f=\mid \mid x\mid -1\mid $ defined on the interval $[-2,2]$. Let $\pi_0(X_t)=\pi_0(f^{-1}((-\infty,t]))$.
Let $(T,h_T)$ be the merge tree associated to the sequence $\T$. Now we obtain the local representation $\varphi^{\Theta_\mathcal{L}}_{T}(e_i)$.

We have $\varphi^{\Theta_\mathcal{L}}_{T}(e_1)=\mid 1+t-1+t\mid =2t$ for $t\in [0,1)$, and $0$ otherwise. Clearly $\varphi^{\Theta_\mathcal{L}}_{T}(e_1)=\varphi^{\Theta_\mathcal{L}}_{T}(e_2)$. Lastly $\varphi^{\Theta_\mathcal{L}}_{T}(r_T)=4\chi_{[2,+\infty)}$.

\subsection{Homological Information}
\label{sec:decorated}

Lastly, we propose a function $\Theta_p$ to combine 
homological information \cite{hatcher} of different dimensions obtaining dendrograms which are closely related to the barcode decorated merge trees defined by \cite{curry2021decorated} and the decorated mapper graphs defined in \cite{curry2023topologically}. 
We consider the topological spaces with $p$-th homology of finite type, that is, their $p$-th homology group is finitely generated, and collect all the spaces with finitely generated $1,\ldots,p$-th homology groups in the set $\FTop_p$.
Consider $\Theta_p: \FTop_p \rightarrow \mathbb{N}\times \ldots\times \mathbb{N}$ defined on a topological space $U$ as the component-wise Betti function $\Theta_p(U)=(\dim(H_0(U;\mathbb{K}),\ldots,\dim(H_p(U;\mathbb{K}))$, with $H_p(U;\mathbb{K})$ being the $p$-th homology group of $U$ with coefficients in the field $\mathbb{K}$. 
Note that, by definition, generators of homology groups of $U$ lie inside a path-connected component.
In this way we are able to track if in a path-connected component there are some kind of holes arising or dying, and thus collecting a more complete set of topological invariants which capture the shape of each path-connected component. From another point of view, at every step along a filtration, we are decomposing homological information of a topological space by means of its path-connected components. This, for instance, could be useful in situations like the one depicted in \Cref{fig:homology}

Note that we clearly have $(T,\varphi^{\Theta_p}_T)\cong (T',\varphi^{\Theta_p}_{T'})$ implying $(T,h_T)\cong (T',h_{T'})$. In fact, considering $\FTop_p \xrightarrow{\Theta_p} \mathbb{N}\times\ldots\times \mathbb{N}\xrightarrow{\pi_1} \mathbb{N}$ (with $\pi_1$ being the projection on the first component), is equivalent to taking $\Theta_0 \equiv 1$ since all the sets on which we evaluate $\Theta_0$ are path-connected. Thus, $\supp (\varphi_T^{\Theta_0}(e))=[h_T(v),h_T(v')]$, for every $e\in E_T$.

\subsection{Normalizing and Truncating Functions}
\label{sec:normalizing_functions}

The functions described in the previous sections do not fit in our framework, as they are away from $0$ at infinity. 
However, one can clearly truncate such functions at some height $K$ in order to obtain a a function with finite integral.
That is, 
given $f:D_{\T}\rightarrow E$, with $\mathcal{M}(\T)=(T,h_T)$, we can always put $f$ to $0\in E$, after some $K\geq \max h_T$, considering $f_{\mid U_\infty}\cdot \chi_{[\max h_T,K]}$.
We call $f_{\mid K}$ the function defined on $U\in \mathcal{U}(D_{\T})$ as:
\[
(f_{\mid K})_{\mid U }:= \begin{cases}
f_{\mid U}\cdot 
\chi_{[\max h_T,K]}\text{ if }U= U_\infty, \\
    f_{\mid U}, \text{ otherwise.}
\end{cases}
\]

Clearly $f_{\mid K}\in L_1(D_{\T},E)$.
Although, in general, this solution is rather artificial and $d_E(f_{\mid K},g_{\mid K})$ depends on $K$, we show that if $f,g$ are definitively equal going upward toward the roots,
then $d_E(f,g):=d_E(f_{\mid K},g_{\mid K})$, for some $K$ big enough.
More formally, suppose we have $f:D_{\T}\rightarrow E$,
$g:D_{\G}\rightarrow E$, with the height functions of the display posets being respectively $h^f$ and $h^g$, and suppose there is $K\in\mathbb{R}$
such that, for all $x>K$:
\begin{equation}\label{eq:equal_at_inf}
f\circ (h^f_{\mid U^f_\infty})^{-1}(x) = g\circ (h^g_{\mid U^g_\infty})^{-1}g(x),    
\end{equation}
then $d_E(f,g):=d_E(f_{\mid K},g_{\mid K})$.

Before proving such result, we highlight that $\Theta_c$, 
$\Theta_\mathcal{L}$, and $\Theta_p$, in many relevant situations, can be modified so that the resulting function  satisfy \Cref{eq:equal_at_inf}. In particular:
\begin{itemize}
    \item $\Theta_c(X_t)$ can be normalized by the cardinality of the considered point cloud $X$, so that it is definitively equal to $1$;
    \item $\Theta_\mathcal{L}(X_t)$ can be normalized by the total measure of the domain of $f$;
    \item $\Theta_p$ satisfies \Cref{eq:equal_at_inf} anytime, for two filtrations $\X$ and $\Y$, the topological spaces $\bigcup_{t\in \R}X_t$ and $\bigcup_{t\in \R}Y_t$ have isomorphic homology groups.
\end{itemize}

Now we prove the main result of the section.

\begin{prop}[Truncation]\label{prop:truncation}
Take $(T,\varphi_T)$ and  $(T',\varphi_{T'})$.
Suppose $r_T$ and $r_{T'}$ are of degree $1$ and there is a splitting $\{(v,r_T)\}\rightarrow \{(v,v'),(v',r_T)\} $ and $\{(w,r_{T'})\}\rightarrow \{(w,w'),(w',r_{T'})\} $ giving the dendrograms $(G,\varphi_{G})$ and $(G',\varphi_{G'})$. Suppose moreover that $\varphi_{G}((v',r_{T}))=\varphi_{G'}((w',r_{G}))$. Then $d_E(T,T')=d_E(sub_{G}(v'),sub_{G'}(w'))$.
\begin{proof}
    See \Cref{sec:proofs_TDA}.
\end{proof}
\end{prop}

 To conclude, consider $(T,h_{T},\varphi_T)$ and $(G,h_{G},\varphi_G)$ being local representations of, respectively, $f$ and $g$. Suppose we can split $e_T=(v,r_T)\in E_T$ into $e'_T=(v,v')$, $e''_T=(v',r_T)$, and $e_G=(w,r_G)\in E_G$ into $e'_G=(w,w')$, $e''_G=(w',r_G)$, so that
 $\varphi_T(e''_T)=\varphi_G(e''_G)$. Note that, this implies $h_T(v')=h_G(w')$. We call $K=h_T(v')=h_G(w')$.

Let $(T',h_{T'})$ and $(G',h_{G'})$ be the merge trees obtained with such splittings.
If we call $\varphi_{T'}$ and $\varphi_{G'}$ the local representations of $f$ on $T'$ and $g$ on $G'$, respectively, we have:
$\parallel \varphi_{T'}(e_T'')-
 \varphi_{G'}(e_G'')
\parallel_{L_1(\mathbb{R},E)}=0$. Thus, we are in the position to apply \Cref{prop:truncation} to $T'$ and $G'$
and truncate $f$ and $g$ from $K$ upward.
By \Cref{prop:truncation}, we are guaranteed that $d_E(f,g):=d_E(f_{\mid K},g_{\mid K})$.

\section{Stability}
\label{sec:stability}

In this section we establish some stability properties for the metric $d_E$ between functions defined on merge trees.

Establishing stability results across the wide range of scenarios we consider would require analyzing numerous specific conditions, each tied to the details of the function-generating pipelines. This lies beyond the scope of the present work.

For this reason, we focus on a general setting (arguably the canonical one for stability results in TDA), which does not require application-specific assumptions, and find a result involving the pipeline described in \Cref{sec:decorated}, which is analogous to the stability properties of the $1$-Wasserstein metric between persistence diagrams.

In order to state the result, we need to introduce some pieces of notation related to the interleaving distance \cite{merge_interl} between merge trees.

\begin{defi}
    For any functor $F$ defined on $\R$, we define $\mathcal{S}_\varepsilon F$ as the functor $\mathcal{S}_\varepsilon F(t)=F(t+\varepsilon)$.
\end{defi}

\begin{defi}[adapted from \cite{merge_interl}]
Take two abstract merge trees $\T$ and $\G$.
Two natural transformations $\alpha:\pi_0(\X)\rightarrow \mathcal{S}_\varepsilon\pi_0(\Y)$, $\beta:\pi_0(\Y)\rightarrow \mathcal{S}_\varepsilon\pi_0(\X)$ are $\varepsilon$-compatible if:

\begin{itemize}
\item $\beta_{t+\varepsilon}\circ \alpha_t = \pi_0(X_{t\leq t+2\varepsilon})$
\item $\alpha_{t+\varepsilon}\circ \beta_t = \pi_0(Y_{t\leq t+2\varepsilon})$.
\end{itemize}

Then, the interleaving distance between $\T$ and $\G$ is:
\[
d_I(\T,\G) = \inf \{\varepsilon>0\mid \exists \alpha,\beta \text{ }\varepsilon\text{-compatible}\}. 
\]

We also say that $\T$ and $\G$ are $d_I(\T,\G)$-interleaved. 
\end{defi}

Given a finite set $S = \{s_1,\ldots,s_n\}$, we define:
\[
\R-S^\varepsilon := \R-\bigcup_{i=1,\ldots,n} [s_i-\varepsilon,s_i+\varepsilon].
\]
Moreover, we set:
\[
\mu_{D_{\T}}(\R-S^\varepsilon):=\mu_{D_{\T}}(h^{-1}_{D_{\T}}(\R-S^\varepsilon)).
\]

\begin{figure}
    	\centering
	    \includegraphics[width = 0.8\textwidth]{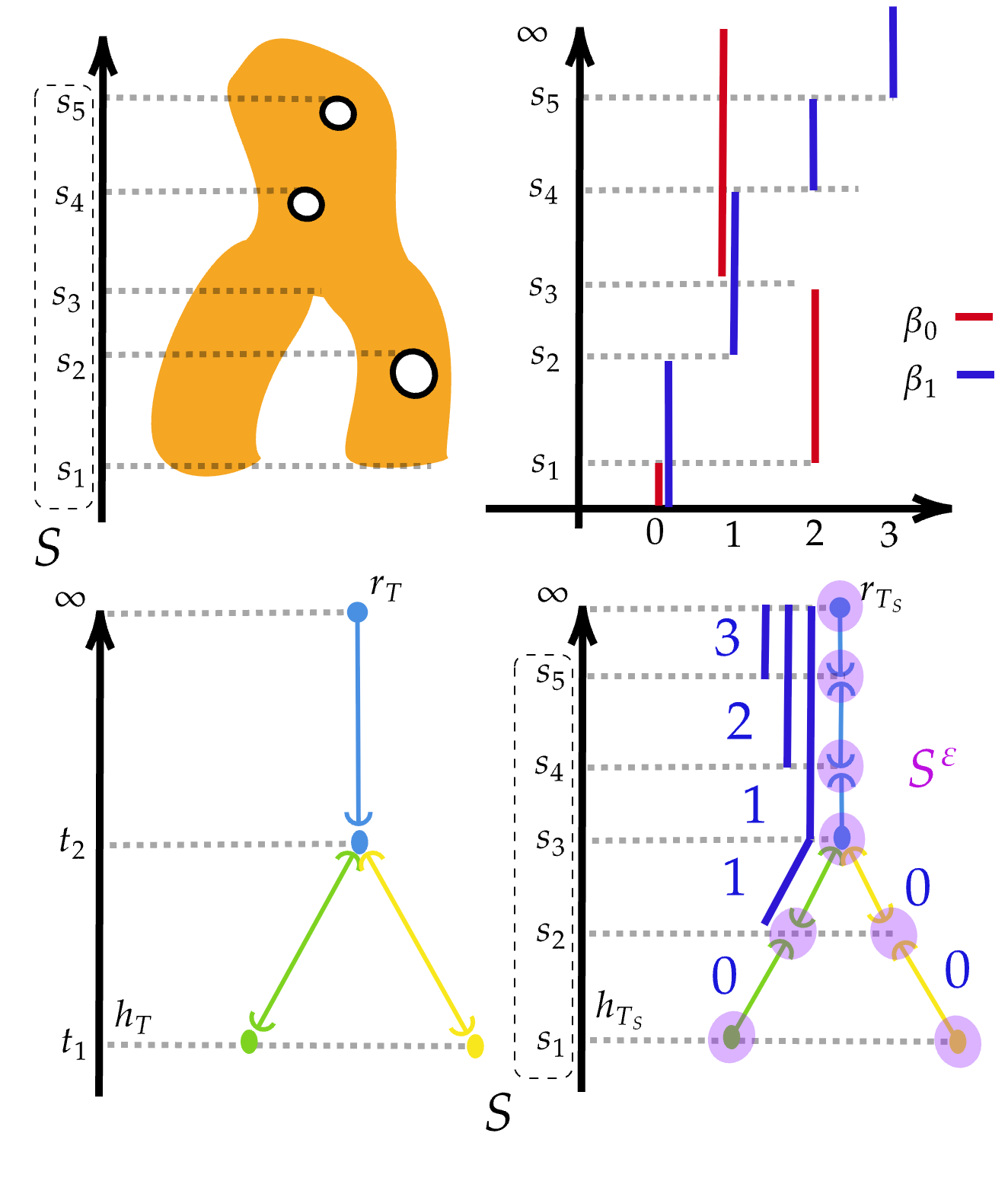}
    	\caption{A visual summary of the setting discussed in \Cref{sec:stability}. In the upper left, we depict an orange topological space equipped with a height function. The critical values of the persistence modules $H_0 \circ \X$ and $H_1 \circ \X$ are shown and collected into the set $S$. The upper right displays the persistent Betti functions $\beta_0$ and $\beta_1$, defined as $\dim(H_0 \circ \X)$ and $\dim(H_1 \circ \X)$, respectively. In the lower left, we show the merge tree $T$ associated with the height function. The lower right illustrates the merge tree $T_S$,  with $T_S\sim T$ induced by the pullback of $\R-S$ via $h_T$. Also highlighted is the set $S^\varepsilon$, which—according to \Cref{teo:fun_stab}—marks the region where the largest discrepancies between the compared functions may occur. Finally, we represent the persistence bars of $H_1 \circ \X$, which determine the function $\varphi^{\Theta_1}$.}
    	\label{fig:stability}
\end{figure}

Next, consider $\T$ abstract merge tree and $S$ finite set which contains the critical values of $\T$. We define the merge tree $T_S \sim_2 \mathcal{M}(\T)$ as the merge tree associated to the a.e. cover of $D_{\T}$ given by $\pi_0(h_{D_{\T}}^{-1}(\R-S))$. See \Cref{fig:stability}.

Consider $\mathcal{M}(\T)=(T,h_T)$ and $\mathcal{M}(\G)=(G,h_G)$. Following \cite{pegoraro2024functional}, we connect the continuous nature of $\G$ and the discrete one of $(G,h_G)$, via following map:
\begin{align*}
L_g:D_{\G}&\rightarrow E_{G} \\
    q &\mapsto \max\{w\in E_{G} \mid  w\leq q\}.
\end{align*}
Leveraging these definition, we set $\phi:= L_g\circ\alpha$.

Using this notation, we can re-formulate \cite[Theorem 1]{pegoraro2024functional} in a way which is more convenient for our purposes.

\begin{teo}[adapted from \cite{pegoraro2024functional}]\label{teo:stab_mapp}
    Suppose there exists an $\varepsilon$-interleaving between $\T$ and $\G$, given by two maps $\alpha$ and $\beta$.
    The map $\phi$ can be used to build a minimizing mapping $M$ between $\mathcal{M}(\T)=(T,h_T)$ and $\mathcal{M}(\G)=(G,h_G)$ such that $cost((a,b))\leq 2\varepsilon,$ for every $(a,b)\in M$.

In particular, this implies the existence of an edit path  such that:
\begin{itemize}
    \item each edge is involved in at most one deletion or shrinking;
    \item all the deleted edges of $T$ satisfy $h_T(parent(v))-h_T(v) \leq 2\varepsilon$, and all the deleted edges of $G$ satisfy $h_G(parent(w))-h_G(w) \leq 2\varepsilon$;
    \item all the shrinking edits, which shrink an edge $(v,parent(v))$ on an edge $(w,parent(w))$ (for some $v\in E_T$ and $w \in E_G$), are of the form $w= \phi(x)$ and they satisfy the following inequalities:
    \begin{align*}
        &\mid h_T(v)-h_G(w)\mid \leq \varepsilon, \\
        &\mid h_T(parent(v))-h_G(parent(w))\mid \leq \varepsilon,
    \end{align*}
    where $parent(v)$ and $parent(w)$ are computed after having applied all the deletions on $T$ and $G$, respectively, and all the possible ghostings. 
\end{itemize}
\end{teo}

We now prove the main result of this section, from which the stability of $\Theta_p$-induced functions will follow.

\begin{teo}\label{teo:fun_stab}
Consider $f\in L_1(D_{\T},(E,d'))$ and $g\in L_1(D_{\G},(E,d'))$, with local representations $(T,h_T,\varphi_T)$ and $(G,h_G,\varphi_G)$, for which we can find $C>0$ satisfying:
\[
  \sup_{t\in D_{\T}} d'(f(t),0), \sup_{t\in D_{\G}} d'(g(t),0) \leq C.
\]
Let $M$ be a mapping found via an $\varepsilon$-interleaving as in \Cref{teo:stab_mapp}.  
Suppose, lastly, that there exist a finite set $S = \{s_1,\ldots,s_n\}$  such that:
\begin{itemize}
    \item $S$ contains the critical values of $\T$ and $\G$;
    \item for every shrinking of an edge $e=(v,parent(v))$ on an edge $e'=(w,parent(w))$ (as before, these edges are obtained after all the deletions and ghostings), for every
    \[
    [t,t+\varepsilon]\subset \supp(e)\cap \supp(e')\cap(\R-S^\varepsilon)
    \]
    and for every $t' \in [t,t+\varepsilon]$,
    we have:
    \[
    d'(\varphi_T(e)(t'),\varphi_G(e')(t')) \leq  \varepsilon'.
    \]
\end{itemize}

Then:
\begin{equation}\label{eq:stability}
\begin{aligned}
 d_E(\varphi_T,\varphi_G)\leq &(\# E_{T_S}+\# E_{G_S}) \cdot 2C\varepsilon + \\
&\varepsilon' \cdot \min \{\mu_{D_{\T}}(\R-S^\varepsilon),\mu_{D_{\G}}(\R-S^\varepsilon)\}.      
\end{aligned}
\end{equation}

\begin{proof}
    Before starting the proof, we stress that we have 3 different metrics appearing: 
    \begin{enumerate}
        \item $d'$, the metric in the editable space $E$;
        \item $d_{L_1}$, the metric in $L_1(\R,E)$;
        \item $d_E$, the edit distance in $(\mathcal{T}_2,E)$.
    \end{enumerate}

    Thanks to \Cref{teo:stab_mapp}, we know that all deletions, which can not be more than $(\# E_T + \# E_{G})\leq (\# E_{T_S}+\# E_{G_S})$, have a cost which is bounded from above by $2C\varepsilon$.
    Thus, we only need to check what happens with shrinkings.
    
    Consider the shrinking of an edge $e=(v,parent(v))$ on an edge $e'=(w,parent(w))$.
    Set $I_e = \supp(\varphi_T(e))$, $I_{e'} = \supp(\varphi_G(e'))$, and $I=I_e\cap I_{e'}$.
    Let $A = \{s_1,\ldots,s_k\}=I\cap S$. 
    
    \begin{rmk}\label{rmk:critical_values}
        Note that the extremes of $I_e$ and $I_{e'}$ are contained in $[s_1-\varepsilon,s_1+\varepsilon]\cup [s_k-\varepsilon,s_k+\varepsilon]$.    
    \end{rmk}
    
     We can write: 
    \[
    d_{L_1}( \varphi_T(e),\varphi_G(e') ) = \int_{\R-A^\varepsilon} d'(\varphi_T(e)(t),\varphi_G(e')(t)) d\mathcal{L}(t) +  \int_{A^\varepsilon} d'(\varphi_T(e)(t),\varphi_G(e')(t)) d\mathcal{L}(t)
    \]
    By hypothesis, we know that:
    \[
    \int_{A^\varepsilon} d'(\varphi_T(e)(t),\varphi_G(e')(t)) d\mathcal{L}(t) \leq 2C\varepsilon k.
    \]
    Similarly:
    \[
    \int_{\R-A^\varepsilon} d'(\varphi_T(e)(t),\varphi_G(e')(t)) d\mathcal{L}(t) \leq \varepsilon' \mathcal{L}(I-A^\varepsilon).
    \]
    Since the edges are involved at most once in shrinking edits, adding up these contributions for all edges in $T$ and $G$, and using the additive properties of the measures $\mu_{D_{\T}}$ and $\mu_{D_{\G}}$ lead to the result. 
\end{proof}
\end{teo}

We briefly comment on \Cref{eq:stability}: 
\begin{itemize}
    \item first of all, such inequality does not give a stability result unless we can control $\varepsilon'$ and guarantee that it goes to $0$ under some hypotheses;
    \item if we can find appropriate conditions to control $\varepsilon'$, then, the first addend in  \Cref{eq:stability} controls the differences between functions where the values of the functions don't go to $0$ (e.g. where we have deletions or differences between the supports between $\varphi_T(e))$ and $\varphi_G(e')$). Locally, these appear on intervals of length at most $2\varepsilon$, and there is at most one such interval for every edge in $T_S$ and $G_S$. The right addend in \Cref{eq:stability}, instead, controls the difference between $f$ and $g$, according to the correspondence given by the interleaving maps $\alpha,\beta$, where we need to have pointwise convergence. Note that, the smaller $\varepsilon$, the bigger $\R-S^\varepsilon$.
\end{itemize}

Now we show that functions of the form $\varphi^{\Theta_p}$, with rather general assumptions, satisfy the conditions to apply \Cref{teo:stab_mapp}.
To formally introduce such assumptions, we need to introduce some further definitions which are standard in TDA.

\begin{defi}
    Let $\Vect_\mathbb{K}$ be the category of vector spaces, with linear maps, over some field $\mathbb{K}$. A functor $F: \R \rightarrow \Vect_\mathbb{K}$ is called a persistence module. A persistence module is called tame if $F(t)$ is always finite dimensional and there is a finite collection of real numbers $\{t_1<t_2<\ldots<t_n\}$ which satisfies:
\begin{itemize}
\item $F(t)=\emptyset$ for all $t<t_1$;
\item if $t,t'\in (t_i,t_{i+1})$, with $t<t'$, then $F(t<t')$ is bijective.
\end{itemize}
The values $\{t_1<t_2<\ldots<t_n\}$ are called critical values of the persistence module. Similarly, a continuous function $f:X\rightarrow \R$ is called $p-$tame ($p\in \N$) if, for all $i=1,\ldots,p$, the persistence module $t\mapsto H_i(f^{-1}(-\infty,t])$ is tame. The set $\Crit(f)$ collects all the critical values of the persistence modules $H_i(\X,\mathbb{K})$ for $i=0,\ldots,p$.
\end{defi}

Note that these definitions extend classical Morse Theory for smooth functions \cite{morse}.

Consider two $p$-tame functions $f,g:U\rightarrow \R$ such that $\parallel f-g \parallel_\infty \leq \varepsilon$.
Let $\X$ and $\Y$ be, respectively, the sublevel set filtrations of $f$ and $g$. By construction, for every $t\in \R$, we have the following commutative diagram:
\[
\begin{tikzcd}
X_t\ar[r]\ar[dr]&X_{t+\varepsilon}\ar[r]\ar[dr]&X_{t+2\varepsilon}\\
Y_t\ar[ur]\ar[r,]&Y_{t+\varepsilon}\ar[r]\ar[ur]&Y_{t+2\varepsilon},
\end{tikzcd}
\]
where all the maps are given by inclusions. Applying $\pi_0$ we obtain the $\varepsilon$-compatible maps:
\[
\alpha:\pi_0(\X)\rightarrow \mathcal{S}_\varepsilon\pi_0(\Y), \qquad \beta:\pi_0(\Y)\rightarrow \mathcal{S}_\varepsilon\pi_0(\X).
\]
Instead, applying $H_p(\cdot,\mathbb{K})$, we obtain:
\begin{equation}
\begin{tikzcd}\label{eq:homology}
H_p(X_t,\mathbb{K})\ar[r]\ar[dr]&H_p(X_{t+\varepsilon},\mathbb{K})\ar[r]\ar[dr]&H_p(X_{t+2\varepsilon},\mathbb{K})\\
H_p(Y_t,\mathbb{K})\ar[ur]\ar[r,]&H_p(Y_{t+\varepsilon},\mathbb{K})\ar[r]\ar[ur]&H_p(Y_{t+2\varepsilon},\mathbb{K}).
\end{tikzcd}
\end{equation}
Now, consider $S=\Crit(f)\cup \Crit(g)$,  the union of the critical values of all the persistence modules $H_i(\X,\mathbb{K})$, for $i=0,\ldots,p$.
Suppose that there exists $t$ such that $[t,t+\varepsilon]\cap S = \emptyset$. This implies that the horizontal arrows in \Cref{eq:homology} are all isomorphism. In particular,
\begin{equation}
 \dim(H_i(X_{t'},\mathbb{K}))=\dim(H_i(Y_{t'},\mathbb{K})),
\end{equation}
for all $t'\in [t,t+\varepsilon]$. Which implies that we can 
apply \Cref{teo:fun_stab} with $\varepsilon'=0$. Thus, we have the following result (whose proof follows immediately from these considerations).

\begin{prop}[Finite Stability]\label{prop:stab_betti}
    Given two $p-$tame functions $f,g:U\rightarrow \R$ such that:
    \[
    \dim(H_i(X_{t},\mathbb{K})),\dim(H_i(Y_{t},\mathbb{K}))\leq C,
    \]  
    for all $t\in\R$ and for all $i=0,\ldots,p$. Then, we have:
     \begin{equation}\label{eq:finitely_stab_p}
    d_E(\varphi_f^{\Theta_p},\varphi_g^{\Theta_p}) \leq (\#\Crit(f)+\#\Crit(g)) 2C\parallel f-g\parallel_\infty.
     \end{equation}
    In particular, if $p=0$, we obtain:
     \begin{equation}\label{eq:finitely_stab}
    d_E(\varphi_f^{\Theta_0},\varphi_g^{\Theta_0}) \leq (\#E_T+\#E_G) 2\parallel f-g\parallel_\infty.         
     \end{equation}
\end{prop}

To conclude, we comment on \Cref{eq:finitely_stab_p} and \Cref{eq:finitely_stab}.
In \cite{pegoraro2024finitely}, \Cref{eq:finitely_stab} is used to define a condition for a metric between merge trees which is named \emph{finite stability}. This is motivated by the fact that, if we have an uniform upper bound on the number of critical values of the considered merge trees, then, the metric gives a Lipschitz operator from functions (with the sup norm) to merge trees. 

Supported by the fact that the $1$-Wasserstein metric between persistence diagrams satisfies:
\[
W_1(D_f,D_g)\leq (\#D_f+\#D_g)\parallel f-g\parallel_\infty,
\]
in \cite{pegoraro2024finitely} it is argued that such stability properties may be preferable to universal stability properties, in which the operator from functions to topological representation is $1$-Lipschitz. Moreover, such comparison is framed in analogy to the bias-variance tradeoff in statistical modeling. 

Building on such considerations, we highlight that, if we represent the persistence module $H_i(\X,\mathbb{K})$ with the persistence diagram $D_f^i$, we have:
\[
\sum_{i=0}^p W_1(D^i_f,D^i_g)\leq (\#\Crit(f)+\#\Crit(g)) \parallel f-g\parallel_\infty.
\]
We can therefore view \Cref{eq:finitely_stab_p} as extending this analogy to the setting where all homological dimensions are considered simultaneously, analogously to what \cite[Theorem 3.18]{curry2021decorated} achieves for the bottleneck distance.

\section{Simulated Scenarios}
\label{sec:simulated_data}

Now we use two simulated data sets to put to work the frameworks defined in \Cref{sec:info_examples}. The algorithm employed to compute the metric is proposed in \cite{pegoraro2023edit}. 

The examples are basic, but suited to assert that dendrograms and the metric $d_E$ capture the information we designed them to grasp. In particular, since examples in  \Cref{sec:ex_dend} and \Cref{sec:ex_merge} already give insights into the role of the 
tree-structured information, we want to isolate and emphasize the key role of weight functions. We also deal with the problem of approximating the metric $d_E$ when the number of leaves in the tree structures in the data set is too big to be handled.
The examples presented concern hierarchical clustering dendrograms and dendrograms representing scalar fields.

In the implementations, dendrograms are always considered with a binary tree structure, obtained by adding negligible edges, that is edges $e$ with arbitrary small $d(\varphi(e),0)$, when the number of children of a vertex exceeds $2$.

\subsection{Pruning}\label{sec:pruning_TDA}
In this section we present a way of approximating the edit distance when the number of leaves of the involved tree structures is too high, taken from \cite{pegoraro2024functional}.

If one defines a proper weight function with values in an editable space $(E,d)$ coherently with the aim of the analysis, then the value $d(\varphi_T(e),0)$ can
be thought as the amount of information carried by the edge 
$e$. The bigger such value is, the more important that edge will be for the dendrogram.
In fact such edges are the ones most relevant in terms of $d_E$.
A sensible way to reduce the computational complexity of the metric $d_E$, losing as little information as possible, is therefore the following.
 Given \(\varepsilon >0\) and a dendrogram $(T,\varphi_T)$, define the following 1-step process:

\begin{itemize}
\item[$(\mathcal{P}_\varepsilon)$] Take a leaf $l$ such that $d(\varphi_{T}(l),0)$ is minimal among all leaves; if two or more leaves have minimal weight, choose \(l\) at random among them.  If   \(d(\varphi_{T}(l),0)<\varepsilon,\)
delete $l$ and ghost its parent if it becomes a degree $2$ vertex after removing $l$. 
\end{itemize}

We set $T_0=T$ and we apply operation (\({\cal P}_\varepsilon\)) to obtain $T_1$. On the result we apply again (\({\cal P}_\varepsilon\)) obtaining $T_2$ and, for \(n > 2,\) we proceed iteratively until we reach the fixed point of the sequence \(\{T_n\}\), which we call $P_\varepsilon(T)$. In this way we define the pruning operator $P_\varepsilon:\mathcal{T}\rightarrow \mathcal{T}$. Note that the fixed point is surely reached in a finite time since the number of leaves of each tree in the sequence is finite and non increasing along the sequence.

Lastly, if we define $\parallel T\parallel $ as $\parallel T\parallel =\sum_{e\in E_T}d(\varphi(e),0)$, we can quantify the (normalized) lost information with what we call \emph{pruning error} ($PE$): 
$(\parallel T\parallel -\parallel P_\varepsilon(T)\parallel )/\parallel T\parallel $.
%\begin{equation*}
%\frac{\mid \mid T\mid \mid -\mid \mid P_\varepsilon(T)\mid \mid }{\mid \mid T\mid \mid }
%\end{equation*}

\begin{figure}	
    \centering
    \begin{subfigure}[c]{0.6\textwidth}    	
    	\includegraphics[width = \textwidth]{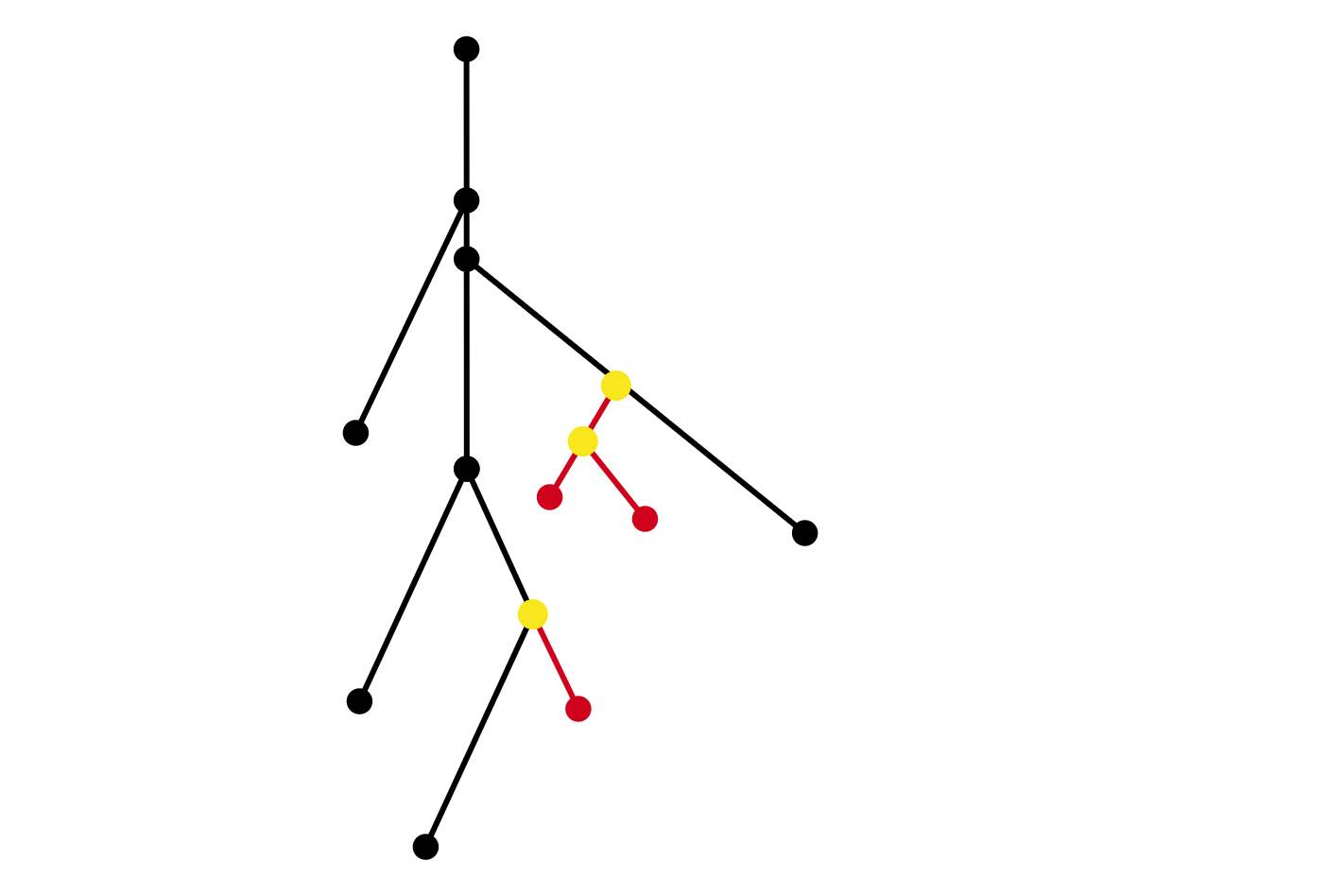}
	\end{subfigure}
%	\hspace{0.5 cm}
\caption{Pruning of a weighted tree: in red the deletions and in yellow the ghostings.}
\label{fig:pruning}
\end{figure}

\subsection{Hierarchical Clustering Dendrograms}\label{sec:clust_example}
We consider a data set of $30$ points clouds in $\mathbb{R}^2$, each with $150$ or $151$ points.
Point clouds are generated according to three different processes and are accordingly divided into three classes.
Each of the first $10$ point clouds is obtained by sampling 
independently two clusters of $75$ points respectively from normal distributions centered in $(5,0)$ and $(-5,0)$, both with $0.5 \cdot Id_{2	\times 2}$ covariance.  
Each of the subsequent $10$ point clouds is obtained by sampling independently $50$ points from each of the following Gaussian
distributions: one centered in $(5,0)$, one in $(-5,0)$
and one in $(-10,0)$. All with covariance $0.5 \cdot Id_{2	\times 2}$.
Lastly, to obtain each of the last $10$ point clouds,  we sample independently $150$ points as done for the first $10$ clouds, that is $75$ independent samples from a Gaussian centered $(5,0)$ and $75$ from one centered in $(-5,0)$, an then, to such samples, we add an outlier placed in $(-10,0)$.

\begin{figure}
    \begin{subfigure}[c]{0.24\textwidth}
    	\centering
    	\includegraphics[width = \textwidth]{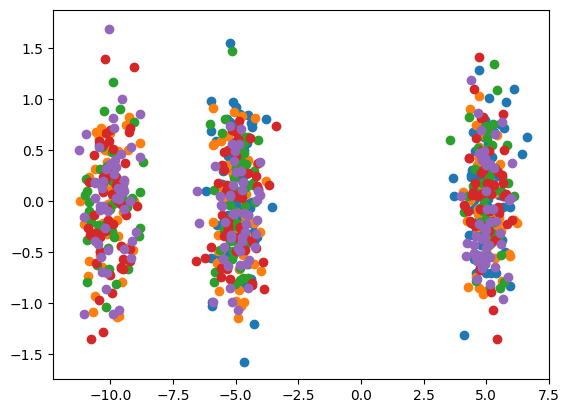}
    	\caption{Data from the second class.}
    	\label{fig:esempio_hier_second_class}
%    	\label{fig:cont_0_pt}
	\end{subfigure}
%	\hspace{0.5 cm}
    \begin{subfigure}[c]{0.24\textwidth}
    	\centering
	    \includegraphics[width = \textwidth]{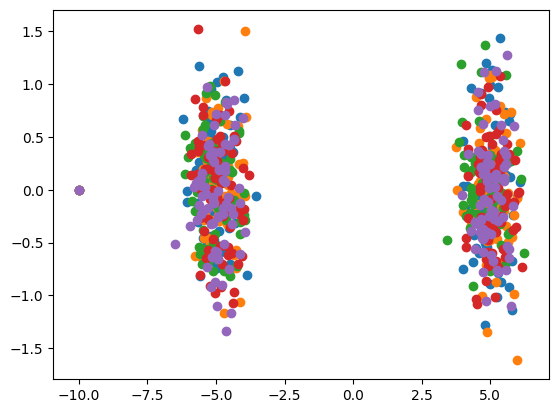}
    	\caption{Data from the third class.}
    	\label{fig:esempio_hier_third_class}
%    	\label{fig:cont_0_dend}
	\end{subfigure}			
    \begin{subfigure}[c]{0.24\textwidth}
    	\centering
    	\includegraphics[width = \textwidth]{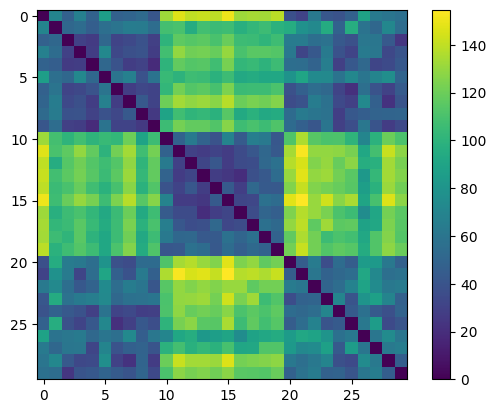}
    	\caption{Pairwise distances with dendrograms.}
    	\label{fig:esempio_hier_dendro}
%    	\label{fig:cont_1_pt}
	\end{subfigure}
%	\hspace{0.5 cm}
    \begin{subfigure}[c]{0.24\textwidth}
    	\centering
	    \includegraphics[width = \textwidth]{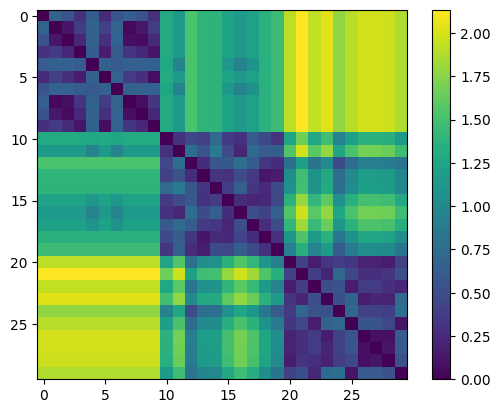}
    	\caption{Pairwise distances with PD.}
    	\label{fig:esempio_hier_PD}
%    	\label{fig:cont_1_dend}
		\end{subfigure}	
\caption{Data and pairwise distance matrices involved in the hierarchical clustering example.}
\label{fig:esempio_hier}
\end{figure}

Some clouds belonging to the second class and third classes are plotted respectively in  \Cref{fig:esempio_hier_second_class} and
 \Cref{fig:esempio_hier_third_class}.
We obtain dendrograms induced by the single linkage hiercarhical clsutering dendrograms, with the normalized cardinality functions induced by $\Theta_c$ and then resort to pruning because of the high number of leaves, but we still expect to be able to easily separate point clouds
belonging to the first and third classes (that is, with two major clusters) from clouds belonging to the second class, which feature three clusters, thanks to the cardinality information function defined in \Cref{sec:dendro_clus}. All dendrograms have been pruned with the same threshold, giving an average pruning error of $0.15$.

We can see in  \Cref{fig:esempio_hier_dendro} that this indeed the case. 
 It is also no surprise that persistence diagrams do not perform equally good in this classification task, as displayed in  \Cref{fig:esempio_hier_PD}. In fact PDs have no information about the importance of the cluster, making it impossible to properly recognize the similarity between data from the first and third class. They are, however, able to distinguish clouds belonging to class two from clouds belonging to class three since the persistence of the homology class associated to the leftmost cluster in clouds belonging to class two is smaller compared to what happens in clouds from the third class. 
The cluster centered in $(-10,0)$ and the one in $(-5,0)$ are in fact closer when the first one is a proper cloud, than when it is a cluster made by a single point.

\subsection{Dendrograms of Functions}\label{sec:merge_example}

This time our aim is to work with dendrograms obtained from functions, adding the (truncated) weight function induced by the Lebesgue measure of the sublevel sets $\Theta_{\mathcal{L}}$ and 
using them to discriminate between two classes in 
a functional data set.

We simulate the data set so that the discriminative information is contained in the size of the sublevel sets and not in the structure of the critical points.
To do so, we reproduce a situation which is very similar to the one shown by \cite{kmean_aling} for the Berkeley Growth Study data, where 
all the variability between groups in a classification task is explained by warping functions.
We fix a sine function defined over a compact $1D$ real interval (with the Lebesgue measure) and we apply to its domain $100$ random non linear warping functions belonging to two different, but balanced, groups. Warpings from the first group are more likely to obtain smaller sublevel sets, while in the second groups we should see larger sublevel sets and so \virgolette{bigger} weight functions defined on the edges. 
Note that, being the Lebesgue measure invariant with the translation of sets, any horizontal shifting of the functions would not change the distances between dendrograms.

The base interval is $I=[0,30]$ and the base function is $f(x)=\text{sin}(x)$. 
The warping functions are drawn in the following way.
Pick $N$ equispaced control points in $I$ and then we draw $N$ samples from a Gaussian distribution truncated to obtain only positive values. We thus have $x_1,\ldots x_N$ control points and $v_1,\ldots,v_N$ random positive numbers. Define $y_i := \sum_{j=1}^iv_j$.
The warping is then obtained interpolating with monotone cubic splines  the pairs $(x_i,y_i)$.
Being the analysis invariant to horizontal shifts in the functions, for all statistical units we  fix $x_0=y_0=0$ for visualization purposes.

\begin{figure}
	\centering
	\begin{subfigure}[c]{0.29\textwidth}
    	\centering
    	\includegraphics[width = \textwidth]{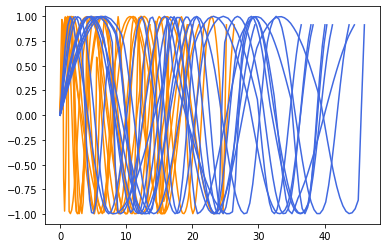}
    	\caption{Subset of functions colored by classes.}
    	\label{fig:esempio_merge_dati}
    \end{subfigure}
%    \hfill
	\begin{subfigure}[c]{0.29\textwidth}
		\centering
		\includegraphics[width = \textwidth]{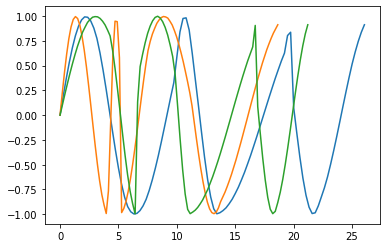}
		\caption{Few functions from the first class.}
	\end{subfigure}
	\begin{subfigure}[c]{0.29\textwidth}
    	\centering
    	\includegraphics[width = \textwidth]{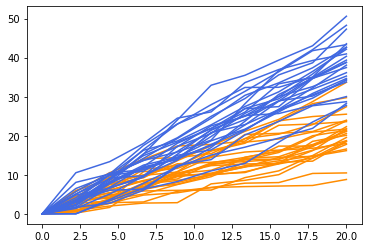}
    	\caption{Subset of warping functions colored by classes.}
    	\label{fig:esempio_merge_warp}
    \end{subfigure}
%    \hfill
    \begin{subfigure}[c]{0.29\textwidth}
    	\centering
    	\includegraphics[width = \textwidth]{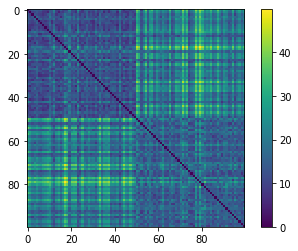}
    	\caption{Pairwise distances with dendrograms.}
    	\label{fig:esempio_merge_dist_dend}
    \end{subfigure}
    \begin{subfigure}[c]{0.29\textwidth}
    	\centering
    	\includegraphics[width = \textwidth]{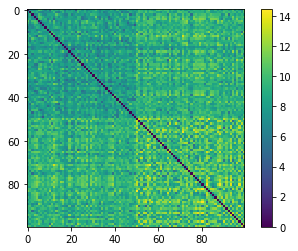}
    	\caption{Pairwise distances with $L_2$ metric.}
    	\label{fig:esempio_merge_dist_L2}
    \end{subfigure}
%    \hfill
    \begin{subfigure}[c]{0.29\textwidth}
    	\centering
    	\includegraphics[width = \textwidth]{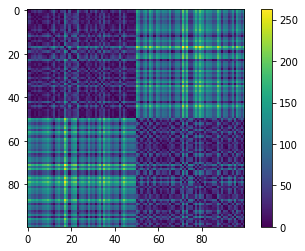}
    	\caption{Pairwise distances with $L_2$ metric on warping functions.}
    	\label{fig:esempio_merge_dist_warp}
    \end{subfigure}
    \begin{subfigure}[c]{0.29\textwidth}
    	\centering
    	\includegraphics[width = \textwidth]{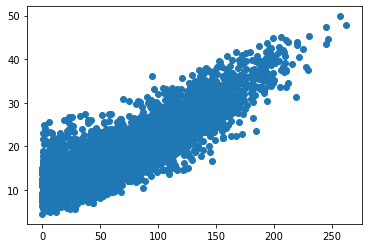}
    	\caption{Correlation between dendrograms and warping functions metric.}
    	\label{fig:esempio_merge_corr_dend}
    \end{subfigure}
%    \hfill
    \begin{subfigure}[c]{0.29\textwidth}
    	\centering
    	\includegraphics[width = \textwidth]{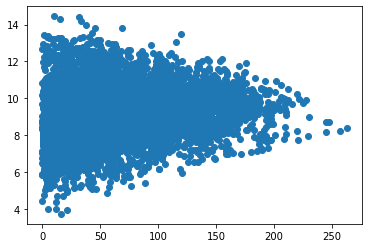}
    	\caption{Correlation between naive $L_2$ and warping functions metric.}
    	\label{fig:esempio_merge_corr_L2}
    \end{subfigure}

\caption{Overview of the example of  \Cref{sec:merge_example}.}
\label{fig:esempio_merge}
\end{figure}

The groups are discriminated by the parameters of the Gaussian distribution from which we sample the positive values $v_i$ to set up the warpings. For the first class we sample $N=10$ positive numbers from a truncated Gaussian with mean $3$ and standard deviation $2$; for the second the mean of the Gaussian is $5$ and the standard deviation is $2$.
Thus we obtain each of the first $50$ functions sampling $10$ values $v_i$ from the truncated Gaussian centered in $3$, building the warping function as explained in the previous lines, and then reparametrizing the sine function accordingly. The following $50$ functions are obtained with the same pipeline but employing a Gaussian centered in $5$. Note that, by construction, all the functions in the data set share the same merge tree. We truncate the functions induced by $\Theta_{\mathcal{L}}$ at height $1$.

 Examples of the warping functions can be seen in  \Cref{fig:esempio_merge_warp}; the resulting functions can be seen in  \Cref{fig:esempio_merge_dati}. 
The key point here is that we want to see if the dendrograms can retrieve the information contained in the warping functions. For this reason we compare the $L_2$ pairwise distances between such functions (see  \Cref{fig:esempio_merge_dist_warp}) and the pairwise distances obtained with dendrograms (see  \Cref{fig:esempio_merge_dist_dend}).
The visual inspection confirms the close relations between
the two sources of information. Moreover, if we vectorize the arrays given by the two matrices (considering only entries above the diagonal) and compute the Fisher correlation, we get a score of $0.85$ (see  \Cref{fig:esempio_merge_corr_dend}).
Instead, a naive approach with the $L_2$ metric applied directly to the data set would capture no information at all, as we can observe from  \Cref{fig:esempio_merge_dist_L2} and the Fisher correlation with the matrix obtained from warping functions is $0.15$ (see  \Cref{fig:esempio_merge_corr_L2}).

Note that, in general, the problem of finding warping functions to align functional data is deeply studied and with no easy  solution (see, for instance, the special issue of the Electronic Journal of Statistics dedicated to phase and amplitude variability - year 2014, volume 8 or \cite{sriva_SRV}) especially for non-linear warping of multidimensional or non-euclidean domains.  
Instead, dendrograms are less sensitive to such dimensionality issues, as they only 
arise in calculating the connected components and measure of the sublevel sets.

\section{Conclusions}
\label{sec:conclusions}

We develop a framework to work with functions defined on different merge trees. 
As motivated in the paper, we argue that these kinds of topological summaries can succeed in situations where persistence diagrams and merge trees alone are not effective. 
They also provide a great level of versatility because of the wide range of additional information that can be extracted from data. 

We introduce a metric structure that we argue is well-suited to the types of objects under consideration. This claim is supported both by qualitative insights and by formal stability results. These results focus on a setting of particular interest in TDA, where merge trees are used to decompose the homological information of a filtration across its path-connected components.
We demonstrate that the stability properties of our metric are analogous to those of the $1$-Wasserstein distance, when the ensemble of all persistence diagrams up to a fixed homological dimension is used as a topological summary.
Finally, we test the proposed framework beyond this purely topological setting, using simulated data to assess its practical effectiveness.

The main drawback of the framework is that the deformation between two functions is not guaranteed to always produce a function at the intermediate steps i.e. the metric space of local representation of functions is embedded in a bigger dendrograms space, but geodesics between points in general are not contained in this subspace.
In future works we would like to investigate when there are geodesics which remains in this subspace and if we can somehow modify this framework so that geodesics are always intrinsic. In case this does not hold true, it may limit the intrinsic statistical tools that can be defined in this space: should Frech\'et means exists, for instance, it is not guaranteed that they are functions.

\section*{Availability of data and materials}
The generative process of the simulated data is described in the paper. The code is available upon request.

\section*{Acknowledgments}
This work was carried out as part of my PhD Thesis, under the supervision of Professor Piercesare Secchi.
I also acknowledge the support of the Wallenberg AI, Autonomous Systems and Software Program (WASP), and of the SciLifeLab and Wallenberg National Program for Data-Driven Life Science (DDLS), which fund the project: Topological Data Analysis of Functional Genome to Find Covariation Signatures.

\appendix

\section*{Outline of the Appendix}

\Cref{sec:main_ideas} briefly motivates the use of merge trees over more traditional TDA's techniques with a pair of examples, introducing also some problems which can be solved by considering functions on merge trees.
In \Cref{sec:example_edit_fn} we showcase why all the machinery we set up to work with functions defined on merge trees does not work with persistence diagrams.  
\Cref{sec:proofs_TDA} contains the proofs of the results in the paper which are not included in the main text.

%%=============================================%%
%% For submissions to Nature Portfolio Journals %%
%% please use the heading ``Extended Data''.   %%
%%=============================================%%

%%=============================================================%%
%% Sample for another appendix section			       %%
%%=============================================================%%

%% \section{Example of another appendix section}\label{secA2}%
%% Appendices may be used for helpful, supporting or essential material that would otherwise 
%% clutter, break up or be distracting to the text. Appendices can consist of sections, figures, 
%% tables and equations etc.

\section{Why Use Trees}
\label{sec:main_ideas}

We want to give some motivation to propel the use of merge trees and functions defined on merge trees over persistence diagrams, in certain situations. We give only two brief examples since a similar topic is already tackled for instance in \cite{mergegrams, smith2022families, kanari2020trees, curry2024trees, curry2021decorated}.

\subsection{Point Clouds}\label{sec:ex_dend}

\begin{figure}[H]
%\begin{figure}[t]
    \begin{subfigure}[c]{0.3\textwidth}
    	\centering
    	\includegraphics[width = \textwidth]{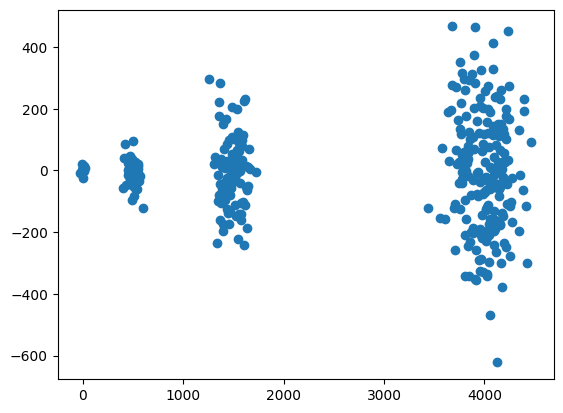}
    	\caption{First point cloud.}
    	\label{fig:cont_0_pt}
	\end{subfigure}
%	\hspace{0.5 cm}
    \begin{subfigure}[c]{0.3\textwidth}
    	\centering
	    \includegraphics[width = \textwidth]{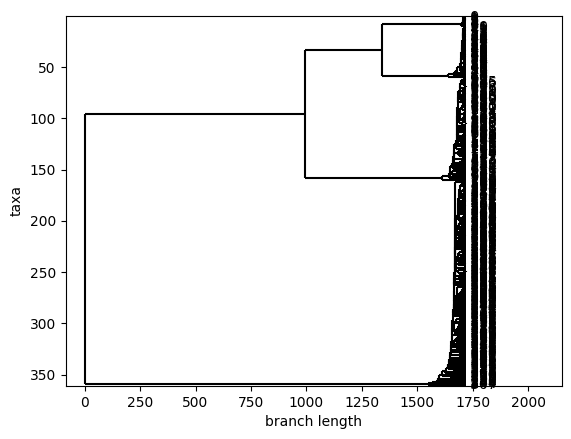}
    	\caption{First dendrogram.}
    	\label{fig:cont_0_dend}
	\end{subfigure}
%	\hspace{0.5 cm}
    \begin{subfigure}[c]{0.3\textwidth}
    	\centering
	    \includegraphics[width = \textwidth]{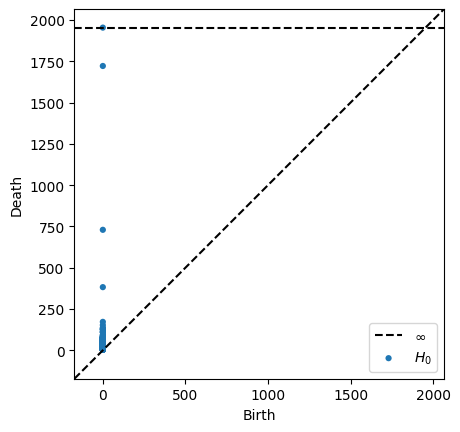}
    	\caption{First PD.}
    	\label{fig:cont_0_PD}
	\end{subfigure}
				
    \begin{subfigure}[c]{0.3\textwidth}
    	\centering
    	\includegraphics[width = \textwidth]{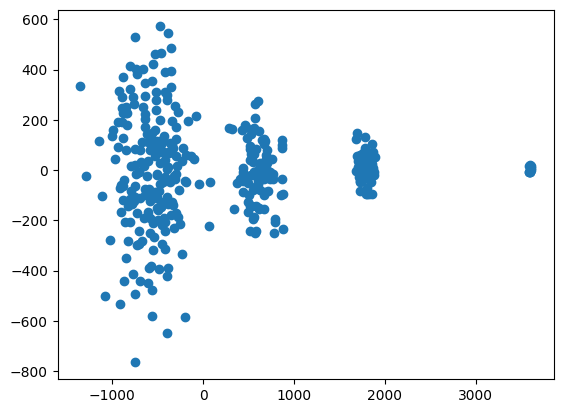}
    	\caption{Second point cloud.}
    	\label{fig:cont_1_pt}
	\end{subfigure}
%	\hspace{0.5 cm}
    \begin{subfigure}[c]{0.3\textwidth}
    	\centering
	    \includegraphics[width = \textwidth]{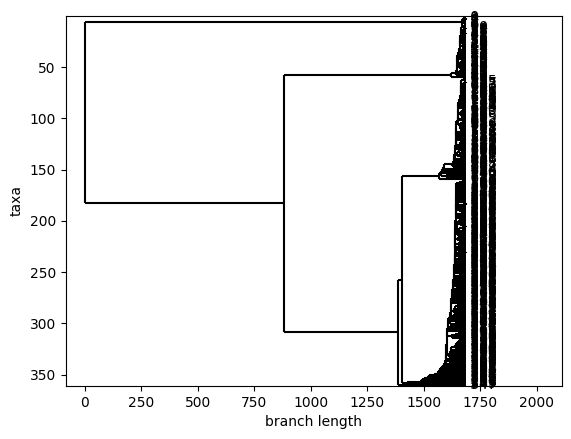}
    	\caption{Second dendrogram.}
    	\label{fig:cont_1_dend}
		\end{subfigure}
%	\hspace{0.5 cm}
    \begin{subfigure}[c]{0.3\textwidth}
    	\centering
	    \includegraphics[width = \textwidth]{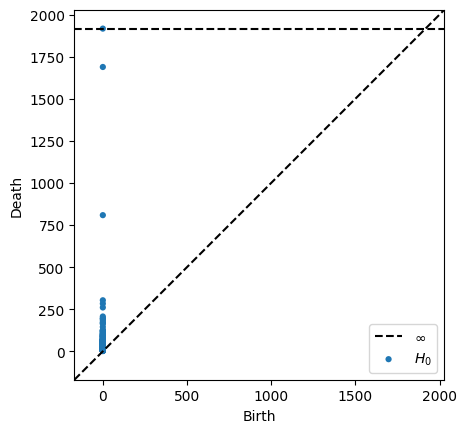}
    	\caption{Second PD.}
    	\label{fig:cont_1_PD}
	\end{subfigure}

\caption{Data clouds, hierarchical clustering dendrograms and PDs involved in the first example.}
\label{fig:example_dend}
\end{figure}

Given a point cloud $C=\{x_1,\ldots, x_n\}$ in 
$\mathbb{R}^n$ there are many ways in which one can build a family of simplicial complexes \cite{PH_survey} whose vertices are given by 
$C$ itself and whose sets of higher dimensional simplices get bigger and bigger.
A standard tool to do so is the Vietoris-Rips filtration of $C$ \cite{PH_survey}, as are 
$\alpha$ filtrations, $\check{C}$ech filtrations etc..

As we are interested only in path-connected components we restrict our attention to $0$ dimensional simplices (points) and $1$ dimensional simplices (edges). With such restrictions, many of the aforementioned filtrations become equivalent and amount to having a family of graphs
$\{C_t\}_{t\geq 0}$ such that the vertex set of $C_t$ is $C$ and the edge between $x_i$ and $x_j$ belongs to $C_t$ if and only if 
$d(x_i,x_j)<t$. Thus, the set of edges of $C_{t'}$ contains the set of edges of $C_t$, with $t\leq t'$; while the set of vertices is always $C$. 
Note, for instance, that the path-connected components of $C_t$ are equivalent to the ones of $X_{t/2}$ with $\X$ being the $\check{C}$ech filtration built in \Cref{sec:defi_intro}.
Along this filtration of graphs, the closest points become connected first and the farthest ones at last. It is thus reasonable to interpret the path-connected components of $C_t$ as clusters of the point cloud $C$. In order to choose the best \virgolette{resolution} 
to look at clusters, i.e. in order to choose $t$ and use $C_t$ to infer the clusters, statisticians look at 
the merge tree $\mathcal{M}(\pi_0(C_t)_{t\geq 0})$, which is called hierarchical clustering dendrogram. More precisely, $\mathcal{M}(\pi_0(C_t)_{t\geq 0})$ is the \emph{single linkage hierarchical clustering dendrogram}. Note that $\T$ is a regular abstract merge tree.

Suppose, instead, that we have the persistence diagram 
obtained from $\{\pi_0(C_t)\}_{t\in\mathbb{R}_{\geq 0}}$. Persistence diagrams are made of points in $\mathbb{R}^2$ whose coordinates $(b,d)$ represent the value of $t$ at which a certain path-connected component appears and the value of $t$ at which that component merges with a component which appeared before $b$. Each point in the point cloud is associated to a path-connected component but, in general, we have no way to distinguish between points of the diagram associated to path-connected components which are proper clusters and points of the diagrams associated to outliers.

Now, consider the single linkage dendrograms and the zero dimensional PDs obtained from point clouds as in  \Cref{fig:example_dend}. The persistence diagrams (in \Cref{fig:cont_0_PD} and \Cref{fig:cont_1_PD}) are very similar, in fact they simply record that there are four major clusters which merge at similar times across the Vietoris-Rips filtrations of the two point clouds. The hierarchical dendrograms, instead, are clearly very different since they show that in the first case (\Cref{fig:cont_0_pt}, \Cref{fig:cont_0_dend}, \Cref{fig:cont_0_PD}) the cluster with most points is the one which is more separated from the others in the point cloud; while in the second case (\Cref{fig:cont_1_pt}, \Cref{fig:cont_1_dend}, \Cref{fig:cont_1_PD}) the two bigger clusters are the first that get merged and the farthest cluster of points on the right could be considered as made by outliers. In many applications it would be important to distinguish between these two scenarios, since the two main clusters get merged at very different heights on the respective dendrograms.

These observations are formalized in 
\cite{curry2024trees}, with the introduction of the \emph{tree realization number} which is a combinatorial description of how many merge trees share a particular persistence diagram. With hierarchical clustering dendrograms with $n$ leaves, such number is $n!$: all leaves are born at height $0$, and so, at the first merging point, each of the $n$ leaves can merge with any of the $n-1$ remaining ones. At the following merging step we have $n-1$ clusters and each one of them can merge with the other $n-2$ etc..

\subsection{Real Valued Functions}\label{sec:ex_merge}

\begin{figure}
    \begin{subfigure}[c]{0.47\textwidth}
    	\centering
    	\includegraphics[width = \textwidth]{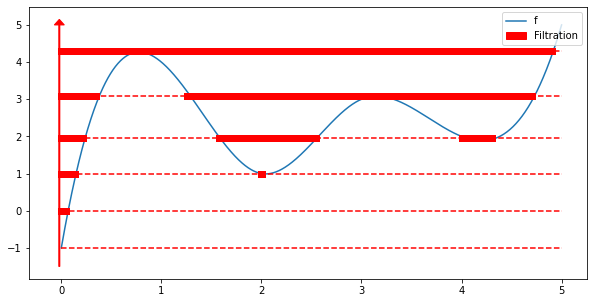}
    	\caption{Sublevel sets of a function}
    	\label{fig:sublvl_TDA}
	\end{subfigure}
%	\hspace{2 cm}
    \begin{subfigure}[c]{0.47\textwidth}
    	\centering
	    \includegraphics[width = \textwidth]{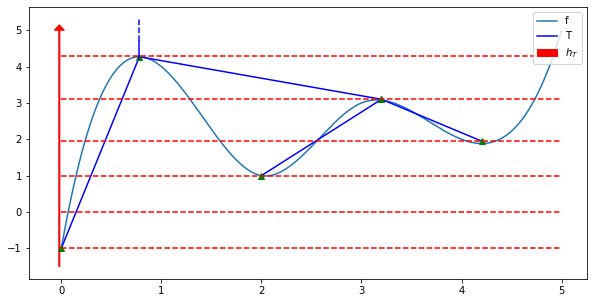}
    	\caption{A function with its associated merge tree.}
    	\label{fig:func_tree_TDA}
		\end{subfigure}
\caption{Merge Trees of Functions}
\end{figure}

Given a continuous function $f:[a,b]\rightarrow \mathbb{R}$ we can extract the merge tree $\mathcal{M}(\T)$, with $\X$ being the sublevel set filtration (see \Cref{sec:defi_intro} and \Cref{sec:easy_example}): we obtain a merge tree that tracks the evolution of the path-connected components of the sublevel sets $f^{-1}((-\infty,t])$. For a visual example see  \Cref{fig:func_tree_TDA}. Moreover, \cite{pegoraro2024functional} shows that $\T$ is a regular merge tree.

\begin{figure}[H]
	\centering
	\begin{subfigure}[t]{0.32\textwidth}
    	\centering
    	\includegraphics[width = \textwidth]{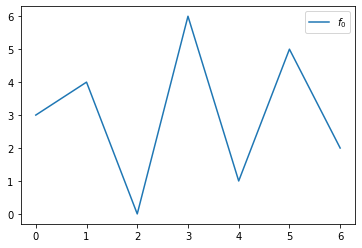}
    	\caption{A function $f_0$.}
    \end{subfigure}
	\begin{subfigure}[t]{0.32\textwidth}
    	\centering
    	\includegraphics[width = \textwidth]{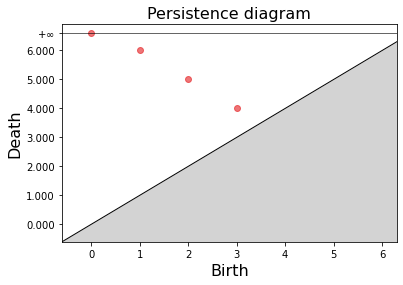}
    	\caption{The $0$-dimensional persistence diagram of $f_0$.}
    \end{subfigure}
	\begin{subfigure}[t]{0.32\textwidth}
    	\centering
    	\includegraphics[width = \textwidth]{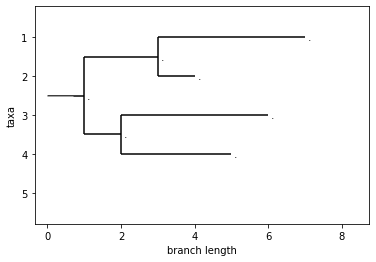}
    	\caption{The merge tree of $f_0$.}
    \end{subfigure}
    
	\begin{subfigure}[t]{0.32\textwidth}
    	\centering
    	\includegraphics[width = \textwidth]{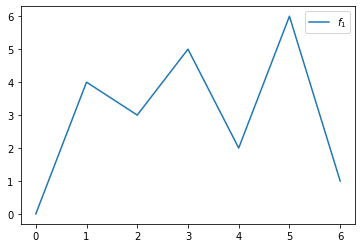}
    	\caption{A function $f_1$.}
    \end{subfigure}
	\begin{subfigure}[t]{0.32\textwidth}
    	\centering
    	\includegraphics[width = \textwidth]{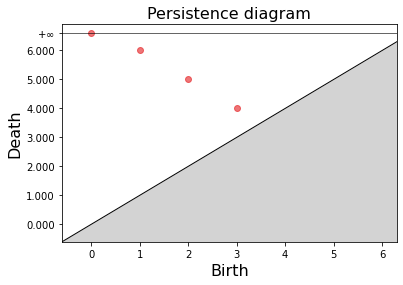}
    	\caption{The $0$-dimensional persistence diagram of $f_1$.}
    \end{subfigure}
	\begin{subfigure}[t]{0.32\textwidth}
    	\centering
    	\includegraphics[width = \textwidth]{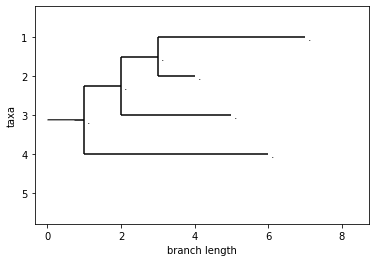}
    	\caption{The merge tree of $f_1$.}
    \end{subfigure}

	\begin{subfigure}[t]{0.32\textwidth}
    	\centering
    	\includegraphics[width = \textwidth]{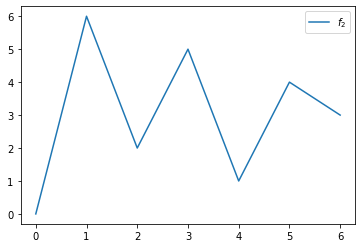}
    	\caption{A function $f_2$.}
    \end{subfigure}
	\begin{subfigure}[t]{0.32\textwidth}
    	\centering
    	\includegraphics[width = \textwidth]{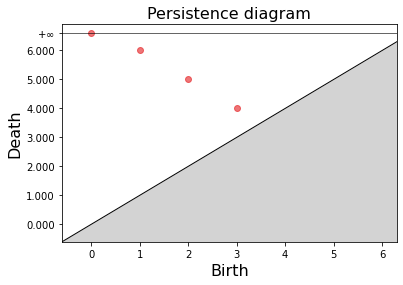}
    	\caption{The $0$-dimensional persistence diagram of $f_2$.}
    \end{subfigure}
	\begin{subfigure}[t]{0.32\textwidth}
    	\centering
    	\includegraphics[width = \textwidth]{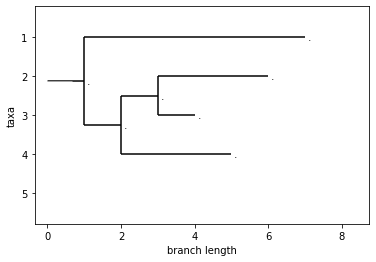}
    	\caption{The merge tree of $f_2$.}
    \end{subfigure}
    
	\begin{subfigure}[t]{0.32\textwidth}
    	\centering
    	\includegraphics[width = \textwidth]{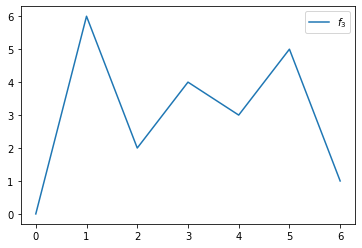}
    	\caption{A function $f_3$.}
    \end{subfigure}
	\begin{subfigure}[t]{0.32\textwidth}
    	\centering
    	\includegraphics[width = \textwidth]{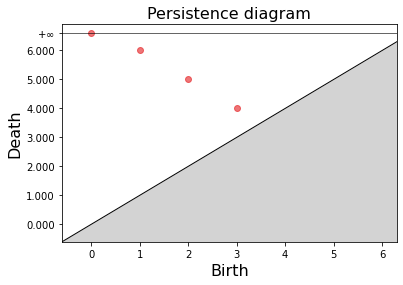}
    	\caption{The $0$-dimensional persistence diagram of $f_3$.}
    \end{subfigure}
	\begin{subfigure}[t]{0.32\textwidth}
    	\centering
    	\includegraphics[width = \textwidth]{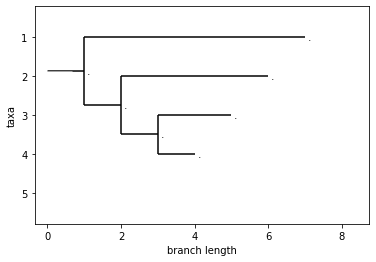}
    	\caption{The merge tree of $f_3$.}
    \end{subfigure}
    
\caption{We compare four functions; they are all associated to the same $PD$ but to different merge trees. Functions are displayed in the first column and on each row we have on the centre the associated $PD$ and on the right the merge tree.}
\label{fig:PD_vs_merge_TDA}
\end{figure}

We use this example to point out two facts.
First PDs may not be able to distinguish functions one may wish to distinguish, as made clear by  \Cref{fig:PD_vs_merge_TDA}.
Second, Proposition 1 of \cite{pegoraro2024functional} 
states that if one changes the parametrization of a function by means of homeomorphisms, then, both the associated merge tree and persistence diagram do not change.
A consequence of such result is that one can shrink or spread the domain of the function $f:[a,b]\rightarrow \mathbb{R}$ with reasonably regular functions, without changing its merge tree (and PD). There are cases in which such property may be useful but surely there are times when one may want to distinguish if an oscillation lasted for a time interval of $10^{-5}$ or $10^5$.
The measure related function defined in \Cref{sec:info_examples} can solve this issue.

\section{Functions on Merge Trees vs Functions on PDs}
\label{sec:example_edit_fn}

\begin{figure}
	\begin{subfigure}[c]{0.47\textwidth}
    	\centering
    	\includegraphics[ width = \textwidth]{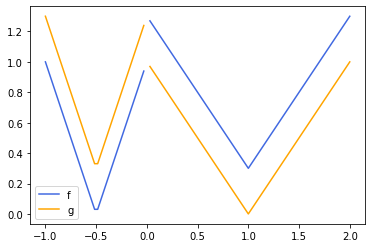}
    	\caption{The functions $f$ and $g$ in \Cref{sec:example_edit_fn}; with $\varepsilon = 0.3$.}
    	\label{fig:fun_functions}
    \end{subfigure}
    	\centering
    \begin{subfigure}[c]{0.47\textwidth}
    	\includegraphics[ width = \textwidth]{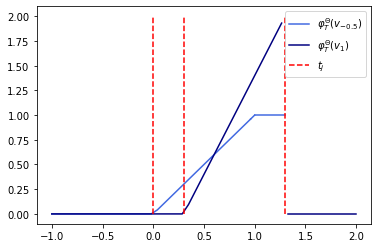}
    	\caption{In the context of the example in \Cref{sec:example_edit_fn}, we plot $\varphi^{\Theta_\mathcal{L}}_T(v_1)$ and $\varphi^{\Theta_\mathcal{L}}_T(v_{-0.5})$. The dotted lines represent critical values.}
    	\label{fig:fun_mult}
    \end{subfigure}
   
\end{figure}

We make one example which shows what could happen if we try to define functions on PDs in the same way we do for merge trees. In particular, the elder rule, via the instability of the persistence pairs, makes it very difficult to add pieces of information to persistence diagrams in a stable way.

Consider the following functions, plotted in \Cref{fig:fun_functions}, defined on $[-1,2]$:

\begin{align*}
f(x)=\mid x-1\mid +\varepsilon &\text{ if } x\geq 0 \\
 f(x)=\mid 2x-1\mid &\text{	if } x< 0
\end{align*}
and 
\begin{align*}
g(x)=\mid x-1\mid &\text{ if } x\geq 0 \\
 g(x)=\mid 2x-1\mid +\varepsilon &\text{ if } x< 0
\end{align*}
for a fixed $\varepsilon>0$.
 
Let $(T,h_T)$ and $(T',h_{T'})$ be the merge trees associated to the sublevel set filtrations of $f$ and $g$; moreover let $\varphi^{\Theta_\mathcal{L}}_T$ and $\varphi^{\Theta_\mathcal{L}}_{T'}$ the two respective local representations of the induced functions with $\Theta_\mathcal{L}$ being the Lebesgue measure on $\mathbb{R}$. 
Note that $\parallel f-g\parallel_\infty = \varepsilon$.
The local minima of the functions are the points $\{-0.5,1\}$, with $f(-0.5)=0$, $f(1)=\varepsilon$, $g(-0.5)=\varepsilon$ and $g(1)=0$. 
Thus the merge trees have isomorphic tree structures: we represent $T$ with the vertex set $\{v_{-0.5},v_{1},v_{0},r_T\}$ and edges $\{(v_{-0.5},v_0),(v_{1},v_0),(v_0,r_{T'})\}$; and $T'$ with vertices $\{v_{-0.5},v_{1},v_{0},r_{T'}\}$ and edges $\{(v_{-0.5},v_{0}),(v_{1},v_{0}),(v_0,r_{T'})\}$. The height functions are the following: $h_T(v_{-0.5})=0$, $h_{T'}(v_{-0.5})=\varepsilon$, $h_{T}(v_{1})=\varepsilon$, $h_{T'}(v_1)=0$ and 
$h_T(v_0)=h_{T'}(v_0)=1+\varepsilon$.
 
Having truncated both functions at height $1+\varepsilon$, the weight functions (see \Cref{fig:fun_mult}) are given by: $\varphi^{\Theta_\mathcal{L}}_T(v_{-0.5})(t)=t\chi_{[0,1)}+\chi_{[1,1+\varepsilon)} $, 
$\varphi^{\Theta_\mathcal{L}}_T(v_{1})(t)=2(t-\varepsilon)\chi_{[\varepsilon,1+\varepsilon)}$ and $\varphi^{\Theta_\mathcal{L}}_{T'}(v_{-0.5})(t)=(t-\varepsilon)\chi_{[\varepsilon,1+\varepsilon)}$ and $\varphi^{\Theta_\mathcal{L}}_{T'}(v_{1})(t)=2t \chi_{[0,1)}+ 2\chi_{[1,1+\varepsilon)}$. 

The zero-dimensional persistence diagram associated to $f$ (we name it $PD_0(f)$) is given by a point with coordinates $(0,+\infty)$, associated to the connected component $[-t/2-0.5,t/2-0.5]$ which is born at $t=0$, and the point $(\varepsilon,1+\varepsilon)$, associated to the component $[1-(t-\varepsilon),1+(t-\varepsilon)]$, born at level $t=\varepsilon$ and \virgolette{dying} at level $t=1+\varepsilon$, due to the elder rule, since it merges an older component, being the other component born at a lower level. 

For the function $g$, the persistence diagram $PD_0(g)$ is made by the same points, but the situation is in some sense \virgolette{reversed}. In fact, 
the point $(0,+\infty)$ is associated to the connected component 
\virgolette{centered} in $1$, which is $[1-t,1+t]$, and the point 
$(\varepsilon,1+\varepsilon)$, is associated to the component
\virgolette{centered} in $0.5$, that is
 $[-(t-\varepsilon)/2-0.5,(t+\varepsilon)/2-0.5]$.

The consequence of this change in the associations between points and the components originating the points of the diagrams is that the information regarding the two components, end up being associated to very different spatial locations in the two diagrams: $(0,+\infty)$ and $(\varepsilon,1+\varepsilon)$. And this holds for every $\varepsilon>0$. Thus it seems very hard to design a way to \virgolette{enrich} $PD_0(f)$ and $PD_0(g)$ with additional information,  originating the \virgolette{enriched diagrams} $D_f$ and $D_g$, respectively, and design a suitable metric $d$, so that $d(D_f ,D_g )\rightarrow 0$ as $\varepsilon\rightarrow 0$.

Instead, if we consider the edit path $\gamma$ which shrinks
$v_{-0.5}\rightarrow v_{-0.5}$
and 
$v_{1}\rightarrow v_{1}$
we have $d_E((T,\varphi^{\Theta_\mathcal{L}}_T),(T',\varphi^{\Theta_\mathcal{L}}_{T'}))\leq cost(\gamma)=3\varepsilon$.

\section{Proofs}\label{sec:proofs_TDA}

\bigskip\noindent
\textbf{\Cref{prop:regular}.}
For every abstract merge tree $\pi_0(\X)$ there is a unique (up to isomorphism) abstract merge tree $R(\pi_0(\X))$ such that:
\begin{enumerate}
\item $\pi_0(\X)\cong_{a.e.} R(\pi_0(\X))$;
\item $R(\pi_0(\X))$ is regular.
\end{enumerate}

\begin{proof}    
Let $\pi_0(\X)$ be an abstract merge tree with critical values $t_1<\ldots<t_n$. Suppose that at $t_j$ changes happen across the critical value. Then we can fix $\varepsilon>0$ such that $t_j+\varepsilon< t_{j+1}$ and define $\X'$ with $X'_t=X_t$ for all $t\neq t_j$ and $X'_{t_j}=X_{t_j+\varepsilon}$. Now we need to define the $\X'$ on maps:
\begin{itemize}
\item if $t=t_j$ and $t<t'\leq t_j+\varepsilon$,  $X_{t<t'}' = (X_{t'\leq t_j+\varepsilon})^{-1}$ which is well defined as $X_{t'\leq t_j+\varepsilon}$ is an isomorphism;
\item if $t'=t_j$,  $X_{t<t'}' = X_{t\leq t_j+\varepsilon}$ which is well defined as $X_{t'\leq t_j+\varepsilon}$ is an isomorphism;
\item otherwise $X'_{t<t'}=X_{t<t'}$.
\end{itemize}

We need to check that $\X'$ is a regular abstract merge tree.
First we have: 
\[
X'_{t,t_j}\circ X'_{t_j,t'}= X_{t\leq t_j+\varepsilon} \circ (X_{t'\leq t_j+\varepsilon})^{-1}=X_{t,t'}=X_{t,t'}'
\] 
if $t'\leq t_j+\varepsilon$, otherwise 
\[
X'_{t,t_j}\circ X'_{t_j,t'}= X_{t\leq t_j+\varepsilon} \circ X_{t_j+\varepsilon\leq t'}= X_{t,t'}=X_{t,t'}'.
\]

The filtration $\X'$ is regular at $t_j$ by construction  as  $X'_{t_j}= X_{t_j+\varepsilon}\cong X'_{t'}$ for $t'\in [t_j,t_j+\varepsilon]$. Always by construction, it is a.e. isomorphic to $\X$: the natural transformation $\varphi:\X\rightarrow \X'$ is given by $\varphi_t = Id:X_t\rightarrow X'_t$ for $t\neq t_j$ and, in fact, it is defined on $\mathbb{R}-\{t_1,\ldots,t_n\}$.

If $t_j$ is the only critical value at which changes in $\X$ happen across the value we are done, otherwise consider $t_k$ such that 
changes in $\X$ happen across $t_k$. The same, by construction, holds also for $\X'$. Thus we can recursively apply the steps proposed up to now on $\X'$ until we obtain an abstract merge tree $R(\T)$ which is regular. This is reached in a finite number of steps since the critical values are a finite set.

Uniqueness (up to isomorphism) follows easily.
\end{proof}

\bigskip\noindent
\textbf{\Cref{prop:equivalence}.}
The following hold:
\begin{enumerate}
\item we can associate to a regular abstract merge tree $R(\T)$, a merge tree without degree $2$ vertices $\mathcal{M}(R(\T))$;
\item we can associate to a merge tree $(T,h_T)$, a regular abstract merge tree $\mathcal{F}((T,h_T))$. Moreover, we have $\mathcal{M}(\mathcal{F}((T,h_T)))\cong_2 (T,h_T)$ and $\mathcal{F}(\mathcal{M}(R(\T))\cong_{a.e.} \T$;
\item given two abstract merge trees $\X$ and $\Y$, $\mathcal{M}(R(\T))\cong \mathcal{M}(R(\G))$ if and only if $\T\cong_{a.e} \G$.
\item given two merge trees $(T,h_T)$ and $(T',h_{T'})$, we have $\mathcal{F}((T,h_T))\cong \mathcal{F}((T,h_T))$ if and only if $(T,h_T)\cong_2 (T',h_{T'})$.
\end{enumerate}

\begin{proof}
\begin{enumerate}
\item WLOG suppose $\T \cong R(\T)$; we build the merge tree $\mathcal{M}(\T)=(T,h_T)$ along the following rules in a recursive fashion  starting from an empty set of vertices $V_T$ and an empty set of edges $E_T$. We simultaneously add points and edges to $T$ and define $h_T$ on the newly added vertices. Let $\{t_i\}_{i=1}^n$ be the critical set of $\T$ and let $\pi_0(X_{t}):=a_{t}:=\{a^{t}_1,\ldots,a^{t}_{n_t}\}$. Call $\psi_t^{t'}:=\pi_0(X_{t\leq t'})$.
Lastly, from now on, we indicate with $\#C$ the cardinality of a finite set $C.$

Considering in increasing order the critical values:

\begin{itemize}
\item for the critical value $t_1$ add to $V_T$ a leaf $a_{t_1}^k$, with height $t_1$, for every element $a_{t_1}^k\in a_{t_1}$;
\item for $t_i$ with $i>1$, for every $a_{t_{i}}^k\in a_{t_i}$ such that $a_{t_{i}}^k\notin \text{Im}(\psi_{t_{i-1}}^{t_{i}}))$, add to $V_T$ a leaf $a_{t_i}^k$ with height $t_{i}$;
\item for $t_i$ with $i>1$, if $a_{t_i}^k=\psi_{t_{i-1}}^{t_{i}}(a_{t_{i-1}}^s)=\psi_{t_{i-1}}^{t_{i}}(a_{t_{i-1}}^r)$, with $a_{t_{i-1}}^s$ and $a_{t_{i-1}}^r$ distinct basis elements in $a_{t_{i-1}}$, add a vertex $a_{t_i}^k$ with height $t_{i}$, and add edges so that the previously added vertices 
\[
v = \arg\max \{h_T(v')\mid v' \in V_T \text{ s.t. }\psi_{t_{v'}}^{t_i}(v')=a_{t_i}^k  \}
\]
and 
\[
w = \arg\max \{h_T(w')\mid w' \in V_T \text{ s.t. }\psi_{t_{w'}}^{t_i}(w')=a_{t_i}^k \}
\]
 connect with the newly added vertex $a_{t_i}^k$.  
\end{itemize}

The last merging happens at height $t_n$ and, by construction, at height $t_n$ there is only one point, which is the root of the tree structure.

These rules define a tree structure with a monotone increasing height function $h_T$. In fact, edges are induced by maps $\psi_t^{t'}$ with $t<t'$ and thus we can have no cycles and the function $h_T$ must be increasing. Moreover, we have $\psi_t^{t_n}(a_i^t)=a_1^{t_n}$ for every $i$ and $t<t_n$ and thus the graph is path-connected.   

\item Now we start from a merge tree $(T,h_T)$ and build an abstract merge tree $\T$ such that $\mathcal{M}(\T)\cong (T,h_T)$. 
%Clearly, 
%for every regular abstract merge tree $\T$ if the tree structure of $\mathcal{M}(\T)$ has an other two vertex, then it must be the root. Moreover, at every critical point the are path components merging with each other or arising, and so \MPnote{root} $r_T$ cannot be of degree $1$.

To build the abstract merge tree, the idea is that we would like to \virgolette{cut} $(T,h_T)$ at every height $t$ and take as many elements in the set of path-connected components as the edges met by the cut.

Let $\{t_1,\ldots, t_n\}$ be the ordered image of $h_T$ in $\mathbb{R}$.

Consider the sets $v_{t_j}=\{v_i^{t_j}\}_{i=1,\ldots, n_{t_j}}=h_T^{-1}(t_j)$. We use the notation $\mathcal{F}((T,h_T))_t:=a_{t}:=\{a^{t}_1,\ldots,a^{t}_{n_t}\}$.
We define $a_{t_1}= v_{t_1}$. 
For every $\varepsilon>0$ such that $t_2-t_1>\varepsilon$, we set $a_{t_1+\varepsilon}=a_{t_1}$ and consequently $\psi_{t}^{t'}=Id$ for every $[t,t']\subset [t_1,t_2)$. Now we build $a_{t_2}$ starting from $a_{t_1}$, using $v_{t_2}$. We need to consider $v_i^{t_2}\in v_{t_2}$. There are the following possibilities:
\begin{itemize}
\item if $v_i^{t_2}$ is a leaf, then we add $v_i^{t_2}$ to $a_{t_1}$;
\item if $v_i^{t_2}$ is an internal vertex with $\#\text{child}(v_i^{t_2})>1$ - i.e. a merging point, we add $v_i^{t_2}$ to $a_{t_1}$ and then remove $\text{child}(v_i^{t_2})=\{v\in V_T \mid v \text{ is a children of }v_i^{t_2}\}$. Note that, by construction, $\text{child}(v_i^{t_2})\subset a_{t_1}$, and, by hypothesis, $\#\text{child}(v_i^{t_2})>1$;
\item if $v_i^{t_2}$ is an internal vertex with $\#\text{child}(v_i^{t_2})=1$ - i.e. a degree $2$ vertex, we don't do anything.  
\end{itemize}
We obtain $a_{t_2}$ from $a_{t_1}$  by applying recursively these rules for every $v_i^{t_2}\in v_{t_2}$. The map $\psi_{t}^{t_2}$, for $t\in [t_1,t_2)$ is then defined by setting $\psi_{t}^{t_2}(a_i^{t_1})=v_i^{t_2}$ if 
$a_i^{t_1}\in \text{child}(v_i^{t_2})$ and $\psi_{t}^{t_2}(a_i^{t_1})=a_i^{t_1}$ otherwise.
To define $a_t$ for $t>t_2$ we recursively repeat for every critical value $t_i$ (in increasing order) the steps of defining $a_{t_i+\varepsilon}$ equal to $a_{t_i}$ for small $\varepsilon>0$, and then, adjusting $a_{t_i}$ according to the tree structure to obtain $a_{i+1}$ and $\psi_{t_i}^{t_{i+1}}$, using $v_{t_{i+1}}$, as explained above.
When reaching $t_n$ we have $v_{t_n}=\{v^{t_n}_1\}$ and we set $a_{t}=v_{t_n}$ for every $t\geq t_n$.

We call this persistent set $\mathcal{F}((T,h_T))$. Note that, by construction:
\begin{itemize}
\item for every $v\in V_T$ we have $v\in a_t$ for $t\in [h_T(v),h_T(\text{parent}(v))$;
\item $\mathcal{F}((T,h_T))$ is regular;
\item $\mathcal{F}((T,h_T))$ is independent from degree $2$ vertices of $(T,h_T)$;
\item $\mathcal{F}((T,h_T))$ is an abstract merge tree.
\end{itemize}

Now we need to check that $(T',h_{T'})=\mathcal{M}(\mathcal{F}((T,h_T)))\cong_2 (T,h_T)$. 
WLOG we suppose $(T,h_T)$ is without degree $2$ vertices and prove 
$(T',h_{T'})=\mathcal{M}(\mathcal{F}((T,h_T)))\cong (T,h_T)$.
Let $\T= \mathcal{F}((T,h_T))$.

As before, for notational convenience, we set $a_t:= \pi_0(X_t)$ and $\psi_{t}^{t'}:=\pi_0(X_{t\leq t'})$. By construction, 
$a_t\subset V_T$ for every $t$. Which implies $V_{T'}\subset V_T$.

Consider now $a_{t_i}$ with $t_i$ critical value. To build $\T$ elements $a_j^{t_{i-1}},a_k^{t_{i-1}}\in a_{t_{i-1}}$ are replaced by $v$ in $a_{t_{i}}$ if and only if they merge with $v$ in the merge tree $(T,h_T)$: $(a_j^{t_{i-1}},v),(a_k^{t_{i-1}},v)$, with $h_T(v)=t_i$. The maps $\psi_{t_i}^{t_{i}}:a_{t_{i-1}}\rightarrow a_{t_{i-1}}$ are defined accordingly to represent that merging mapping $a_j^{t_{i-1}} \mapsto v$ and $a_k^{t_{i-1}} \mapsto v$.  So an element $v'$ stays in $a_{t}$ until the edge $(v',parent(v'))$ meets another edge in $(T,h_T)$, and then is replaces by $parent(v')$. As a consequence, we have 
$a_j^{t_{i-1}},a_k^{t_{i-1}},v\in V_{T'}$ and 
$(a_j^{t_{i-1}},v),(a_k^{t_{i-1}},v)\in E_{T'}$.

Since $(T,h_T)$ has no degree $2$ vertices then 1) $V_T = \bigcup_{i=1,\ldots,n}  a_{t_i}$ 2)  $V_{T}=V_{T'}$ 3) $id:V_T\rightarrow V_{T'}$ is an isomorphism of merge trees.

Now we consider $\T$ regular abstract merge tree and prove $\mathcal{F}(\mathcal{M}(R(\T))\cong \T$. 
Consider $t_i$ critical value, $\varepsilon >0$ such that $t_{i-1}<t_i-\varepsilon$ and let $v_{t_i}=\{v \in \pi_0(X_{t_i})\mid \#\pi_0(X_{t_i-\varepsilon <t_i})^{-1}(v)\neq 1\}$.
By construction, $v_{t_i}\subset V_T$, for every $t_i$ critical value, with $(T,h_T)=\mathcal{M}(R(\T)$.

For every $v\in \pi_0(X_t)$, for any $t\in \mathbb{R}$ there is $v_j^{t_i}\in v_{t_i}$ for some $t_i$, such that $\pi_0(X_{t_i\leq t})(v_j^{t_i})=v$. Moreover the following element is 
well defined:
\[
s(v):=\max \{ w \in v_{t_i}, t_i \text{ critical value } \mid \pi_0(X_{t_i\leq t})(w)=v\} 
\]

By construction we have $v = \pi_0(X_{t_i\leq t})(s(v))$.

Let $\G=\mathcal{F}(\mathcal{M}(R(\T))$.
Define $\alpha_t:\pi_0(X_t)\rightarrow \pi_0(Y_t)$ given by $v = \pi_0(X_{t_i\leq t})(s(v)) \mapsto s(v)$.
It is an isomorphism of abstract merge trees.

\item if $\T\cong_{a.e.} \G$, then $R(\T)\cong R(\G)$ and then the merge trees $\mathcal{M}(R(\T)$ and $\mathcal{M}(R(\G)$ differ just by a change in the names of the vertices. If $\mathcal{M}(R(\T)\cong \mathcal{M}(R(\G)$ then 
$\mathcal{F}(\mathcal{M}(R(\T))\cong \mathcal{F}(\mathcal{M}(R(\G))\cong R(\T)\cong R(\G)$.

\item the proof is analogous to the one of the previous point, with regularity condition on abstract merge trees being replaced by being without degree $2$ vertices for merge trees.
\end{enumerate}
\end{proof}

\bigskip\noindent
\textbf{\Cref{prop:metric}.}
The display poset $D_{\T}$ of any abstract merge tree can be given a pseudo-metric structure with the following formula:
\[
d((a,t),(b,t'))= (\tilde{t} - t)  + (\tilde{t} - t')  
\]
with $\tilde{t}= \inf \{h(p)\mid p \in \CA(\{(a,t),(b,t')\})\}$. If $\T$ is regular, then, $d$ is a metric.

\begin{proof}
    First note that even if 
$\inf \CA(Q)$, with $Q\subset D_{\T}$ and $\sup h(Q)<\infty$,  may be a set with more than one element, $\inf\{h(p)\text{ with }p\in D_{\T}\mid p\geq Q\}$ is uniquely defined.
Moreover, consider $p=(b,t_b),q=(c,t_c)\in \inf \CA(Q)$. For every $(a,t)\in \CA(Q)$ we know $\pi_0(X_{t_b\leq t})(b)=\pi_0(X_{t_c\leq t})(c)=a$. Clearly $t_b$ and $t_c$ must be critical values otherwise we can consider $p'>p$ and $q'>q$ with $q',p'\leq Q$, which is absurd. But the same holds if 
$t_b\neq t_c$: suppose $t_b< t_c\leq h(Q)$ then $p' = (\pi_0(X_{t_b<t_b+\varepsilon})(b),t_b+\varepsilon)$
(with $\varepsilon>0$ small enough) satisfies $p<p'$ and $p'\leq Q$, which is absurd. Thus $t_b= t_c=t_i$ critical value.

The map $d:D_{\T} \times D_{\T}\rightarrow \mathbb{R}_{\geq 0}$ is symmetric. 
For what have said before $d(p,q)=0$ if and only if $p,q\in \inf \CA(\{p,q\})$ and 
$h(p)=h(q)=t_i$ critical value. This is equivalent to
 $p=(b,t_i),q=(c,t_i)\in D_{\T}$ such that $\pi_0(X_{t_i<t_i+\varepsilon})(b)=\pi_0(X_{t_i<t_i+\varepsilon})(c)$ for every $\varepsilon>0$.

Thus, if $\T$ is regular we have  
$d((b,t_i),(c,t_i))=0$ if and only if 
$p=q$; in fact $X_{t_i<t_i+\varepsilon}$ is an isomorphism for $\varepsilon>0$ small enough.

Now we check the triangle inequality.
Let $p_1,p_2,p_3\in D_{\T}$. 
And let $t_i=h(p_i)$, $Q_{ij}=(p_i,p_j)$, $q_{ij}=\inf\{h(p)\text{ with }p\in \CA( Q_{ij})\}$ and  $q=\inf\{h(p)\text{ with }p\in \CA(\{p_1,p_2,p_3\})\}$.

Consider $P_1=\CA(\{p_1\})$. 
Clearly $\inf \CA(\{p_1,p_2\})\subset P_1$ and 
$\inf \CA(\{p_1,p_3\})\subset P_1$.
Thus either (1) $q_{13}\leq q_{12}$ (and $q_{23}= q_{12}$) or (2)
$q_{12}< q_{13}$ (and $q_{13}= q_{23}$) hold.

In case (1) holds: 
\[
q_{12}-t_1=q_{12}-q_{13}+q_{13}-t_1 \leq
q_{12}-q_{13}+q_{13}-t_1 + 2q_{13}-2t_3 = 
q_{13}-t_1 + q_{13}-t_3 + q_{23}-t_3  
\]
Thus:
\[
q_{12}-t_1+q_{12}-t_2\leq   
q_{13}-t_1 + q_{13}-t_3 + q_{23}-t_3 + q_{23}-t_2 
\]

The proof in case (2) holds is analogous. 
\end{proof}

\bigskip\noindent
\textbf{\Cref{prop:lattice}.}
The set $\Cov(\T)$ is a lattice. It is a poset
with the relation $\mathcal{O}<\mathcal{O}'$ if $\mathcal{O}$ is a refinement of $\mathcal{O}'$ and for every pair of elements $\mathcal{O}$, $\mathcal{O}'$ there is a unique least upper bound $\mathcal{O}\vee \mathcal{O}'$ and a unique greater lower bound $\mathcal{O}\wedge \mathcal{O}'$. The operations are defined as follows:
\[
\mathcal{O}\vee \mathcal{O}':= \pi_0\left( \bigcup_{U\in\mathcal{O}' \text{ or }U\in\mathcal{O}} U\right)
\]

\[
\mathcal{O}\wedge \mathcal{O}':= \{U\cap U' \mid U'\in\mathcal{O}' \text{ and }U\in\mathcal{O}\}. 
\]
\begin{proof}    
    Let's start with $\mathcal{O}\vee \mathcal{O}'$. It is clearly an a.e. covering. Moreover $ \bigcup_{U\in\mathcal{O}' \text{ or }U\in\mathcal{O}} U$ is clearly contained in $ \bigcup_{U\in\mathcal{U}(D_{\T})} U$, and by functoriality we have that the set $\pi_0\left( \bigcup_{U\in\mathcal{O}' \text{ or }U\in\mathcal{O}} U\right)$ is included in $\pi_0\left( \bigcup_{U\in\mathcal{U}(D_{\T})} U\right).$ And the latter is equal to $\mathcal{U}(D_{\T})$. Thus $\mathcal{O}\vee \mathcal{O}'$ is regular and clearly is refined by $\mathcal{O}$ and $\mathcal{O}'$. Lastly, consider any $\mathcal{O},\mathcal{O}'<\mathcal{O}''$.  Since the sets of $\mathcal{O}''$ are disjoint and path-connected (by construction), then any $U''\in \mathcal{O}''$ contains all the sets of $\mathcal{O}$ and $\mathcal{O}'$ it intersects. Thus it contains a path-connected component of their union.

    Now we turn to $\mathcal{O}\wedge \mathcal{O}'$. All the sets in $\mathcal{O}\wedge \mathcal{O}'$ are disjiont, open and path-connected. And they form an a.e. cover of $D_{\T}$ - otherwise a positive-measure set would be left out by $\mathcal{O}$ or $\mathcal{O}'$. Thus  $\mathcal{O}\wedge \mathcal{O}'$ is a regular a.e. covering which refines $\mathcal{O}$ and $\mathcal{O}'$.
    Consider $\mathcal{O}''$ such that $\mathcal{O},\mathcal{O}'>\mathcal{O}''$. Take $U''\in \mathcal{O}''$. By construction there are $U\in \mathcal{O}$ and $U'\in \mathcal{O}'$ with $U''\subset U',U$. Thus $U''\subset U\cap U'$. So $\mathcal{O}''<\mathcal{O}\wedge \mathcal{O}'$.
\end{proof}

\bigskip\noindent
\textbf{\Cref{prop:truncation}.}
Take $(T,\varphi_T)$ and  $(T',\varphi_{T'})$.
Suppose $r_T$ and $r_{T'}$ are of degree $1$ and there is a splitting $\{(v,r_T)\}\rightarrow \{(v,v'),(v',r_T)\} $ and $\{(w,r_{T'})\}\rightarrow \{(w,w'),(w',r_{T'})\} $ giving the dendrograms $(G,\varphi_{G})$ and $(G',\varphi_{G'})$. Suppose moreover that $\varphi_{G}((v',r_{T}))=\varphi_{G'}((w',r_{G}))$. Then $d_E(T,T')=d_E(sub_{G}(v'),sub_{G'}(w'))$.

\begin{proof}

Consider a minimizing mapping $M \in M_2(G,G')$ (see \Cref{sec:mappings}). 

Apply the deletions and ghostings described by $M$ both on $G$ and on $G'$ obtaining, respectively, the merge trees $G_M$ and $G'_M$. Note that, since $M \in M_2(G,G')$, 
if $v'$ is paired with another edge by $M$, then $v$ must be deleted, otherwise $v'$ would be of degree $2$ after the deletions, and thus it should be ghosted. Similarly, if $v'$ is deleted then $v$ must be deleted as well as: otherwise 
$v'$ would be of degree $2$ after all the other deletions, and thus it cannot be deleted by $M$.
The same of course applies to $w'$ and $w$.

Let $\varphi_{T'}(v')= f'$, $\varphi_{T'}(v)= f$, $\varphi_{G'}(w')= g'$  and $\varphi_{G'}(w)=g$. By construction $f'=g'$ and thus the cost of the pair $(v',w')$, if it can be added a mapping, is zero. Since $M \in M_2(G,G')$, we know $(v',w')\notin M$. Our goal is to build a mapping containing $(v',w')$ with cost equal to the cost of $M$.

\begin{itemize}
    \item First suppose that $v$ is not deleted and $w$ is not deleted.
    
Then there exist $a = v_1 \leq \ldots \leq v_k \leq v'$ and $b = w_1 \leq \ldots \leq w_r \leq w'$ vertices in $T$ and $T'$ respectively such that:
\begin{itemize}
    \item $(a,b)\in M$;
    \item $v_i$ and $w_i$ are ghosted for all $i$ (as $M \in M_2(G,G')$ they cannot be deleted).
\end{itemize}

Let $f_i = \varphi_T(v_i)$ and $g_i = \varphi_{T'}(w_i)$ for all $i$.

We have seen that if $a=v'$, then $v$ must be deleted, which is absurd. Similarly we cannot have $b=w'$. 
Thus: $a= v_1 \leq \ldots \leq v_k = v <v'$
and $b= w_1 \leq \ldots \leq w_r = w<w'$. 
Using $f'=g'$ and property (P4), we have:
\[
d( f' + \sum_{i=1}^{k} f_i , g' + \sum_{i=1}^{r} g_i) = d(\sum_{i=1}^{k} f_i,\sum_{i=1}^r g_i).  
\]

This implies that we can add the pair $(v',w')$ to $M$ without increasing the cost of $M$.

\item Suppose now $v$ (and so $v'$) is deleted but $w$ is not. For what we have said before, it means that $w'$ is ghosted and $w$ is either paired or deleted. 
This implies that the root of $G'_M$ is of degree $1$. But then also the root of $G_M$ must be of degree $1$. In other words, there exist $a = v_1 \leq \ldots \leq v_k < v$ and $b = w_1 \leq \ldots \leq w_r \leq w$ vertices in $T$ and $T'$ respectively such that:
\begin{itemize}
    \item $(a,b)\in M$;
    \item $v_i$ and $w_i$ are ghosted for all $i$ (as $M \in M_2(G,G')$ they cannot be deleted).
\end{itemize}

But this also implies that, after all the other deletions $v$ is of degree $2$, and the same for $v'$. But this cannot happen.
\item Suppose lastly, $v$ and $w$ are all deleted. Which implies that also $v'$ and $w'$ are deleted. In this case we can add $(v,w)$ and $(v',w')$ to $M$ decreasing its cost.
\end{itemize}
Thus we can always add $(v',w')$ to $M$ and since the cost of such pair is zero, we have $d_E(T,T')=d_E(G,G')=d_E(sub_{G}(v'),sub_{G'}(w'))$.
\end{proof}

\bibliography{references}

\end{document}